%% file: gridley.tex
\documentclass[10pt]{amsart}
\usepackage{amsmath, amssymb,graphicx,amscd,color,dbnsymb}
\theoremstyle{plain}

\input{macros}

\begin{document}

\title{Floer Homology and Knot Complements}
\author{Jacob Rasmussen}
\address{Harvard University Dept. of Mathematics, Cambridge, MA 02138}
\email{jacobr@math.harvard.edu}
\thanks{The author was partially supported by a Department of Defense graduate fellowship.}

\begin{abstract}
We use the \OS theory of Floer homology to define an invariant of
knot complements in three-manifolds. This invariant takes the form of
a filtered chain complex, which we call \( \cfr \). It carries
information about the \OS Floer homology of large integral surgeries on the
knot. Using the exact triangle, we derive information about other
surgeries on knots, and about the maps on Floer homology induced by
certain surgery cobordisms. We define a certain class
of {\it perfect}  knots in \(S^3\) for which \( \cfr \) has a
particularly simple
form. For these knots, formal properties of the \OS  theory enable
us to make a complete calculation of the Floer homology. It turns out
that most small knots are perfect. 
\end{abstract}
%\doublespacing
%\maketitle
\include{introduction}

\include{osreview}
\include{heegaard}

\include{filtration}

\include{theorem1}
\include{quasiperfect}

\include{triangle}

\include{alternating}

\include{appendix}

\bibliographystyle{plain}
\bibliography{mybib}

\end{document}

%% file: macros.tex
\theoremstyle{plain}
\newtheorem{lem}{Lemma}[section]
\newtheorem{prop}{Proposition}[section]
\newtheorem{thrm}{Theorem}[section]
\newtheorem{cor}{Corollary}[section]
\newtheorem{principle}{Principle}[section]
\newtheorem{defn}{Definition}[section]
\newtheorem{conv}{Convention}[section]
\newtheorem{tm}{\noindent Theorem}

\newcommand{\Ozsvath}{{Ozsv{\'a}th  }}
\newcommand{\Szabo}{{Szab{\'o}  } }

\newcommand{\MM}{{\mathcal{M}}}

\newcommand{\R}{{\bf{R}}}

\newcommand{\C}{{\bf{C}}}
\newcommand{\Z}{{\bf{Z}}}
\newcommand{\Q}{{\bf{Q}}}
\newcommand{\SSS}{{\Sigma}}

\newcommand{\into}{{\hookrightarrow}}

\newcommand{\dlim}{{\varinjlim \ts}}
\newcommand{\ilim}{{\varprojlim \ts}}
\newcommand{\coker}{{{\text { \rm coker} \ }}}
\newcommand{\im}{{{\text {im} \ }}}
\newcommand{\rank}{{{\text {rank} \ }}}

\newcommand{\sign}{{{\text {sign} }}}

\newcommand{\hfp}{{{HF^+}}}
\newcommand{\hfm}{{{HF^-}}}
\newcommand{\hfred}{{{HF^{red}}}}
\newcommand{\hfpm}{{{HF^{\pm } }}}
\newcommand{\hfi}{{{HF^{\infty} }}}
\newcommand{\hfhat}{{{\widehat{HF}}}}
\newcommand{\cfp}{{{CF^+}}}
\newcommand{\cfm}{{{CF^-}}}
\newcommand{\cfpm}{{{CF^{\pm } }}}
\newcommand{\cfhat}{{{\widehat{CF}}}}
\newcommand{\cfi}{{{CF^{\infty} }}}
\newcommand{\cfis}{{{CF_s^{\infty} }}}
\newcommand{\cfia}{{{CF_a^{\infty} }}}
\newcommand{\cfs}{{{\widehat{CF} _s}}}
\newcommand{\hfs}{{{\widehat{HF} _s}}}
\newcommand{\hfps}{{{HF_s^+}}}
\newcommand{\cfps}{{{CF_s^+}}}
\newcommand{\hfms}{{{HF_s^-}}}
\newcommand{\cfms}{{{CF_s^-}}}
\newcommand{\cfpms}{{{CF_s^\pm}}}
\newcommand{\cfr}{{{\widehat{CF} _r}}}
\newcommand{\cfpr}{{{CF_r^+}}}
\newcommand{\cfmr}{{{CF_r^-}}}

\newcommand{\cfir}{{{CF_r^\infty}}}
\newcommand{\cfa}{{{\widehat{CF} _a}}}

\newcommand{\cfpa}{{{CF_a^+}}}
\newcommand{\hfpa}{{{HF_a^+}}}
\newcommand{\cfma}{{{CF_a^-}}}

\newcommand{\hfmt}{{{\widetilde{HF}^-}}}
\newcommand{\cfpma}{{{CF_a^\pm}}}

\newcommand{\Zu}{{{\Z [u^{-1}]}}}
\newcommand{\Qu}{{{\Q [u^{-1}]}}}
\newcommand{\Zt}{{{\Z [T,T^{-1}]}}}
\newcommand{\Qt}{{{\Q [T,T^{-1}]}}}
\newcommand{\spi}{{{\mathbf s}}}
\newcommand{\vspi}{{{\mathfrak t}}}
\newcommand{\spinc}{{{\mathrm{Spin}^c}}}
\newcommand{\gr}{{{\mathrm{gr} }}}

\newcommand{\DDD}{{{\mathcal D }}}
\newcommand{\PPP}{{{\mathcal P }}}
\newcommand{\MMM}{{\widehat{\mathcal M }}}
\newcommand{\AAA}{{{\mathcal A }}}
\newcommand{\CCC}{{{\mathcal C }}}
\newcommand{\TTT}{{{\bf T }}}
\newcommand{\Tbar}{{{\overline{\TTT}}}}
\newcommand{\jjj}{{{\mathfrak j }}}

\newcommand{\OS}{{{Ozsv\'ath-Szab\'o \ts }}}
\newcommand{\scx}{{{\mathcal S }}}
\newcommand{\fmi}{{{f_-^\infty}}}

\newcommand{\bfy}{{{{\bf y}}}}
\newcommand{\bfx}{{{{\bf x}}}}
\newcommand{\bfw}{{{{\bf w}}}}
\newcommand{\ybar}{{{\overline{\bf y}}}}
\newcommand{\bfz}{{{{\bf z}}}}
\newcommand{\zbar}{{{\overline{\bf z}}}}
\newcommand{\haty}{{{\hat{{\bf y}}}}}
\newcommand{\hatz}{{{\hat{{\bf z}}}}}
\newcommand{\TAB}{{{ \TTT_{\hat{\alpha}} \cap \TTT_{\hat{\beta}}}}}
\newcommand{\TABO}{{{ \TTT_{\hat{\alpha_1}} \cap \TTT_{\hat{\beta}}}}}
\newcommand{\TABE}{{{ \TTT_{{\alpha}} \cap \TTT_{{\beta}}}}}

\newcommand{\balpha}{{{\mathbf \alpha}}}
\newcommand{\bbeta}{{{\mathbf \beta}}}
\newcommand{\Ko}{{{K^o}}}
\newcommand{\bbar}{{{\overline{\beta}}}}
\newcommand{\bbbar}{{{\overline{\bbeta}}}}
\newcommand{\ebar}{{{\overline{\epsilon}}}}
\newcommand{\Nbar}{{{\overline{N}}}}
\newcommand{\tmod}{{{\mathcal{T}}}}
\newcommand{\Ebar}{{{\overline{E}}}}
\newcommand{\spibar}{{{\overline{\spi}}}}
\newcommand{\mutid}{{{\tilde{\mu}}}}

\newcommand{\hatf}{{{\hat{f}}}}
\newcommand{\eqy}{{{\mathcal{Y}}}}
\newcommand{\E}{{{\mathcal H}}}

\newcommand{\ts}{{{\thinspace}}}
\newcommand{\emm}{{{\mathfrak m}}}
\newcommand{\Nsum}{{{N_1 \#_{m} N_2}}}
\renewcommand{\:}{{{\colon}}}

%% file: introduction.tex
\maketitle

\section{Introduction}
\label{Sec:Introduction}

This thesis studies the \OS Floer homology groups \( \hfpm \)
for  three-manifolds obtained by surgery on a knot. 
 These groups were introduced by \Ozsvath and \Szabo in \cite{OS1} and
  are conjectured to be isomorphic to the
Seiberg-Witten Floer homology groups \( HF^{to} \) and \( HF^{from} \)
described by Kronheimer and Mrowka in \cite{MKLect}. They 
 possess many of the known or expected properties of the
Seiberg-Witten Floer homologies, and have been used to define four-manifold 
invariants analogous to the
Seiberg-Witten invariants (\cite{OS3}, \cite{OS4}). 
It would be a mistake, however, to think of the \OS theory as
being nothing more than a convenient technical alternative to
Seiberg-Witten Floer homology. The \OS invariants offer  a new, and in many
respects quite different, perspective on Seiberg-Witten theory. 
As we hope to illustrate, they are substantially more computable than
the corresponding Seiberg-Witten objects. This computability
is a consequence of formal properties which have no obvious
gauge theoretic counterparts. In what follows, we will describe these
formal properties and explain how they can be used to calculate certain 
\OS Floer homologies and four-manifold invariants.

The basic question  we hope  to address is ``What is the Floer homology
of a knot?'' In other words, we seek a single object which encodes
information about the Floer homologies of closed manifolds obtained
from the knot complement --- either by Dehn surgery, or, more
generally, by gluing two knot complements together. That such an
object should exist was suggested by known results (see {\it e.g.}
\cite{ExactT}, \cite{KKR}, \cite{MOY}) in both monopole and instanton
Floer homology. Although the \OS theory has not provided a
complete solution to our question, it does offer a partial solution
very different from   anything envisioned by Seiberg-Witten theory. 
This ``knot Floer homology'' takes the form of a filtered chain complex,
whose homology is \( \cfhat \) of the ambient three-manifold \(Y\). 
That this might be the case
 was suggested by \cite{2bridge}, which used the formal properties
of the \OS Floer homology to compute \( \hfp \) for surgeries on a
particularly simple class of knots --- the two-bridge knots. For these
knots, the complex \( \cfhat \) has a very simple and symmetrical
form. This thesis represents an attempt to find a sense
in which this form generalizes to other,  more complicated knots. 

The invariants and many of the results described here have been
independently discovered by \Ozsvath and \Szabo in \cite{OS7},
\cite{OS8}, and \cite{OS10}. 
In some places, most notably for alternating knots, their
results are better than those given here. 

 We now outline the ideas and main results of the thesis. 
Throughout, a knot \(K\) will be a null-homologous curve in a
three-manifold \(Y\). Often, it is more convenient to think of the
knot in terms of its complement \(N= Y-K\), which is a three-manifold
with torus boundary, together with a special curve (the meridian) on
that boundary. We write \( K(p,q) \) and \(N(p,q) \) interchangably to
denote the result of \(p/q \) Dehn surgery on \(K\). 
\subsection{The reduced complex of a knot}
One of the insights which  \OS  theory  provides us is that it is
much easier and more natural to think about {\it large} surgeries.
 Given a knot \(K\), we study \( \hfp (K(n,1)) \) for \( n \gg 0
\). As described in sections~\ref{Sec:Heegaard} and
\ref{Sec:Filtration}, we are  very naturally led to a certain complex
\( \cfs (K) \). This complex comes with a filtration, which we call
the {\it Alexander filtration}. (The Euler characteristics of its
filtered subquotients are the coefficients of the Alexander polynomial
of \(K\).) Using
this filtration, we introduce a refined version of \( \cfs (K) \),
in which we replace each filtered subquotient by its homology. We denote
the resulting complex by \( \cfr (K) \). Our main result is 
\begin{tm}
\label{Thm:ReducedInvariance}
The filtered complex \( \cfr (K) \) is an invariant of \(K\). 
\end{tm}

By itself, the homology of \( \cfr (K) \) is not very interesting: it is always
isomorphic to \( \cfhat (Y) \). If we understand the filtration on \(
\cfr (K)\), however, we can use it to compute \( \hfhat (K(n,1),
\spi_k) \)  (\(n\gg 0 \))
for any \( \spinc \) structure \( \spi _k \). As described
below, the 
stable complex  can  often  be used to compute  \( \hfp (K(n,1), \spi_k) \)
 for all integer values of \(n\). 
% In addition, it behaves nicely under connected sum: 
%\( \cfr(K_1 \# K_2 ) \cong \cfr (K_1) \otimes \cfr (K_2) \). 
Thus, we like to think of the stable complex as being at least a partial answer to the question ``What is the Floer homology of a knot complement?''
\subsection{Applications of the exact triangle}

By itself, the stable complex is only useful for understanding the
results of large-\(n\) surgery on \(K\). In applications, however, we
are usually interested in \(K(n,1)\) for small values of \(n\). For
the stable complex to be useful to us, we need some method of relating
large-\(n\) surgeries to other  surgeries on \(K\).
This method is provided by the exact triangles of \cite{OS2}. 

For the remainder of this introduction, we restrict our attention to
the case where \(K\) is a knot in \(S^3\).  In this situation, we
define invariants \(h_k(K) \) for \(k \in \Z\). Essentially, \(h_k (K) \)
is the rank of the map \( \hfm (S^3) \to \hfred (K(0,1)) \) induced by
the surgery cobordism. If we know  the \(h_k (K) \)'s and the groups
\( \hfp (K(n,1), \spi _k) \) for \(n \gg 0 \), we can use the exact
triangle to compute the  \OS Floer homology of any integer surgery on \(K\). 

{\it A priori}, the \(h_k(K) \)'s are defined using a cobordism. We
show, however, that they can be computed purely in terms of the
complex \( \cfps (K) \). This provides us with an effective means
of finding the \(h_k\)'s from the stable complex of \(K\). It also
enables us to prove some general theorems  about their behavior.
Using the fact that \(\hfp (S^3) \cong \Z\), we define a knot invariant \(s(K) \) and show that \(h_k (K) 
> 0 \) when  \(|k| < s(K) \).
\subsection{Perfect knots}
We say that \(K\) is {\it perfect} if the Alexander filtration on \( \cfr (K) \) is the same as the filtration induced by the homological grading. 
 In this case, we have the following theorem, which generalizes the
calculation of \( \hfp \) for two-bridge knots given in  \cite{2bridge}:
\begin{tm}
\label{Thm:QPCalc}
If \(K\) is perfect, the groups \(
\hfp (K_n, \spi_k) \) are determined by the Alexander polynomial of \(K\) and the invariant \(s(K) \). 
\end{tm}

 The exact form of \(  \hfp (K_n, \spi_k) \) is described in
 Theorem~\ref{Thm:PerfectKnots}. For the moment, we remark that the
 method employed is quite different from the calculations of the
 instanton and monopole Floer homology for Seifert fibred-spaces in
 \cite{FSSFHF} and \cite{MOY}, which rely on having a chain complex
 in which all generators have the same \( \Z/2\) grading (so \(d\) is
 necessarily trivial.) The complex with which we compute is typically
 much larger than its homology. 

If \(K\) is perfect, the invariant \(s(K) \) is algorithmically
computable.  For every perfect knot for which we have carried out this
computation, \(s(K) = \sigma (K) /2\), where \( \sigma (K) \) is the
classical knot signature. For nonperfect knots, however, the two are
generally different.
\subsection{Alternating knots}
Although perfection is a strong condition, it is
 satisfied by a surprisingly large number of knots. In particular, we
have 
\begin{tm}
\label{Thm:Alternating=QP}
Any small alternating knot is perfect.
%where \(\sigma \) is the ordinary knot signature. 
\end{tm}
 
Here {\it smallness } is a technical condition which we will  describe in
section~\ref{Sec:Alternating}. It is satisfied by all but one
alternating knot with crossing number \( \leq 10\). 
Combining this result with some hand calculations for nonalternating
knots, we  can show that there are only two nonperfect knots with
\(9\) or fewer crossings. These knots are numbers \(8_{19}\) (the
\((3,4) \) torus knot) and \(9_{42}\) in Rolfsen's tables.  

 In fact, the  hypothesis of smallness is unnecessary. In
 \cite{OS8}, \Ozsvath and \Szabo have shown that {\it all} alternating
 knots are perfect and have \( s(K) = \sigma (K) /2 \). 
\subsection{The Stable Complex as a Categorification}
If we wish, we can put aside \( \cfr\)'s connection with gauge theory,
and simply view it as an invariant of knots in \(S^3\). From this
point of view, the stable complex is best thought of as a
generalization of the Alexander polynomial: it is a filtered complex
with homology \(\Z\) whose filtered Euler characteristic is \( \Delta
_K(t) \). Many properties of the Alexander polynomial
carry over to the stable complex. For example, the Alexander
polynomial is symmetric under inversion: \( \Delta _K(t) = \Delta
_K(t^{-1}) \). The reduced stable complex also has a such a symmetry:
it is the analog of the  conjugation symmetry in Seiberg-Witten
theory. The Alexander polynomial is 
multiplicative under connected sum: \( \Delta (K_1 \#K_2) = \Delta
(K_1)\Delta (K_2) \). Likewise, \( \cfr(K_1 \# K_2 ) \cong \cfr (K_1) 
\otimes \cfr (K_2) \). The Alexander polynomial is defined by a skein
relation; the stable complex satisfies (but does not seem to be
determined by) an analogous {\it skein exact triangle.} Finally, the
degree of \( \Delta _K(t) \) gives a lower bound for the genus of
\(K\). The same is true for the highest degree in which \( \cfr (K) \)
is nontrivial; this is the adjunction inequality. In fact, if one
believes the conjecture relating \OS and Seiberg-Witten Floer
homologies, work of  Kronheimer and Mrowka \cite{Kr} implies that 
 this degree should be precisely {\it equal } to the genus. 

It is interesting to compare these properties of the stable complex
with Khovanov's categorification of the Jones polynomial
\cite{Khovanov}. Khovanov's invariant is a filtered sequence of
homology groups whose filtered Euler characteristic gives the Jones
polynomial of \(K\); it has recently been shown by Lee \cite{ESL2}
that these groups can be given the structure of a complex with
homology \(\Z\oplus \Z\). The similarity becomes even more striking
when one considers the  reduced Khovanov homology introduced by
Khovanov in \cite{Khovanov2}. In many instances, the rank of this group is
isomorphic to that of \( \cfr (K) \). There is also a quantity which
resembles the invariant \(s(K) \). Unfortunately, the two groups are
not always the same: one example is given by the 
 \((4,5)\) torus knot, for which \( \cfr \) has rank \(7\), but the
reduced Khovanov homology has rank \(9\). It is an interesting problem
to find some explanation for why these two groups should often, but
not always, be similar. 

\subsection{Acknowledgements}
The author would like to thank Peter Kronheimer for his
constant advice, support, and encouragement
 throughout the past five years. He would also like to thank Peter
\Ozsvath and Zoltan \Szabo for their interest in this work and their
willingness to share their own, and Tom Graber, Nathan Dunfield, and Kim
Fr{\o}yshov for many helpful discussions.

%% file: osreview.tex
\section{The \OS Floer homology}

%At first sight, the computational advantages offered by the \OS 
%Floer homology are not entirely obvious. Although the generators of
%their chain  complex can be determined  mechanically, there are
%often a lot of them --- many more than is usual for Seiberg-Witten
%theory. Moreover, it is not so clear how one could ever hope to work
%out the differentials between all these generators. (My own first
%experience with the subject was an attempt to compute \( \hfhat \) for
%the Poincare sphere. I ended up with a complex with \(27\) generators,
%and no idea how to compute the differentials between them!) Despite
%these apparent obstacles, the \OS  theory is actually very easy to
%compute with. We hope to convey this fact in what follows. 

In this section, we give  a quick review of some basics from
\cite{OS1} and \cite{OS2}. Of course, this is not a substitute
 for these papers, and we encourage the reader to look at them
as well. (Especially the very instructive examples in section 8 of
\cite{OS1}.) We  focus on concrete, two dimensional interpretations of
the objects involved. Let \(Y\) be a closed three-manifold, which we
assume (at least for the moment) to be a rational homology sphere. 
\subsection{Heegaard splittings}

The \OS Floer
 homology of a closed three-manifold \(Y\) is defined using a Heegaard
splitting for \(Y\), {\it i.e.} a  choice of a surface \(
\Sigma \subset Y\) which divides \(Y\) into two  handlebodies. We
can describe such a splitting by starting with a surface \(\Sigma \)
of genus \(g\) and drawing two systems of {\it attaching handles} on
it. Each system  is a set of \(g \) disjoint,
smoothly embedded circles, such that if we surger \( \Sigma \) along
them, the resulting surface is connected (and thus homeomorphic to
\(S^2\).) To recover  \(Y\), we thicken \( \Sigma \)
and glue in a two-handle \(D^2 \times I \) along each attaching
circle. We then fill in the two remaining \(S^2\) boundary components
 with copies of \(B^3\).

\begin{figure}
\includegraphics{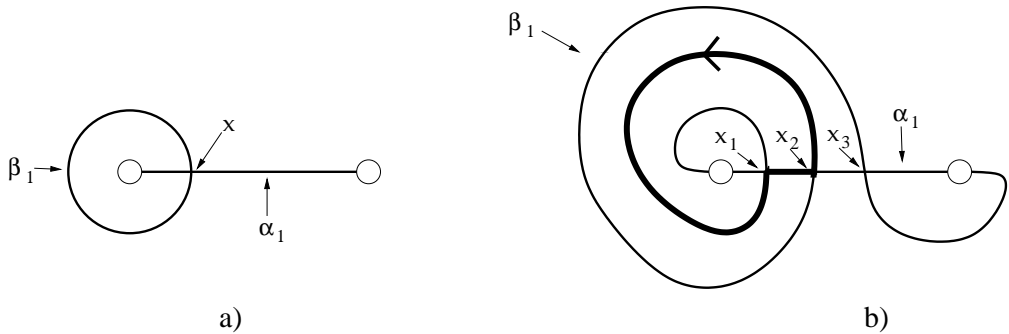}
\caption{\label{Fig:GenusOne} Heegaard splittings of {\it a)} \(S^3\)
  and {\it b)} the lens space \(L(3,1)\). }
\end{figure}

Some simple genus \(1\) Heegaard splittings are shown in
Figure~\ref{Fig:GenusOne}. These (and all other pictures of Heegaard
splittings in this paper) are drawn using the following method. Think of
\(S^2\) as the plane of the paper. To represent the surface \( \Sigma \),
we draw \(g\) pairs of small disks in the plane and surger each pair. A curve
on \( \Sigma \) which goes ``into'' one disk in a pair
 comes ``out'' of the other. (This is the same convention used to
 represent one-handles in Kirby calculus, but one dimension down.)
We label the curves in the first system of attaching handles \( \alpha
_1, \ldots, \alpha _g \), and the curves in the second system \( \beta
_1, \beta _2, \ldots, \beta _g \). Usually, we arrange things so that each \(
\alpha _i \)  is a straight line which joins the two small disks in a
pair. 

\begin{figure}
\includegraphics{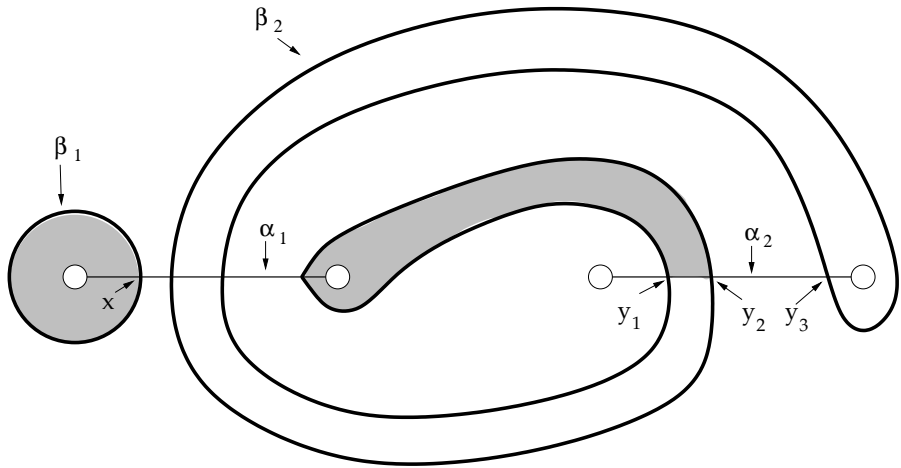}
\caption{\label{Fig:GenusTwo}Another splitting of \(S^3\). The shaded
  region is the domain of a differential.}
\end{figure}

Any three manifold is represented by many different Heegaard splittings. 
For example, the slightly more complicated splitting shown in  
 Figure~\ref{Fig:GenusTwo} also represents \(S^3\).

\subsection{Generators}
Let  \((\Sigma, \alpha, \beta) \) be a Heegaard splitting of \(Y\).
 The system of attaching circles \( \alpha \)  defines a 
Lagrangian torus
\begin{equation*}
 \TTT_\alpha =  \{\alpha _1 \times \alpha
  _2 \times \ldots \times  \alpha _g\} \subset s^g\Sigma 
\end{equation*}
 where \( s^g \Sigma \) is the \(g\)th symmetric product of \( \Sigma\). 
 The simplest \OS homology,
\( \hfhat (Y) \), is a Lagrangian Floer homology defined
by the pair \((\TTT_\alpha, \TTT_\beta)\) inside of the  symplectic
 manifold \(s^g
 \Sigma\). The generators of the complex \( \cfhat (Y)
 \) are the points in \( \TTT_\alpha \cap \TTT_\beta \). The
 differential between two generators \(\bfy \) and \( \bfz\) is 
  defined by counting ``holomorphic disks'' joining \(\bfy \) to
 \(\bfz\); that is, psuedoholomorphic maps \(\phi \: D^2 \to s^g\Sigma\)
 such that \(\phi (-1) = \bfy\), \( \phi (1) = \bfz\), \( \phi
 (C_-)\subset \TTT_\alpha \), and \( \phi (C_+) \subset
 \TTT_\beta\). (\(C_\pm \) are the components of \( \partial D^2 \)
 lying above and below the real axis.) 
We denote the set of homotopy
 classes of such disks by \(\pi_2(\bfy, \bfz) \). For \(g > 1 \), this
 set is either empty or an affine space modeled on \( \pi_2
 (s^g\Sigma ) \cong \Z\). In order to define \( \cfhat \), 
 we must specify which homotopy class we want to count holomorphic
 disks in, as described in  section~\ref{SubSec:Diffs} below. 

In \cite{OS1}, \Ozsvath and \Szabo showed that the homology group \(
\hfhat (Y) \) has two remarkable properties. First, it does not depend
on the choice of Heegaard splitting used for \(Y\). Second, all
of the objects involved can be described in terms of the Heegaard
diagram \( (\Sigma, \alpha, \beta) \). 
For example, the points in 
\( \TTT_\alpha \cap \TTT_{\beta} \)
correspond to unordered \(g\)-tuples of intersection points \(\{x_1,
x_2, \ldots x_g\} \) between the \( \alpha _i \)'s and the
\(\beta_j\)'s, such that every \(\alpha _i \) and \( \beta _j \)
contains exactly one \(x_k\).
Thus if we use the diagram of Figure~\ref{Fig:GenusOne}a, we see
that  \( \cfhat (S^3) \) is generated by the single
intersection point \(x\), so \( \hfhat (S^3)
\cong \Z \). If we use the diagram of Figure~\ref{Fig:GenusTwo},
however, \( \cfhat \) is generated by the three pairs \( \{x,y_1\}, 
 \{x,y_2\}\), and \(  \{x,y_3\}\). Since  \(
 \hfhat (S^3) \cong \Z \), there must be a nontrivial differential in
 this complex.

\subsection{\(\epsilon\)-grading and basepoints}

Given generators \( \bfy, \bfz \in \TTT_\alpha \cap \TTT_\beta \),
 there is a topological obstruction to the existance of a disk in \(
 \pi_2 (\bfy, \bfz) \). As an example,  consider the points \(x_1 \) and
 \(x_2\) of Figure~\ref{Fig:GenusOne}b. The heavy line traces out a
 curve which goes from \(x_1 \) to \(x_2 \) along \( \TTT_\alpha =
 \alpha_1 \) and then from \(x_2 \) to \(x_1 \) along \(\TTT_\beta =
 \beta _1 \). This curve  represents a nontrivial class \(m\)
 in the homology of \(s^1\Sigma = T^2\), so it cannot bound a disk. We
 could try to rectify this problem by choosing different paths from
 \(x_1\) to \(x_2 \) along \( \alpha _1 \) and from \(x_2 \) to
 \(x_1\) along \( \beta _1\). This has the effect of replacing the
 homology class \(m\) with \(m + a[\alpha _1] + b[\beta_1] \) for some
 \(a, b \in \Z \). It is easy to see that  the resulting homology class
 is always nonzero, so \( \pi_2(x_1, x_2) \) is empty. 

In general, the obstruction can be  expressed as a map
\begin{equation*}
 \epsilon \: \TTT_\alpha \cap \TTT_\beta \to \text{Affine} (H_1(Y))
\end{equation*}
which we call the \( \epsilon \)-grading. To define this grading, we
use the notion of a {\it system of paths} joining \(\bfy = \{y_1,
\ldots, y_g\}\)  to \( \bfz = \{z_1, \ldots, z_g\} \)
 on \( \Sigma \). Such a system is a set of polygons with a total of
\(2g \) edges, which map to \(\Sigma \) in the following way: vertices
go alternately to \(y_i\)'s and \(z_j\)'s, and edges go alternately to
 \( \alpha_i\)'s and \(\beta_j\)'s. Every \(y_i\) and \(z_j\) is the
 image of exactly one vertex and every \(\alpha_i\) and \(\beta
 _j\) contains the image of exactly one edge. If \(C(\bfy, \bfz) \) is
 such a system, then \( \epsilon (\bfy) - \epsilon(\bfz) \) is the
 image of \( C(\bfy, \bfz) \) in 
\begin{equation*}
\frac{H_1(s^g\Sigma)}{H_1(\TTT_\alpha)
\oplus H_1(\TTT _\beta)}   \cong \frac{H_1(\Sigma)}{\langle \alpha _i,
  \beta_j\rangle} \cong H_1(Y).
\end{equation*}

The \( \epsilon \)-grading is always easy to compute. Indeed, as we
 describe in the next section, it is essentially just a
 geometric realization of Fox calculus. It divides the points in \(
 \TTT_\alpha \cap \TTT_\beta \) into equivalence classes, which we
 call \( \epsilon\)-classes. Two points \( \bfy,\bfz \in \TABE 
\) are in the same
 equivalence class if and only if there is a system
of paths \(C(\bfy, \bfz) \) which bounds in \(
\Sigma \). When  \(b_1(Y) = 0 \), this system is unique.

The set of \(\epsilon\)-classes is an affine space modeled on \(
H_1(Y) \cong H^2(Y) \). The choice of a {\it basepoint}
\( z \in \Sigma - \alpha - \beta \) defines a
 correspondence between \(\epsilon\)-classes and \( \spinc \)
 structures on \( Y\). The \( \spinc \) structure associated to an
 \(\epsilon\)-class \(E\) is denoted by \( \spi _z (E) \). To change
 the \( \spinc \) structure we are considering, we can vary either the
 \( \epsilon \)-class, the basepoint, or both. 
This fact is  a very useful feature of the  \OS  Floer
homology. 

\subsection{Domains and differentials}
\label{SubSec:Diffs}
Suppose \( \bfy \) and \( \bfz \) are two generators in the same \(
\epsilon \)-class, so they are connected by a system of curves
\(C(\bfy, \bfz) \) which bounds in \(\Sigma \). Then the set of 2-chains
\(\DDD\) in \( \Sigma \) 
with  \( \partial \DDD = C(\bfy, \bfz)\) is isomorphic to \(
H_2(\Sigma) \cong \Z \). In fact, Lemma 3.6 of \cite{OS1} establishes
a natural correspondance between such \( \DDD \)'s and the
elements of \( \pi_2(\bfy, \bfz)\). Roughly speaking, the chain \(\DDD
(\phi)\) is the image of a \(g\)-fold branched cover of \(\phi \)
induced by the \((g!)\)-fold cover \( \Sigma^g \to s^g \Sigma
\). (This explains the definition of a system of paths: it is just a
\(g\)-fold cover of \( \partial \phi \)). 
  We call the chain corresponding to 
\( \phi \) the {\it domain} of \( \phi\) and denote it by \(\DDD (\phi)\).
Adding a copy of \( \Sigma \) to \( \DDD (\phi) \) corresponds to
connect summing \( \phi\) with the generator of \( \pi_2 (s^g\Sigma)
\). (A more detailed treatment of domains is given in the appendix on
differentials.) 

For \(x \in \Sigma - \alpha - \beta \), we denote by \(
 n_x(\phi) \) the multiplicity of \( \DDD(\phi) \) over the component
 of \( \Sigma - \alpha - \beta \) containing \( \phi \). If \( \phi
 \in \pi_2(\bfy, \bfz) \) has a holomorphic representative, the 
fact that holomorphic maps are orientation preserving shows that 
 \(n_x( \phi ) \geq 0 \) for all \(x\). This is a useful tool
 for showing that a class \( \phi \) does not admit any holomorphic
 representatives.

To each \( \phi \in \pi_2(\bfy, \bfz) \) we associate the {\it
  Maslov index} \( \mu (\phi) \), which is the formal dimension of the
  space of holomorphic disks in the class \( \phi \). The parity of \(\mu
  (\phi) \)  is determined by the intersection number: it is even if \(
  \bfy \)  and \( \bfz\)
both have the same sign of intersection, and odd if they do
  not. In the appendix, we give 
 a combinatorial formula for computing \(\mu (\phi) \)
 from \(\DDD (\phi) \).  It is often useful to know how the Maslov
 index of different elements of \( \pi_2(\bfy, \bfz) \) are related:
  if \( \DDD (\phi') = \DDD (\phi) + n [\Sigma] \), then \( \mu (\phi') = 
\mu (\phi) +  2n \).

Fix an \( \epsilon\)-class \(E\) and a basepoint \(z\). Then the
 complex \( \cfhat (Y, \spi_z(E)) \) is generated by the elements of
 \(E\). The grading and the differentials are defined as follows:
for \(\bfy, \bfz \in E\),
let \( \phi_0 (\bfy, \bfz) \in \pi_2(\bfy, \bfz) \) be the class
 with \(n_z(\phi(\bfy, \bfz)) = 0 \). Then 
\begin{equation*}
 \gr \ts \bfy - \gr \ts \bfz = \mu (\phi_0(\bfy, \bfz)).
\end{equation*}
 If \( \mu(\phi_0(\bfy,\bfz)) = 1\), 
the moduli space \( \MM (\phi_0(\bfy, \bfz)) \) of pseudoholomorphic disks
 in the class \( \phi_0(\bfy, \bfz) \) is generically one dimensional
 and endowed with a free action of \(\R\). 
The \( \bfz \) component of
 \(d(\bfy) \) is defined to be the number of points in the quotient
\( {\MMM} (\phi_0(\bfy, \bfz)) \), counted with
sign. (Throughout this section, we gloss over the  technical
issues of compactness, transversality and orientability for these
moduli spaces. These subjects are important, but their treatment in
the \OS theory is essentially the same as in Lagrangian Floer homology.)

In certain cases the number of points in \( \MMM(\phi) \) depends only on
the topology of the domain \( \DDD(\phi) \). If \(\DDD \) is a region
such that \(\# \MMM(\phi) = \pm 1 \) for any \( \phi \) with \(
\DDD(\phi) = \DDD \), we say that \(\DDD \) is the {\it domain of a
  differential}. The shaded region in Figure~\ref{Fig:GenusTwo} is an
example of such a domain. It defines a differential from \(\{x,y_1\} \)
  to \(\{x, y_2\} \).  

%\begin{figure}
%\caption{\label{Fig:Domains}} Some domains of differentials. 
%\end{figure}

\subsection{\( \cfi\), \( \cfp\), and \( \cfm\)} 

In addition to \( \cfhat \), there are 
 also complexes \(\cfp \), \( \cfm \) and \( \cfi \) whose
generators are obtained by ``stacking'' copies of the generators of
\(\cfhat \). To be precise, the generators of these complexes are of
the form \( \{ [\bfy,n] \ts | \ts \bfy \in \TTT_\alpha \cap \TTT
_\beta \} \), where \(n \geq 0 \) for \( \cfp \), \(n < 0 \) for \( \cfm \)
and \( n \in \Z \) for \( \cfi \).

The differential on \( \cfi (Y) \) is defined in the following way.
If \( \bfy, \bfz \in \TTT_\alpha \cap \TTT_\beta \) have different
  signs, we let \( \phi_{\bfz} \in \pi_2(\bfy, \bfz) \) be the class
  with \(\mu (\phi_\bfz) = 1\). Then
\begin{equation*}
d([\bfy, n] ) = \sum _\bfz \# \MMM(\phi_\bfz) \cdot
 [\bfz, n - n_z(\phi _ \bfz)].
\end{equation*}
\( \cfi (Y) \) has the following basic properties:
\begin{enumerate}
\item It is an \( \text{Affine} (\Z )\) graded chain
  complex, with 
\begin{equation*}
 \gr \ts [\bfy, n] - \gr  \ts [\bfz, m] = \mu (\phi) - 2
  n_z(\phi) +2(n-m) .
\end{equation*}
\item It is translation invariant: {\it i.e.} the map \(u\: \cfi (Y) \to \cfi
(Y)\) is an isomorphism. This gives \( \cfi (Y) \) and \( \cfpm (Y) \)
  the structure of \( \Z[u] \) modules. 
\item The group \( \cfm (Y) \) is a subcomplex, and \( \cfi (Y) / \cfm
  (Y) = \cfp (Y) \). 
\end{enumerate}
The last item implies that there is a long exact sequence
\begin{equation*}
\begin{CD}
@>>> \hfp (Y)  @>>> \hfm (Y) @>>> \hfi (Y) @>>> \hfp (Y) @>>>. 
\end{CD}
\end{equation*}
If \(Y\) is a  homology sphere, the group \(\hfi (Y) \) is
always isomorphic to \( \Z[u, u^{-1}] \). Using this fact and
the long exact sequence,  it is easy to compute \(\hfp \) from \( \hfm
\) and vice-versa. 

Intuitively, we  think of the group \( \hfhat (Y) \) as
being the ordinary homology of a space with an \(S^1\) action, and
\(\hfp (Y) \) as being its equivariant homology. (See \cite{Man} for a
very elegant Seiberg-Witten realization of this idea.) It is easy to see
that the Gysin sequence and spectral sequence relating ordinary
and equivariant homology have analogues which relate \( \hfhat \) to
\(\hfp\).

\subsection{Manifolds with \(b_1 > 0 \) and twisted coefficients}

\begin{figure}
\includegraphics{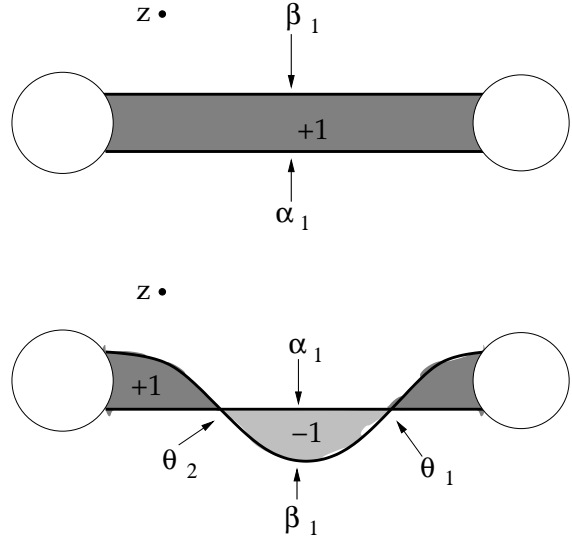}
\caption{\label{Fig:S1S2} Two Heegaard splittings of \(S^1 \times
  S^2\). The multiplicities of the periodic domains are shown.}
\end{figure}

For manifolds with \(b_1 > 0 \), it is no longer true that
\(\pi_2(\bfy, \bfz) \cong \Z \) if \( \bfy \) and \( \bfz \) are in
the same \( \epsilon\)-class. Instead, \(\pi_2(\bfy, \bfz) \cong
\Z\oplus H_2(Y) \). This fact is reflected by the presence of {\it
  periodic domains}: for each class \(x \in H_2(Y) \), there is a
domain \(\PPP (x) \subset \Sigma \) such that \( \partial \PPP (x) \) is
a sum of the \( \alpha _i \)'s and \( \beta _j\)'s. If \( \phi
\in \pi_2(\bfy, \bfz)\), there is a new disk 
\( \phi +x \in \pi_2(\bfy, \bfz)\) with domain \(\DDD(\phi)+ \PPP(x)\).

 To recover the class \(x\)
from \(\PPP (x)\), we put \(\PPP (x) 
 \subset \Sigma \subset Y \) and ``cap off''
each of the \(\alpha_i\)'s and \( \beta _j\)'s in \( \partial \PPP (x)\)
with a disk in the appropriate handlebody. This does not uniquely
determine \(\PPP(x)\), since we can always add copies of \( \Sigma \)
to \( \PPP(x) \) and get something representing the same homology
class. To remedy this problem, we normalize  by requiring that
\(n_z(\PPP(x)) = 0  \). 

In order to define the Floer homology of \( Y\), we must place some
``admissability'' conditions on the behavior of the periodic
domains. This fact is best illustrated by the two Heegaard splittings for
\(S^1 \times S^2\) shown in Figure~\ref{Fig:S1S2}. For the first
splitting, \(\TTT_\alpha \cap \TTT_\beta \) is empty! Since we would
like \( \hfhat (S^1\times S^2) \) to be the homology of the circle
of reducible Seiberg-Witten solutions, this is not very
satisfactory. The solution is to require that all periodic domains
have both positive and negative coefficients, as in  the
second splitting. Now there are two generators \( \theta_1 \) and
\(\theta_2 \) and two separate domains which
contribute to \(d(\bfy, \bfz )\): one  positively, and the
other negatively. Thus we get a complex which looks exactly like the
usual Morse complex for \(H_*(S^1)\).

The periodic domains can be used to define twisted versions of the
homologies discussed above, with  coefficients in \(H^1(Y) \cong
H_2(Y) \). To do this, we fix identifications \(\Phi : \pi_2(\bfy, \bfz)
\cong \Z \oplus H_2(Y)\) which are consistent, in the sense that if
\(\phi_1 \in \pi_2(\bfy, \bfz) \)  and 
\(\phi_2 \in \pi_2(\bfz, {\bf w}) \), then
\( \Phi (\phi_1 + \phi _2) = \Phi (\phi _1) + \Phi (\phi _2) \).
 Then there is a complex \( \underline{\cfhat}(Y)\)
 generated by \((\TABE) \times H_2(Y) \) with differential
\begin{equation*}
d(\bfy,a) = \sum _{\{\phi \in \pi_2(\bfy, \bfz) \ts | \ts n_z(\phi) = 0,
 \ts   \mu (\phi)=1\}} \# \MMM(\phi) \cdot  (\bfz, a + \Phi(\phi)).
\end{equation*} 
For example, if we use the splitting of Figure~\ref{Fig:S1S2}b, \(
\underline{\cfhat}(S^1 \times S^2) \) is isomorphic to the complex
computing \(H_*(S^1) \) with twisted coefficients in \(H^1(S^1)
\). This is the complex which computes the homology of the universal
abelian cover of \( S^1 \), so \( \underline {\cfhat} (S^1 \times S^2)
  \cong \Z \). 

The complexes  \(
\underline{\cfpm}(Y) \) and \( \underline{\cfi} (Y) \) are defined
in an analogous manner.  

\subsection{Maps and holomorphic triangles}

\begin{figure}
\includegraphics{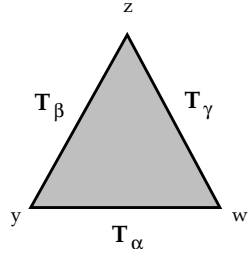}
\caption{\label{Fig:TriangleForm} A model triangle.}
\end{figure}

Maps in the \OS theory are defined by counting holomorphic triangles
in {\it Heegaard triple diagrams}, which have three sets of attaching
handles instead of the usual two. If \( (\Sigma, \alpha, \beta,
\gamma)\) is such a diagram, it determines a cobordism \(W\)  from \(Y_1
\coprod Y_2 \) to \(Y_3 \), where \(Y_1 \) has
Heegaard splitting \( (\Sigma,\alpha, \beta)\), \(Y_2 \) has
 splitting \( (\Sigma,\beta, \gamma)\), and \(Y_3 \) has
 splitting \( (\Sigma,\alpha, \gamma)\). For each \( \spinc \)
 structure \( \spi \) on \(W\), there is a map 
\begin{equation*}
 \hat{F}_{W, \spi} \: \hfhat (Y_1) \otimes \hfhat (Y_2) \to \hfhat (Y_3) 
\end{equation*}
 defined by counting holomorphic triangles in this diagram. 

\begin{figure}
\includegraphics{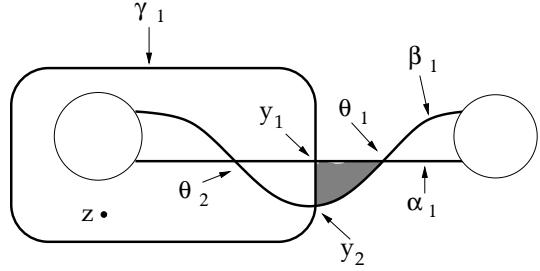}
\caption{\label{Fig:Identity} A Heegaard triple diagram for the
  product cobordism from \(S^3\) to itself.}
\end{figure}

More precisely, let \( \pi _2 (\bfy, \bfz, \bfw ) \) be the set of
homotopy classes of maps of the triangle shown in
Figure~\ref{Fig:TriangleForm} to \(s^g \Sigma \) which take
 vertices and edges  to the points and tori with which they
are labeled. Any \( \psi \in \pi_2(\bfy, \bfz, \bfw ) \) determines a
\( \spinc \) structure \(\spi _z (\psi ) \) on \( W\) which restricts
to the \( \spinc \) structures \( \spi _z (\bfy) \), \( \spi _z (\bfz)
\), \( \spi_z (\bfw) \) on \(Y_1\), \(Y_2\),  \(Y_3\)
respectively. As in the case of differentials, each triangle \( \psi \in
\pi_2(\bfy, \bfz, \bfw ) \) has a domain \( \DDD (\psi) \) which is a
2-chain in \( \Sigma \). Triangles which admit holomorphic
 representatives must have positive domains. 

The map \( \hat{F}_{W, \spi} \) is induced by a map 
\begin{equation*}
\hat{f}_{W, \spi} \: \cfhat (Y_1) \otimes \cfhat (Y_2) \to  \cfhat(Y_3)
\end{equation*}
defined by 
\begin{equation*}
\hat{f}_{W, \spi} (\bfy \otimes \bfz ) = \sum _ {\psi \in A} \# \MMM (\psi)
\cdot \bfw
\end{equation*}
where \(A = \{ \psi \in \pi_2(\bfy, \bfz, \bfw) \ts | \ts n_z(\psi ) =
\mu (\psi) = 0 , \ts \spi_z(\psi) = \spi \} \). There are similarly
defined maps on \( \cfpm \) and \( \cfi \). 

In practice, we are usually
 interested in cobordisms with two boundary components, such as
 surgery cobordisms. These can be represented by Heegaard triple
 diagrams for which \(Y_2 = \# ^k (S^1 \times S^2) \). The original
 cobordism is recovered by ``filling in'' \(Y_2\) with \(
\# ^k (S^1 \times D^3)\).  The map 
\( \hat{g} : \cfhat (Y_1) \to \cfhat (Y_3) \) defined by such a
 cobordism is given by \( \hat{g}(\bfy) = \hat{f} ( \bfy 
 \otimes  {\bf \theta} ) \), where \( {\bf \theta } \) is the generator of \(
 \hfhat (\# ^g (S^1 \times S^2)) \) with the highest homological
 grading. (The reason  is that \(
 {\bf \theta} \) is the relative invariant of \( \# ^k (S^1 \times D^3)\)).
This construction can be useful even when the cobordism in question is
the product cobordism. In fact, this is how \Ozsvath and \Szabo prove
invariance of their Floer homology under handleslides.

As a simple example, consider the genus one triple diagram shown in
Figure~\ref{Fig:Identity}, which represents the product cobordism from
\(Y_1 = S^3\) to \( Y_3 = S^3\). Thus it induces a map \( \hat{f} \: \cfhat
(Y_1) \to \cfhat (Y_3) \).
\( \cfhat (Y_1) \) is generated by \( y_1 \), while \( \cfhat
(Y_3) \) is generated by \( y_2\). 
 \( \hfhat (S^1 \times S^2)\) is generated by
\(\theta _1 \) and \( \theta _2 \), where \( \theta _1 \) has the
higher absolute grading. (We can tell this because the 
differentials go from \(\theta _1 \) to \( \theta _2 \).)  Thus
\( \hat {f}(\bfy) \) is defined by counting triangles in \(\pi_2(y_1,
 \theta _1, y_2)\). There is a unique such triangle; its domain is
shaded in the figure. Using the Riemann mapping theorem, it is easy to
see that this triangle has a unique holomorphic representative. Thus
\(\hat{f} \) is an isomorphism (as it should be). Most of the
holomorphic triangles we need to count will have domains that look
like this one or the union of several disjoint copies of it.

%% file: heegaard.tex
\section{Heegaard splittings and the Alexander grading}
\label{Sec:Heegaard}

This section contains background material on Heegaard splittings for
knot complements, the
\(\epsilon\)-grading, and its relation to Fox calculus and the
Alexander polynomial. We first consider large \(n\)-surgery on
knots in \(S^3\) and then generalize to knots in arbitrary
three-manifolds. Much of this material may be found in 
 \cite{OS1} and \cite{OS2}, although we believe the emphasis on Fox
 calculus is new.

\subsection{ Heegaard splittings and \( \pi_1 \)}

Let \( Y \) be a three-manifold. 
To any Heegaard splitting of \(Y\) there is naturally associated a
presentation of \( \pi_1 (Y) \).  To see this, note that 
\(\pi_1 \) of the handlebody obtained by attaching two-handles along
the \(\alpha\)'s is a free group with \(g\)
generators, and each two-handle \(\beta _i \) gives a relator. More
specifically, we choose as generators of \(\pi_1 (Y) \) a set of
loops \(x_1,x_2,\ldots,x_g\) on \(\Sigma\), where \(x_i\) intersects 
\(\alpha _i \) once with intersection number \(+1\)
 and misses all the other \(\alpha _j\). Then the relator
corresponding to \( \beta _i \)  may be found by traversing \(\beta
_i\) and recording its successive intersections with the \(\alpha
_j\)'s. Each time \( \beta_i \) intersects \( \alpha _j \), we append
\( x_j^{\pm 1} \) to the relator, where \( \pm 1 \) is the sign of the 
intersection.
There a natural correspondence between allowable moves on a
Heegaard splitting and Tietze moves on the associated
presentation. Indeed,  removing a pair of
intersection points by isotopy corresponds to cancelling consecutive
appearances of \(x_i \) and \(x^{-1}_i \) in some word, 
stabilization corresponds to adding a new generator and 
relator,  sliding the handle \(\alpha _i \) over \( \alpha _j \)
corresponds to making the substitution \(x_i' = x_j x_i x_j ^{-1} \),
and  sliding \( \beta _i \) over \( \beta _j \) corresponds to
conjugating \(w_j \) by \(w_j \). 

 Three-manifolds with boundary
also have Heegaard splittings, but with fewer  \(\beta\)
circles than \( \alpha \) circles.
 For example, if \(\partial N = T^2 \), \(N\) has a Heegaard
splitting with one more \( \alpha \) circle than \(\beta \) circles.
 To do Dehn filling on such
a manifold, we  attach a final two-handle along some curve 
\( \beta _g \)  disjoint from \(\beta_1, \beta _2, \ldots \beta _{g-1}
\). This gives  a Heegaard splitting of the resulting closed
manifold. 

It is clear that the correspondence between Heegaard splittings and
presentations of \( \pi _1 \) holds for the case of manifolds with
boundary as well. 

\subsection{Heegaard splittings of knot complements}
\label{Subsec:HeegaardKnot}

Let  \(K\) be 
a knot in \(S^3\). We denote by \(\Ko \) the manifold obtained by
removing a regular neighborhood of \(K\) from \(S^3\). 
In the next two subsections,
we describe this special case in some detail, both to provide background for
later sections and to motivate our treatment of a general three-manifold with
torus boundary. First, we discuss Heegaard splittings of \(\Ko \) and how
to find them. 

We can use a bridge presentation of a knot to find a Heegaard
splitting of its exterior. Recall that \(K\) is said to be
a \(g\)-bridge knot if it has an embedding in \(\R^3\) 
whose \(z\) coordinate has \(g\) maxima.  Any such knot admits a {\it bridge
  presentation} with \(g\) bridges, {\it i.e.} a planar diagram composed
of \(2g\) segments \(a_1,a_2,\ldots a_g , b_1, b_2, \ldots, b_g \)
such that: {\it i)} all of  the \(a_i\)'s and  all of the \(b_i\)'s 
are disjoint, and {\it ii)} the \(b_i\)'s always overcross the \(a_i\)'s. 
Given such a presentation of \(K\), we can obtain a genus \(g\) 
Heegaard splitting for \(\Ko \) as follows: let \( \Sigma \) be the 
surface obtained by  joining the two endpoints of each \(a_i \)
by a tube, and let \( \alpha _i \) be a circle which first traverses
\(a_i\) and then goes over the new tube and back to its starting
point. Finally, let \( \beta _i \) be the boundary of a regular
neighborhood of \(b_i \) in the original diagram. The \(\beta _i\)'s
are linearly dependent; indeed it is easy to see that \(b_1 +b_2 +
\ldots + b_g =0 \) in \(H_1(\Sigma ) \). As a consequence, the
three-manifold obtained by attaching two-handles along any \(g-1\) of
the \( \beta _i \)'s is the same as the one obtained by attaching
two-handles along all \(g\) of them. By convention, we choose to omit
\( \beta _g \), although we could just as
well have skipped any of the others. 

To see why this construction works, consider
the plane of the bridge diagram, which separates \(S^3\) into two
balls. We push the underbridges a little below the plane while leaving
the overbridges in it.  The intersection of \(\Ko \) with the
lower ball is obtained by removing tubular neighborhoods of the
underbridges from the ball. This space is homeomorphic to a handlebody
with boundary \( \Sigma \) and attaching handles \( \alpha _i \). 
To get \(\Ko \), we glue in two-handles along the \(\beta _i \)
(leaving little tubes around each overbridge), and then fill in the 
\(S^2\) boundary component with a ball. 

An example of such a Heegaard splitting is shown in
Figure~\ref{Fig:(3,4)Torus}. In drawing these pictures, we adopt the 
convention of Section 8 of \cite{OS1} and show only the part of \(
\Sigma \) which lies in the plane of the diagram. We do not care much
about the orientation of the \(\beta_i\)'s, but we always orient
the \(\alpha_i\)'s {\it consistantly}, so that \( \beta _i \cdot
(\alpha_1+\alpha _2 +\ldots + \alpha _g) = 0 \) for each \(i \). 
 
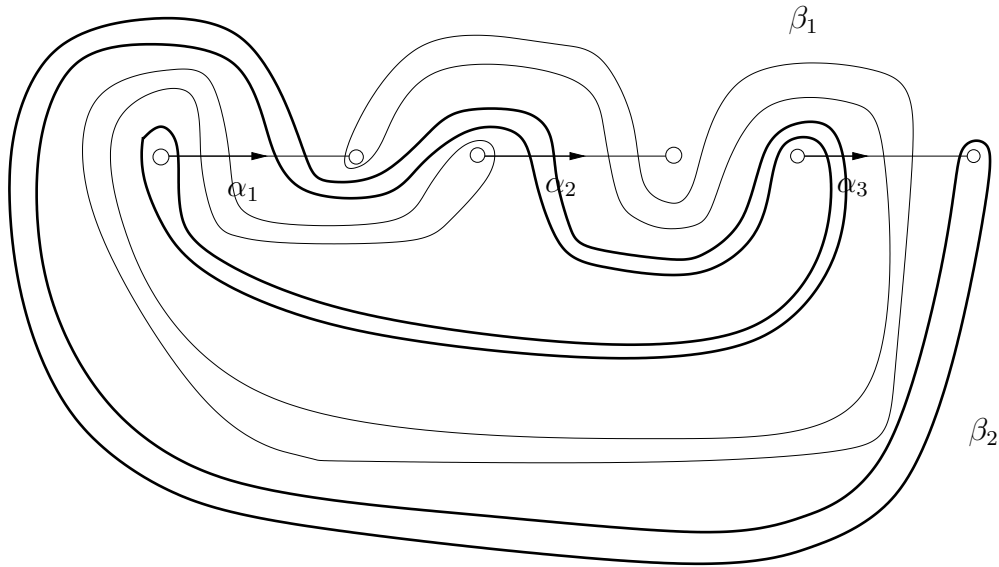
\begin{figure}
\input{figs/34torus.pstex_t}
\caption{\label{Fig:(3,4)Torus} A Heegaard splitting for the
  exterior of the right-handed \( (3,4) \) torus knot. The small
  circles indicate where the tubes are attached and also represent the
  \(x_i\)'s.}
\end{figure}

When our Heegaard splitting comes from a bridge
decomposition, the associated presentation of \(\pi_1 \) is just the
Wirtinger presentation. In particular, all the relators are of the
form \(x_iwx_j^{-1}w^{-1}
\) for some word \(w \). (Our convention for orienting the \( \alpha
_i \)
ensures that \(x_i\) and \(x_j\) have opposite exponents.)
% For example, the relators for Figure~\ref{Fig:(3,4)Torus} are
% \(x_1(x_3^{-1}x_2^{-1}x_1^{-1})x_3^{-1}(x_1x_2x_3) \) and 
%\( x_2(x_1^{-1}x_3^{-1}x_2^{-1})x_1^{-1}(x_2x_3x_1).\)
 Abelianizing, we see that the \(x_i\) are all homologous
to other and generate \( H_1 (\Ko) \cong \Z \). In fact, it is
immediate from our construction of the Heegaard splitting that the
\(x_i\)'s are all meridians of \(K\).  

Given a curve \( \gamma \) in \( \Sigma \), we  often want to
compute its homology class in \(H_1 (\Ko) \). If we let \(\alpha =
\alpha_1 + \alpha_2 + \ldots + \alpha_g \in H_1(\Sigma)\), 
it is easy to see that
\([\gamma] = (\gamma \cdot \alpha ) m\).
For example, if we want to find a longitude of \(K\)  in \( \Sigma \),
we can start with the curve \( \lambda \) obtained by connecting the
\(b_i \)'s to each other by arcs going over the handles. We clearly
have \( m \cdot \lambda = 1 \), and \( \ell = \lambda - (\lambda \cdot
\alpha )m \) is null-homologous in \( \Ko\), so it must be a
longitude.

Any planar diagram of \(K\) gives a bridge decomposition, but some
such decompositions are more suitable for our purposes than others. 
In general, the smaller the bridge number of our presentation
of \(K\), the simpler the complex \(\cfhat(K_n) \) will be. It is thus
 in our interest to be able to find diagrams of a knot with
minimal bridge number. (An exception to this rule may be found in  
section~\ref{Sec:Alternating}, when we discuss alternating
knots.)  For two-bridge knots, this problem was 
solved by Schubert in \cite{Schu}. Unfortunately, there is no such
nice description of \(n\)-bridge knots, even for \(n=3\).

 Starting from  a
diagram of \(K \) with \(n\) obvious maxima and minima, there is a
 straightforward algorithm to find a bridge decomposition by
 successively ``unbraiding'' \(K\). (It is a good exercise to derive
 the splitting of Figure~\ref{Fig:(3,4)Torus} using this method.) In 
practice, however, this method requires a lot of patience and chalk.
It is usually much easier to use the converse, which implies that
every Heegaard splitting of the form described above is  a
splitting of some knot complement.
Given a knot \(K\), we  write down the Wirtinger presentation
for \(\pi_1 (\Ko) \) and use a computer algebra system to simplify
the presentation, eliminating generators while keeping the relators in
the form \(x_iwx_j^{-1}w^{-1}\). Once we have reduced the presentation
to as few generators as possible, we can just try to draw a Heegaard
splitting which gives the simplified presentation. {\it A priori},
of course, we only know that we have drawn a knot complement with the
same \( \pi_1 \) as our original knot. But under very mild hypotheses
(such as \(K\) being hyperbolic), this is enough to ensure that the
knot complement we have drawn corresponds either to \(K \) or its
mirror image. 

\subsection{The Alexander grading for knots}
\label{Sec:AlexGrading}

We now study the effect of doing \(n\) surgery on \(\Ko \) to get the
closed manifold \(K(n,1) \). 
 The final attaching circle \( \beta _g \) will be a curve in \(\Sigma -
\beta _1 - \beta _2 - \ldots - \beta _{g-1} \) homologous to \( n m +
l \). There are many such curves; we choose the one obtained by 
taking the union of \( \ell \) and \(n \) parallel copies of \(x_g
\) and smoothing the intersections to get a simple closed curve
 \( \beta _g \). This procedure was introduced in section 8 of
 \cite{OS1}; we refer to it as {\it twisting up} around \(\alpha_g \).
Of course, our choice of \(x_g \) was arbitrary --- we could
just as well have twisted up around any of the other \(\alpha _i \).  

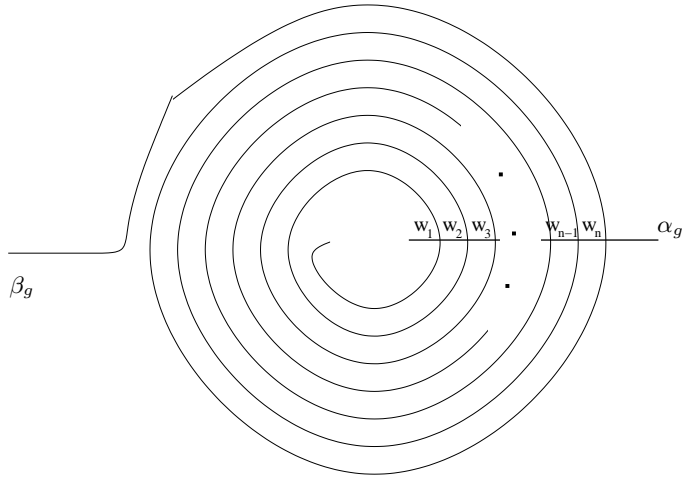
\begin{figure}
\input{figs/spirallh.pstex_t}
\caption{\label{Fig:SpiralRegion} The spiral region of \( \beta _g
  \). Note that the handedness of the spiral is determined by the fact
  that we are doing \( +n \) surgery (and the convention that \( m
  \cdot \lambda = +1 \).)}
\end{figure}

Although \( \beta _g \) can be  complicated, the only part of it
which will be relevant to us is the spiral region containing the \(n\) 
parallel copies of the meridian, which is illustrated
 in Figure~\ref{Fig:SpiralRegion}. We label the intersections
between \( \alpha _g \) and \(\beta _g \) by \(w_i \). To do this, we
must specify the direction of increasing \(i\), or equivalently, an
orientation of \( \ell \). To this end, we fix once and for all a
generator of \( H_1(K(0,1)) \). Given a Heegaard splitting of \(\Ko\),
this choice  determines an orientation on \(m\), and thus on \(\ell\).

Our first class is to understand the \( \epsilon\)-grading 
 associated to this Heegaard diagram.
The following observation is trivial, but quite important.
\begin{lem}
For \(n\) sufficiently large, there is an \(\epsilon\)-class all of
whose elements contain some \(w_i\).
\end{lem}

\begin{proof}
The number of \(\epsilon\)-classes is \(n\), but the
 number of intersection points \({\bf y} \) which do not contain
any \(w_i \) is bounded independent of \(n\).
\end{proof}

We will always  work with such an \(\epsilon\)-class. Thus
we need only consider elements of \( \TABE \) which have
 the form \( {\bf y} = 
\{ \hat{\bf y} , w_i \} \), where \(\hat{\bf y} \) is an intersection
point of the tori \(\TTT _{\hat{\alpha} }\) and \(\TTT _{\hat{\beta}} \)
defined by the \((g-1)\)-fold symmetric products of \( \alpha _1,
\alpha _2, \ldots, \alpha _{g-1} \) and \( \beta _1, \beta _2, \ldots,
\beta _{g-1} \) in \( s^{g-1}\Sigma \). 

When \(n\) is large, most \(\epsilon \)-classes
contains precisely one point of the form \( \{ \hat{\bf y} , w_i \} \)
for each \(\haty \in \TAB \). Indeed, it is easy to see from
Figure~\ref{Fig:SpiralRegion} that 
\begin{equation*}
  \epsilon( \{\haty, w_i\}) - \epsilon (\{ \haty , w_j \}) = (j-i)[m],
\end{equation*}
so a given \( \epsilon\)-class can contain no more than one such
point. On the other hand
\begin{equation*}
  \epsilon( \{\haty, w_i\}) - \epsilon (\{ \hatz , w_i \}) = k[m],
\end{equation*}
for some \(k\), since \([m]\) generates \(H_1(K_n) \). Then by the
additivity of \(\epsilon \), 
\begin{equation*}
  \epsilon( \{\haty, w_i\}) - \epsilon (\{ \hatz , w_{i-k} \}) = k[m]
  -k[m] = 0,
\end{equation*}
so \( \{\haty, w_i\} \) and \( \{ \hatz , w_{i-k} \} \) are in the
same \( \epsilon\)-class. (Our labeling is modulo \(n\), so that if
\(i-k \leq 0 \), we interpret \(w_{i-k}\) as \(w_{i-k+n} \).) 

Now it is easy to see that if  \( \{ \hat{\bf y} , w_{A(\haty )} \}
\) are the points in one \( \epsilon \)-class, 
 \( \{ \hat{\bf y} , w_{A(\haty )+k} \}\) will be the points in
another. Thus most \(\epsilon \)-classes look very similar --- they
are just translates of each other. When the values of \(A(\haty )+k \)
get close to \(1 \) or \(n \), there will be some ``wrapping'' from
one side of the spiral to another, but when \(n\) is large, most \(
\epsilon\)-classes will have values of \(A(\haty )+k \) far away from
\(1\) and \(n\). 
To summarize, we have the following

\begin{prop}
\label{Prop:eClasses}
There is a number \(M\) independent of \(n\)  so that all but \(M\)
\(\epsilon\)-classes are of the form
\( \{\{\haty, w_{k+A(\haty )}\} \ts | \ts \haty \in
T_{\hat{\alpha}} \cap T_{\hat{\beta}} \} \)
 for some value of \(k\). 
\end{prop}

\noindent The function \( A\) is most naturally thought of as an 
affine-valued grading
\begin{equation*}
A:T _{\hat{\alpha} } \cap T_{\hat{\beta}}\to \text{\it Affine} (\Z),
\end{equation*}
which we refer to as the {\it Alexander grading}. Our choice of name
is explained by

\begin{prop}
\label{Prop:AlexP}
Let \( {\sign} (\haty ) \) denote the sign of the intersection between
\( T_{\hat{\alpha}} \) and \( T_{\hat{\beta}} \) at \( \haty \). Then
\begin{equation*}
\sum _{\haty \in \TAB} {\text sign}  (\haty ) \ t^{A(\haty )}
\end{equation*}
represents the Alexander polynomial \(\Delta _K (t)\)  of \(K\). 
\end{prop}

{\bf Remark:} The Alexander polynomial is defined only up to a factor
of \( \pm t^k \), reflecting the fact that \( A \) is only an affine
grading. The symmetry of  \(\Delta _K (t)\) gives us a natural
choice of representative, however, namely the one for which \( \Delta
_K (t) = \Delta _K (t^{-1}) \) and \( \Delta _K (1) = 1 \). We will
always assume \( \Delta _K (t) \) satisfies these properties. To
indicate the weaker condition \(P(t) = \pm t^k \Delta _K (t) \), we
write ``\(P(t)\) represents \( \Delta _K(t) \),''  or simply \( P (t) 
\sim \Delta _K (t) \). 

Our choice of a distinguished representative for the
Alexander polynomial gives us a natural lift of \(A\) to a \(\Z \)
valued map; namely the one for which the sum above is actually equal
to \( \Delta _K (t) \). From now on we use this lift,
which we continue to denote by \(A\). 

Proposition~\ref{Prop:AlexP}  is a corollary of the following principle,
which is very helpful in making computations:

\begin{principle}
The process of computing the elements of \( \TAB \) and their
Alexander gradings is identical to the process of computing the
Alexander polynomial of \(K\) by Fox calculus. 
\end{principle}

More precisely,  our choice of Heegaard splitting naturally
gives us a presentation \(P = \langle x_1, x_2, \ldots, x_g \ts | \ts w_1,
w_2, \ldots, w_{g-1} \rangle \) of \(\pi _1 (\Ko ) \). 
Recall from \cite{CrowFox} that to compute the Alexander polynomial of
\(K\) from \(P\), we take the determinant
\begin{equation*}
(d_{x_i} w_j) _{1 \leq i,j \leq g-1}.
\end{equation*}
of the matrix of free differentials. If we take the free
differentials, expand the determinant, and multiply out {\it without
ever combining terms}, the monomials of the resulting Laurent
polynomial will naturally correspond to the points of \( \TAB
\), with their exponents giving the Alexander grading and their signs
giving the sign of intersection.

\subsection{Fox Calculus and Closed Manifolds}

We now extend the ideas of the preceding section to more general
three-manifolds. First, we describe the relationship between Fox
calculus and the \(\epsilon\)-grading on a closed manifold \( Y\).
 We choose a
Heegaard splitting  \( (\Sigma, \balpha, \bbeta ) \) of \(Y\)
 and take \(P\) to be the associated presentation of \(\pi _1
(Y) \).  

\begin{principle}
\label{Prin:Alex0} The process of computing the points of \(\TTT_{\alpha}
\cap \TTT_{\beta} \) and their \( \epsilon\)-gradings is identical to the
process of computing the \(0\)-th Alexander ideal associated to the
presentation \(P\).
\end{principle}

Actually, we need to use a slight variation of the usual version of Fox
calculus: we take coefficients in the group ring \( \Z [H_1(Y)] \)
rather than in  \( \Z [H_1(Y) / {\text{\it Torsion}}] \). (The first
ring is usually used because one needs a UFD  to define the
greatest common divisor.  Since we will never
 take \(\gcd\)'s, there is no problem with using the larger ring.)

The first step in establishing the principle is to observe that the
\(\epsilon\)-grading satifies a slightly stronger version of additivity.

\begin{lem}
\label{Lem:EpsAdd}
For points \({\bf y} = \{y_1,y_2,\ldots, y_{g}\} \) and 
\({\bf y}' = \{y_1',y_2,\ldots, y_{g}\} \) in \( \TABE \), the
 difference 
\( \epsilon(\{y_1,y_2,\ldots, y_{g}\}) - \epsilon(\{y_1',y_2,\ldots, y_g\}) \)
depends only on \(y_1 \) and \(y_1'\). 
\end{lem}

\begin{proof}
Suppose \(y_1 \) is an intersection point of \( \alpha _i \) and \(
\beta _j \). Then \(y_1'\) is another intersection point of \( \alpha
_i \) and \(\beta _j \), since these are the only curves which do not
contain one of \( y_2, \ldots, y_{g} \). To compute 
\( \epsilon(\{y_1,y_2,\ldots, y_{g}\}) - \epsilon(\{y_1',y_2,\ldots, y_g\})
\), we consider the system of loops in \( \Sigma \) obtained by
going from  \(y_1 \) to \(y_1'\) along \( \alpha _i \) and 
from \( y_1'\) to \(y_1\) along \( \beta _j \), and then joining each
\(y_i\), \( i \geq 2 \) to itself by the trivial loop. The homology
class of this system is just the class of the first loop,
which depends only on \( y_1 \) and \(y_1 '\).  
\end{proof}

Thus for each \(i \) and \(j\) we have a well-defined affine grading
\begin{equation*}
 \epsilon _{ij} :  \alpha _i \cap \beta_j \to {\text{\it Affine}}(H_1(N))
\end{equation*} 
given by the homology class of this loop.
 Repeated application of the lemma shows that if \(\{y_1,y_2, \ldots
y_{g}\} \) and \(\{y_1',y_2', \ldots y_{g}'\} \) share the same
 combinatorics, {\it i.e.} \(y_k \) and \(y_k'\) are both intersection
 points of the same \( \alpha_{i_k} \) and \( \beta _{j_k} \), then 
\begin{equation*}
A(\{y_1,y_2,\ldots, y_{g}\}) - A(\{y_1',y_2',\ldots, y_{g}'\}) =
\sum _{k=1}^{g} \epsilon_{i_kj_k}(y_k) - \epsilon_{i_kj_k}(y_k')
\end{equation*}

We now relate \( \epsilon_{ij} \) to the free differential.

\begin{lem}
\label{Lem:FreeDiff}
\begin{equation*}
d_{x_i} w_j \sim \sum _{y \in \alpha _i \cap \beta _j}
 \text{sign}(y) \ {\epsilon_{ij} (y)}
\end{equation*} 
\end{lem}

\begin{proof}
Each intersection \(y \in \alpha _i \cap \beta _j \) corresponds to
an appearance of \( x_i ^{ \text{sign}(y)}\) in \(w_j \), and thus to
a monomial \(\sign (y) f(y) \) in 
\( d_{x_i} w_j \). We claim that \(f(y) \sim \epsilon_{ij} (y)
\). Indeed, suppose that  
\(y_1\) and \( y_2\) are two elements of \( \alpha _i \cap \beta _j
\).  
 Assume first that both intersections are positive, so
that \(w_j \) has a segment of the form \(x_i x^{a_1}_{i_1}
x^{a_2}_{i_2} \cdots x^{a_n}_{i_n} x_i \), where the two appearances
of \(x_i \) correspond to \(y_1 \) and \( y_2 \) respectively. Then
\(f(y_2) - f(y_1) \) is the image of \(x_i x^{a_1}_{i_1}
x^{a_2}_{i_2} \cdots x^{a_n}_{i_n}\) in \(H_1 (N) \). 
 On the 
other hand, \(\epsilon_{ij}(y_2) - \epsilon_{ij}(y_1) \) is the
homology class of the loop which starts just after the first \(y_1 \),
travels along \( \beta _j \) until just before the second \(y_2\), and
then returns to its starting point along \( \alpha _i \), going up once
through \( \alpha _i \) in the process. Thus this loop is also
represented by  
\(x_i x^{a_1}_{i_1}x^{a_2}_{i_2} \cdots x^{a_n}_{i_n}\).
This proves the lemma when both intersections are
positive. The other cases are similar.
\end{proof}

%\begin{figure}
%\caption{\label{Fig:FreeDiff} ???????}
%\end{figure}

The \(0\)-th Alexander ideal  is generated by 
 the determinant of the matrix of free differentials:
\begin{equation*}
(d_{x_i} w_j) _{1 \leq i,j \leq g}.
\end{equation*}

\begin{prop}
\label{Prop:MatrixDet}
\begin{equation*}
\det (d_{x_i} w_j)  \sim \sum _{{\bf y} \in \TABE} {\text sign} \
({\bf y} ) \  {\epsilon({\bf y} )}
\end{equation*}
\end{prop}

\begin{proof}
We expand the determinant as a sum 
\begin{equation*}
\sum _{\sigma \in S_{g}} \sign (\sigma ) \  d_1 w_{\sigma (1)} d_2 
w_{\sigma (2)} \cdots d_{g} w_{\sigma (g)}
\end{equation*}
 and multiply out  without ever
 combining terms. We claim that the monomials in the resulting polynomial
 correspond precisely to the points of \(\TABE\). Indeed, to specify a
 point of \( \TABE \), we must first choose a partition of 
the \( \alpha _i \)'s and \(\beta_j \)'s into \(g \) sets, each
 containing one \( \alpha _i \) and one \( \beta _j \), which
 corresponds to picking a permutation \( \sigma  \in S_g\).
  Then for each pair \(
 \{ \alpha _i , \beta_{\sigma(i)} \} \) we must choose a point in \( \alpha
 _i \cap \beta _{\sigma (i)} \). This corresponds to picking a monomial out
 of each term \( d_{x_i} w_{\sigma(i)} \) appearing in the product, or
 equivalently, a single monomial from the expanded product. Thus 
 each \( \bfy \in \TABE \) is associated to some monomial \( \pm
 f(\bfy) \)  in the expansion of the determinant. 

We claim that for \( \bfy, \bfz \in
 \TABE \), \( f(\bfy ) - f(\bfz) = \epsilon (\bfy) - \epsilon (\bfz)
 \). 
 Lemmas~\ref{Lem:EpsAdd} and
\ref{Lem:FreeDiff} show that this is true for intersection points
\( \haty \) and \( \hatz \) which have the same combinatorics.
To check it in general, we consider the loop representing 
  \( \epsilon(\{y_1,y_2,y_3,\ldots, y_{g}\})
- \epsilon(\{y_1', y_2 ', y_3, \ldots, y_{g} \}) \)
and argue, as in the proof of Lemma~\ref{Lem:FreeDiff}, that this loop
 also represents the difference appearing in the free differential.
 This proves the claim when \(g-2 \) of the points
in \({\bf y} \) and \( {\bf y}' \) are the same, and the
general result follows by repeated application of this fact. (If
necessary, we introduce some extra pairs of intersection points to
ensure that each product in the determinant is nonempty.)

%\begin{figure}
%\caption{\label{Fig:Alex2} Path for computing \(
% \epsilon(\{y_1,y_2,y_3,\ldots, y_{g}\}- \epsilon(\{y_1', y_2 ', y_3, \ldots,
% y_{g}\} \) }
%\end{figure}

Finally, we check that the sign of the monomial associated to \( \bfy
\) is the sign of the intersection. 
 Locally, each \( {\bf y} \in \TABE \) looks like a product of 
intersections \( \alpha _i \cap \beta _{\sigma (i) } \) in \( \Sigma
\). Now 
\begin{equation*}
(\alpha _1 \times \alpha _2 \times \ldots \times \alpha _{g} ) \cdot 
(\beta _{\sigma(1)} \times \beta _{\sigma (2)} \times \ldots \times 
\beta _{\sigma(g)}) = 
(\alpha _1 \cdot \beta _{\sigma(1)})(\alpha _2 \cdot \beta _{\sigma
  (2)})
\cdots (\alpha _{g} \cdot \beta _{\sigma (g)})
\end{equation*}
is the sign of the monomial corresponding to \(\bfy \) in 
\(d_1 w_{\sigma (1)} d_2 w_{\sigma (2)} \cdots d_{g} w_{\sigma
  (g)} \). Since the product orientation on 
\( (\beta _{\sigma(1)} \times \beta _{\sigma (2)} \times \ldots \times 
\beta _{\sigma(g)}) \) is \( \sign (\sigma ) \) times the usual
orientation, the signs are correct.
\end{proof}

Thus we have verified Principle~\ref{Prin:Alex0}. An analogue of
Proposition~\ref{Prop:AlexP} holds as well, but is rather
uninteresting, since it amounts to the well-known fact that \( \chi
(\cfhat (N) ) \) is \(1\) if \(b_1 (N) = 0 \) and \(0 \) otherwise. 

\subsection{The Alexander grading in general}
\label{Ssec:GenAlex}

We now describe the analogues of Propositions~\ref{Prop:eClasses} and
 \ref{Prop:AlexP} for a knot in a general
 three-manifold. We prefer to think about the knot complement \(N\),
 which is a manifold with torus boundary. Any three-manifold with
 torus boundary has a class \( l \in H_1 (\partial N) \) which bounds
 in \(N\); if \(N\) is a knot complement, \(l\) may be chosen to be
 primitive in \( H_1 (\partial _N) \). There are two possible choices
 of such an \( l\) --- we fix one  of them. 
Next, we choose a class   
\( m \in H_1 (\partial N) \) with 
\( m \cdot l = 1\) to play the role of the meridian. Unlike the choice
 of \(l\), which is just a sign convention, this choice of \(m\) is 
  very important. Everything that follows  will depend on
 which \(m\) we pick. We write \(N(p,q) \) for the closed manifold
 obtained by doing Dehn filling on \(pm + ql \). We will often need
 to consider \(N(1,0) \), which we denote by \(\Nbar \). 

We now proceed much as we did for knots in \(S^3\). We choose a Heegaard
splitting \( (\Sigma, \balpha, \bbeta ) \) of \(N\)  and represent
\(m\) and \(l\)
by curves \( \emm \) and \( \ell \) 
on \( \Sigma \) disjoint from the \( \beta\)'s. The most
significant difference from the knot case is that \(\emm \) may have many
intersections with the \( \alpha _i\)'s, rather than a single one. If we wish,
however, we can always reduce to the latter case by increasing the
genus of our Heegaard splitting. 

\begin{lem}
\label{Lem:Meridian}
Any Heegaard splitting of \(N \) is equivalent to one in which \(\emm \)
intersects a single \( \alpha _i \) geometrically once. 
\end{lem}

\begin{proof}
Stabilize the Heegard splitting by connect summing \( \Sigma \) with 
the genus 1 Heegaard splitting for \(S^3 \). Call the new attaching
handles \( \alpha _{g+1} \) and \( \beta _g \). Then we can slide \(
m\) over \( \beta _g \) so it has a single intersection with \( \alpha
_{g+1} \). Now remove the intersections of \( m \) with the other \(
\alpha _i \) by sliding them over \( \alpha _{g+1}\). 
\end{proof}

\noindent Despite this fact, we will continue to work in the more
general setting, since it is needed for the proof of
 Theorem~\ref{Thm:ReducedInvariance}. 

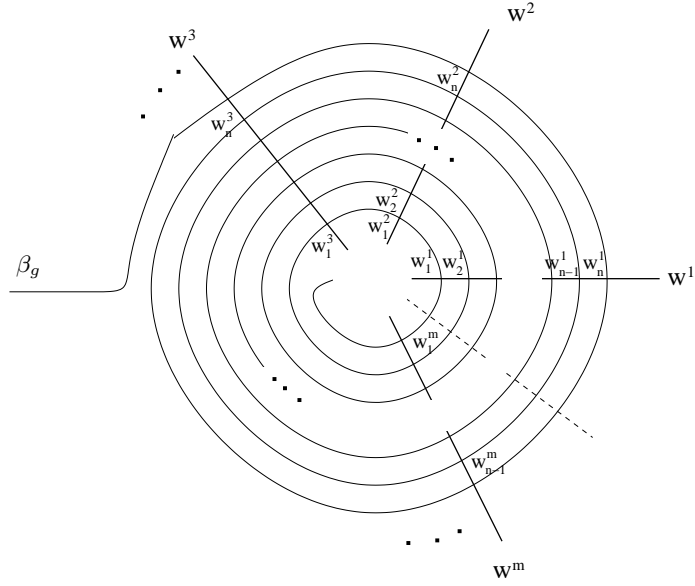
\begin{figure}
\input{figs/genspirallh.pstex_t}
\caption{\label{Fig:GeneralSpiral} 
The spiral region for a general Heegaard splitting}
\end{figure}

To get a Heegaard splitting of \(N(n,1)\), we let 
  \( \beta _g \) be the curve obtained by taking the union of \(\ell \) with
\(n \) parallel copies of \(\emm \) and then resolving. 
As in the knot case, we take \(n \gg 0 \) and 
focus on the spiral region shown in
Figure~\ref{Fig:GeneralSpiral}. We label the  segments of
the \( \alpha \) curves passing through the spiral \(w^1,w^2, \ldots,
w^m\). We will need to know which \( \alpha _j \) contains \(w_i\), so we
  fix numbers \(m_i\) so that \(w^i \subset \alpha _{m_i} \). 
Finally, we label the individual intersection
  points on by \(w^i_j \), using our chosen orientation on \( l \) to
  determine the direction of increasing \(i\).

When \(n\) is large, most points in \( \TTT_{\alpha} \cap \TTT_{\beta} \)
 will be of the form \( \bfy = \{ \haty, w^i_j \} \) for
 \( \haty \in \TTT_{\widehat{\alpha _{m_i}}} \cap \TTT_{\hat{\beta}}
 \), where  \( T_{\widehat{\alpha _{m_i}}} \) is the torus
 in \( s^{g-1} \Sigma \) determined by the symmetric product of 
\(\alpha_1, \alpha _2, \ldots, \hat{\alpha_i},\ldots \alpha_ g\).
We often  think of the point \( \bfy \) as lying over \( w^i_j\),
 and write \( w(\bfy) = w^i_j\).

To each \( \bfy = \{ \haty, w^i_j \}\),
 there is an associated point \( \overline{\bfy} = \{ \haty, w^i \}
 \) in \( \TTT_{\alpha} \cap \TTT_{\bbbar } \), where \((\Sigma, \balpha,
 \bbbar ) \) is the Heegaard splitting for \(\Nbar \) obtained by
 taking \( \bbar _g  = m \). We let \(\ebar \) denote the
 \(\epsilon \)-grading on the \(\{\haty, w^i \} \) associated to this
 Heegaard splitting. It is easy to see
 that \( \ebar (\overline{\bfy}) \) is the image of \( \epsilon (
\bfy ) \) under the projection
 \begin{equation*}
 H_1 (N) / n\langle m
 \rangle \to H_1(N)/ \langle m \rangle \cong H_1(\Nbar). 
\end{equation*}

There are \(\Z /n \) \( \epsilon\)-classes over each \( \ebar\)-class,
and if \( \{\haty, w^i \} \) is in an \( \ebar\)-class, one of the 
\( \{\haty, w^i_j \} \) is in each of the \( \epsilon \)-classes above
it.

\begin{prop}
\label{Prop:GoodEClass}
There is an number \(M\) independent of \(n\) such that for any \(
\ebar\)-class \(\Ebar \), all but \(M\) of the \( \epsilon\)-classes above
\( \Ebar\) are of the form
\begin{equation*}
\{ \{ \haty, w^i_{k+A(\haty, w^i)} \} \ts | \ts \{ \haty, w^i\} \in \Ebar \}
\end{equation*}
\end{prop}

\begin{proof}
The number of points in \( \TABE \) not of the form 
\(  \{ \haty, w^i_{k+A(\haty, w^i)} \} \) is bounded independent of
\(n\). The number of \( \epsilon\)-classes above \(\Ebar\) in which
the points ``wrap'' (so that some of them are on the left end of the
spiral while others are on the right) is also bounded independent
of \(n\). Finally, \( \epsilon (\{\haty, w^i_j ) - \epsilon
(\{\haty, w^i_l \}) = (l-j)m \), so if 
\(\{ \haty, w^i_{A(\haty, w^i)} \} \) are all in  one \(
\epsilon\)-class, so are \(\{ \haty, w^i_{A(\haty, w^i)} \} \).
\end{proof}

We call a class of this form as a {\it good \(\epsilon\)-class}.
The map \(A: \Ebar \ \to {\text{\it Affine}}(\Z) \)  is the
Alexander grading. The sign of the Alexander grading depends
on our choice of orientation for \(l\) --- reversing the orientation
gives the opposite sign. (This reflects the fact that our
representation of the Alexander polynomial depends on a choice of
basis for \(H_1 \).) To summarize, we have
% Note that it only makes sense to talk about the
%difference in Alexander grading of two points \(\{ \haty, w^i \} \)
%and \(\{\hatz, w^j \} \) if they are in the same \(\ebar\)-class.

\begin{cor}
\label{Cor:InSpiral}
Any two points in a good \( \epsilon \)-class are joined by
null-homologous system of paths whose \( \beta _g \) component is
 supported inside the
spiral region. The difference in their Alexander gradings is just the
number of times the \( \beta _g \) component intersects the dashed
line in Figure~\ref{Fig:GeneralSpiral}. 
\end{cor}

\begin{proof}
The corresponding points in \( T_{\alpha} \cap T_{\bbar} \) are in
the same \( \ebar \)-class, so they can be joined by a null-homologous
 system of paths
along \( \alpha _1, \alpha _2, \ldots, \alpha _g, \beta _1, \beta _2,
\ldots, \beta _{g-1}, \emm \). The path for the original points is the
same, but with the \(\emm\) component lifted up to the spiral region.
\end{proof}

\begin{prop} The process of finding the \(\ebar \)-classes and 
Alexander gradings of all the points \(\{ \haty, w^i\} \) is the same
as that of computing
\begin{equation*}
\det (d_i w_j) _{1\leq i,j \leq g }
\end{equation*}
with coefficients in \( \Z[H_1 (N)] \). 
\end{prop}

\begin{proof}
If we reduce coefficients to \( \Z [H_1(\Nbar)] \),
 Principle~\ref{Prin:Alex0} shows that
 the  determinant computes the \( \ebar\)-gradings of the \( \{\haty,
 w^i \} \). With in \( \Z[H_1 (N) ] \), nearly
 everything is the same, but
 the difference in the group coefficients of the monomials
 corresponding to points  \( \{\haty,w^i \} \) and \( \{\hatz, w^j \}
 \) will be the class in 
\begin{equation*}
 H_1 (\Sigma)/\langle \alpha_1, \alpha _2,
 \ldots, \alpha _g, \beta _1, \beta _2, \ldots, \beta _{g-1} \rangle \cong H_1
 (N) 
\end{equation*}
 of a system of loops joining \( \{\haty,w^i_k \} \) to 
\( \{\hatz, w^j_k \}\), where we require the part of the loop on \(\beta_g \)
to stay inside the spiral. (This is why we do not have to quotient out
 by \(\beta _g \)). Thus if \( \{\haty,w^i \} \) and \( \{\hatz, w^j \}
 \) are in the same \( \ebar\) class, Corollary~\ref{Cor:InSpiral} shows
 that the difference in their group
 coefficients is just \(m \) times the difference in their Alexander
 gradings.  
\end{proof}

This proposition can be rephrased by saying that there is a lift of
the \( \epsilon \)-grading on \(\Nbar \) to a function
\begin{equation*}
\eta \: \TTT _{\alpha} \cap \TTT_{\bbar} \to {\text{Affine}}(H_1(N))
\end{equation*}
We call \( \eta \) the global Alexander grading. 
Note that it only makes sense to talk about the difference of
Alexander grading of two points as an integer when their their global
Alexander gradings reduce to the same class in \( H_1(\Nbar) \), {\it
  i.e.} when they belong to the same \( \ebar \)-class.  

\begin{cor}
\label{Cor:AlexP}
\begin{equation*}
\sum _{{\bf y} \in T_{\alpha} \cap T_{\bbar}} {\text{sign}}  ({\bf y})
 \ \eta({\bf y})
\end{equation*}
represents the Alexander polynomial of \(N\). 
\end{cor}

\begin{proof}
The previous proposition implies that 
\begin{equation*}
\sum _{{\bf y} \in T_{\alpha} \cap T_{\bbar}} {\text{sign}}  ({\bf y})
 \ \eta({\bf y}) = \det (d_i w_j) _{1\leq i,j \leq g }
\end{equation*}
The usual arguments from Fox calculus show that this determinant is
invariant under the operations of stabilization, sliding one \(
\alpha_i \) over another, and sliding \(m\) over \( \beta _i\).
 Thus by  Lemma~\ref{Lem:Meridian} it suffice to prove the result
 when \( m \) has a unique
intersection with  \( \alpha _g \). In this case the determinant reduces to
\begin{equation*} 
 \det (d_{x_i}w_j)_{1\leq i,j \leq g-1}
\end{equation*}
Now to compute the Alexander polynomial we usually need to take the
\(\gcd \) of all the polynomials 
\begin{equation*}
P_k = \det (d_{x_i}w_j)_{j\neq k}
\end{equation*}
When the presentation comes from a Heegaard splitting, however,
all of the \(P_k\)'s divide each other, so we are done. 
(See, {\it e.g},  Theorem 5.1 of
\cite{McMullen}.)
\end{proof}

%% file: figs/34torus.pstex_t
\begin{picture}(0,0)%
\includegraphics{figs/34torus.pstex}%
\end{picture}%
\setlength{\unitlength}{1460sp}%
\begingroup\makeatletter\ifx\SetFigFont\undefined%
\gdef\SetFigFont#1#2#3#4#5{%
  \reset@font\fontsize{#1}{#2pt}%
  \fontfamily{#3}\fontseries{#4}\fontshape{#5}%
  \selectfont}%
\fi\endgroup%
\begin{picture}(16734,9493)(-1415,-10931)
\put(12676,-4561){\makebox(0,0)[lb]{\smash{\SetFigFont{12}{14.4}{\rmdefault}{\mddefault}{\updefault}{\color[rgb]{0,0,0}\( \alpha _3\)}%
}}}
\put(7726,-4561){\makebox(0,0)[lb]{\smash{\SetFigFont{12}{14.4}{\rmdefault}{\mddefault}{\updefault}{\color[rgb]{0,0,0}\( \alpha _2 \)}%
}}}
\put(11862,-1786){\makebox(0,0)[lb]{\smash{\SetFigFont{12}{14.4}{\rmdefault}{\mddefault}{\updefault}{\color[rgb]{0,0,0}\( \beta_1 \)}%
}}}
\put(14922,-8791){\makebox(0,0)[lb]{\smash{\SetFigFont{12}{14.4}{\rmdefault}{\mddefault}{\updefault}{\color[rgb]{0,0,0}\( \beta _2\)}%
}}}
\put(2326,-4636){\makebox(0,0)[lb]{\smash{\SetFigFont{12}{14.4}{\rmdefault}{\mddefault}{\updefault}{\color[rgb]{0,0,0}\(\alpha _1\)}%
}}}
\end{picture}

%% file: figs/spirallh.pstex_t
\begin{picture}(0,0)%
\includegraphics{figs/spirallh.pstex}%
\end{picture}%
\setlength{\unitlength}{2171sp}%
\begingroup\makeatletter\ifx\SetFigFont\undefined%
\gdef\SetFigFont#1#2#3#4#5{%
  \reset@font\fontsize{#1}{#2pt}%
  \fontfamily{#3}\fontseries{#4}\fontshape{#5}%
  \selectfont}%
\fi\endgroup%
\begin{picture}(7898,5387)(-1425,-6043)
\put(6436,-3271){\makebox(0,0)[lb]{\smash{\SetFigFont{9}{10.8}{\rmdefault}{\mddefault}{\updefault}{\color[rgb]{0,0,0}\(\alpha_g\)}%
}}}
\put(-959,-3961){\makebox(0,0)[lb]{\smash{\SetFigFont{9}{10.8}{\rmdefault}{\mddefault}{\updefault}{\color[rgb]{0,0,0}\(\beta_g\)}%
}}}
\end{picture}

%% file: figs/genspirallh.pstex_t
\begin{picture}(0,0)%
\includegraphics{figs/genspirallh.pstex}%
\end{picture}%
\setlength{\unitlength}{2171sp}%
\begingroup\makeatletter\ifx\SetFigFont\undefined%
\gdef\SetFigFont#1#2#3#4#5{%
  \reset@font\fontsize{#1}{#2pt}%
  \fontfamily{#3}\fontseries{#4}\fontshape{#5}%
  \selectfont}%
\fi\endgroup%
\begin{picture}(8283,6585)(-1425,-6630)
\put(-899,-3136){\makebox(0,0)[lb]{\smash{\SetFigFont{9}{10.8}{\rmdefault}{\mddefault}{\updefault}{\color[rgb]{0,0,0}\(\beta_g\)}%
}}}
\end{picture}

%% file: filtration.tex
\section{The Alexander Filtration}
\label{Sec:Filtration}

In this section, we consider the \OS Floer homology of the manifolds
\(N(n,1) \) when \(n\) is very large. To be specific, fix a \(
\spinc \) structure \(\spi \) on \(N\) and denote by \(\spi _k\) the
\( \spinc \) structures on \(N(n,1) \) which restrict to it. We study
 the groups \(\hfhat (N(n,1), \spi _k ) \) and the relations between
 them. It turns out that for \(n\) sufficiently large, the chain complexes
  \( \cfhat (N(n,1), \spi _k) \) are all generated by the same set of
  intersection points.   This fact enables us
to describe all of the \( \cfhat (\spi _k) \) in terms of a single
complex \( \cfs \), which we refer to as the stable complex. We show
that the
Alexander grading is a filtration on the stable complex, and we use
this fact to refine \( \cfs \) to a new complex \( \cfr \) --- the
reduced stable complex of \(N\). Our first main theorem is that \(
\cfr \) is actually an invariant of \(N\) and \(m\).  

Throughout this section, we continue to work with the basic framework
we set up in section~\ref{Ssec:GenAlex}. Since we made quite a few
choices in doing so, we pause to review them here.
First, we made homological choices: the classes \(m,l \in H^2 (\partial N)
\). These choices are extrinsic: the invariants we construct in this
section will depend on them. Second, we made  geometrical choices: the
Heegaard splitting \( (\Sigma, \alpha, \beta) \), the geometric
representative \(\emm\) of \(m\), and the labeling on the \(w^i\)'s. We
abbreviate all of these geometrical choices by the symbol \( \E \). 
Our
invariants do not depend \( \E\), but many of steps we take
in constructing them will. 

\subsection{\(\spinc \) structures}

We begin by establishing some  conventions about \( \spinc \)
structures on \(N\) and its Dehn fillings. Recall that
 a \( \spinc \) structure
on any  filling of \(N\) restricts to a \(\spinc \) structure on
\(N\). Conversely,  any \( \spinc \) structure
\(\spi \) on \( N \) extends to \(n\) different \( \spinc \)
structures \( \spi_1, \spi_2, \ldots, \spi_n \) on \(N(n,1) \) with
\( \spi_i - \spi _j = (i-j) PD([m]) \) in \(H^2 (N(n,1)) \). 

We would like to describe these processes of restriction and extension
in terms of Heegaard diagrams. To this end, we briefly recall the
description of \( \spinc \) structures given in \cite{Turaev} and 
\cite{OS1}. Let \(Y\)
be a three-manifold (possibly with boundary). There is a one-to-one
correspondence between \( \spinc \) structures on  \(Y\) and  homotopy
classes of non-vanishing vector fields on \(Y\), if \(Y\) has boundary, 
or \( Y- B^3 \), if \(Y\) is closed. To restrict a
\( \spinc \) structure from \(Y\) to a  codimension \(0\) submanifold \(Y_1\),
 we simply  restrict the corresponding vector fields. 

If \( (\Sigma, \alpha, \beta ) \) is a Heegaard
splitting for  a closed manifold \(Y \) and 
\(z \in \Sigma - \balpha - \bbeta \) is a
basepoint, any \( \bfy \in \TTT_{\alpha} \cap \TTT_{\beta} \) determines a
\( \spinc \) structure \( \spi _z (\bfy ) \).
 To construct this \( \spinc \) structure, we
equip \(Y\) with a Morse function \(f\) which gives the Heegaard
splitting. The point \( \bfy  \)
corresponds to a  \(g\)-tuple of flowlines from the index 2 critical points
to the index 1 critical points, such that every critical point is in
exactly one flow. In addition, the basepoint \(z\) determines a
flowline from the
index 3 critical point to the index 0 critical point. The vector field
 \( \nabla f \) is non-vanishing on the complement of a 
tubular neighborhood of these flowlines and extends to  a
non-vanishing vector field \({\bf v} \) on all of \(Y\).
 The \( \spinc \) structure determined by \({\bf v} \) is \( \spi _z
 (\bfy) \). 

Similarly, if \(Y_1 \) is a manifold with torus boundary and a
 Heegaard splitting \( (\Sigma, \alpha, \beta ) \),
 we can define a \( \spinc \) structure on \(Y_1\)
by choosing a basepoint \(z \in \Sigma - \balpha - \bbeta \) and
 points \( \haty \in T_{\hat{\alpha_i}}
\cap T_{\beta} \), \( y_i \in \alpha _i - \bbeta \). 
The construction is essentially the same as in the closed case,
 except now there is no
index 0 critical point, and the Morse function attains its minimum
along the boundary. In addition to the flowlines between the critical
 points determined by \( \haty \), we remove neighborhoods of 
the flowlines from the index 3 critical
point to the basepoint and from the index 2 critical point
corresponding to \( \alpha _i \) to \(y_i \). As in the closed case,
 it is easy to see that \( \nabla f \) extends from this complement to
 all of \(Y_1\), and that the homotopy class of the resulting vector
 field is independent of the extension. We denote the resulting \(
 \spinc \) structure by \( \spi _z (\haty, y_i)\). It is easy to
 see that \( \spi _z (\haty, y_i)\) is unchanged by 
 isotopies of \(z\) which avoid \( \balpha \)
 and \( \bbeta \).  In addition, 
if \( Y_1 \) is obtained by omitting a 2-handle
 from a Heegaard splitting of \( Y\), 
 the restriction of \( \spi _z (\bfy )  \) to \(Y_1\) 
is  \(\spi _z (\haty, y_i ) \).
Finally,  note that basepoints which are distinct in the splitting for
 \(Y \) may restrict to isotopic basepoints in the splitting for \(Y_1 \). 

\vskip0.05in
We now return to our particular three-manifold \(N\) and its
Heegaard splitting. We fix once and for all the following 
\begin{conv}
The basepoint \(z\) in the Heegaard splitting of \(N\) will always lie
in the region between \(w_1 \) and \(w_m \) which contains \(\emm\). When
we study Dehn fillings of \(N\), we will only consider basepoints
which restrict to this \(z\). 
\end{conv}

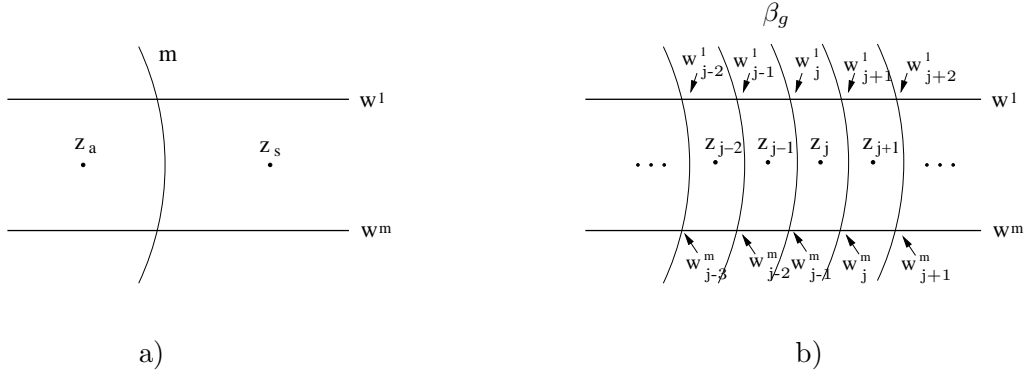
\begin{figure}
\input{figs/basepoints.pstex_t}
\caption{\label{Fig:Basepoints} Possible choices of basepoint on {\it a)}
\(\Nbar\) and {\it b)}  \(N(n,1)\).}
\end{figure}

For example, there are two isotopy classes of  points which restrict
 to \(z\)  in 
 our Heegaard splitting for \(\Nbar \); they are the points \(z_s \)
 and \( z_a \) shown in 
Figure~\ref{Fig:Basepoints}a. Since \( \spinc \)  structures on \(
\Nbar \)  are in 
\(1\)-\(1\) correspondence with \( \spinc \) structures on \(N\), we
 expect that 
\( z_s \) and \(z_a \) will induce the 
same \( \spinc \) structures on \(\Nbar \). This is indeed the case,
 since the longitude \(\ell\) is a curve which intersects \(\emm\) once
and misses all of the other \( \beta _i\). By Lemma 2.12 of
\cite{OS1}, \( \spi_{z_s} - \spi_{z_a} = PD ([l]) = 0. \)

We now suppose that we are given a \( \spinc \) structure \( \spi \) on
\(N\), and let \( \spibar \) be the \( \spinc \) structure on \(
\Nbar \) which restricts to \( \spi \).
There is a unique \( \ebar \) class \( \Ebar \)  on \( \Nbar \) with 
with \(\spi _z (\Ebar) = \spibar \). By Proposition~\ref{Prop:GoodEClass},
if \(n\) is sufficiently large, we can find a good \(\epsilon \)-class
 \(E \) on \( N(n,1) \) which lies above \(\Ebar \). 

\begin{lem}
The \(n\) \( \spinc \) structures on \(N(n,1) \) which restrict to \(\spi
\) are  \(\spi _{z_i} (E)\) for \( 1 \leq i \leq n \), where
the \(z_i\) are as shown in Figure~\ref{Fig:Basepoints}b. 
\end{lem}

\begin{proof}
Suppose \(\bfy = \{\haty, w^i_j\} \in E \), and that \( \overline{\bfy} = 
\{\haty, w^i\}\) is the corresponding point in \(\Ebar \). Then
clearly \( \spi_z (\overline{\bfy}) \) and \( \spi_{z_i} (\bfy ) \)
restrict to the same \( \spinc \) structure on \(N\). By the
definition of \(E_{\spi} \), this \( \spinc \) structure is \(\spi \).
Since \(\emm \) intersects \(\ell + n\emm \)
 exactly once and misses all the \(
\beta _i \), \( \spi _{z_i} - \spi _{z_{i+1}} = PD(m) \).
\end{proof}

 From
now on, we fix a  \( \spinc \) structure \( \spi \) on
\(N \) and a  good \( \epsilon\)-class \(E\) above
\( \Ebar \). 
We would like to label the \(\spi _{z_i} (E)\) in some way which does
not depend on which \(E  \) we picked. The easiest way to do this is
to relabel the \(w^i_j\)'s and \(z_j\)'s according to our choice of
\(E\).  Recall that the different possible choices of \(E\) all have the same
Alexander grading \(A\). We fix an
integer-valued lift of \(A\)  ---  if \(N\) is a  knot complement in
\(S^3\), we
use the canonical one. Once we have chosen \(E\), we 
add some constant to  the lower index of  all the \(w^i_j\)'s so that
the Alexander grading
of the point \({\bf y} = \{\haty, w^i_j\} \in E \) is  \(j\). We
label the \(z_j\)'s so that \(z_j \) is in the region below and to the
right of \(w^1_j \), as shown in Figure~\ref{Fig:Basepoints}. 
With these conventions, it is not difficult to see that \( \spi_j 
= \spi_{z_j}(E) \) is independent of our choice of \(E\).
\vskip0.1in
%\noindent{\bf Remark:} So far we have not specified how we tell which
%side of \(m\) \(z_s \) is on. To do this in a consistant fashion, we
%choose a generator of \( H^2(N(n,1) ) \) and require that \(z_a - z_s
%\) be a positive multiple of this generator.
%\vskip0.1in
\noindent{\bf Example:} If \(N = \Ko \) is a knot exterior, there is only one
\(w^i = w \), and the signed number of points above \(w_j \) is just the
coefficient of \(t^j \) in \( \Delta_K (t) \). In \cite{2bridge}, it
is shown that if \(K\) is a two-bridge knot with Heegaard splitting
coming from the bridge presentation, all the points above \(w_j \)
have the same sign, so their number is  determined by the
Alexander polynomial. This is an unusual state of affairs:
almost any Heegaard splitting coming from an \(n\)-bridge diagram with
\(n \geq 3 \) will have more points in \(E \) than the minimum number
dictated by \( \Delta_K (t) \).

\subsection{The stable complex}
We  now begin our study of  \( \cfp (\TTT _\alpha \cap
 \TTT _{\beta }, \spi _k) \). We tend to think of these complexes as
 being generated by a single good \( \epsilon \)-class  with a varying
 basepoint \(z_k \), so we use a different notation from
 that of \cite{OS1}. For the moment, we  write
\( \cfp (\E,n, E, \spi_k) \) to denote the complex which computes \(
 \hfp (N(n,1), \spi _k ) \) from our standard Heegaard splitting,
 using the good \(\epsilon\)-class \(E\) and the basepoint \(z_k\).  
This notation is a bit  cumbersome, but we will soon show that it can
 be shortened. 

We have already seen that to each point \( \bfy = \{ \haty, w^i_j \} \in E \)
there is a corresponding point   \(\overline{\bfy} = \{ \haty, w^i \} 
\in E_{\spi}\). In fact, we have  

\begin{lem}
\label{Lem:SpiralCrush}
For \(k\gg 0 \), there is an isomorphism \( \cfp (\E,n,E, \spi_k)
\cong \cfp (\E, \Nbar,z_s, \spibar ) \).
\end{lem}

\begin{proof}
Doing \(n\) surgery on the knot \(K\) induces a cobordism from \(
\Nbar \) to \(N(n,1) \). Let \(f\colon \cfp (\E, \Nbar,z_s, \spibar ) \to 
\cfp (\E,n,E, \spi_k) \) be the map induced by this
cobordism. 
In the relevant Heegaard diagram, 
there is a unique triangle \( \psi _{\bfy}\) joining \( \bfy \) to \(
\overline{\bfy}\) supported in the spiral region: all other triangles
involving \( \bfy \) have domains which include regions outside the
spiral. It is easy to see that \( \mu(\phi) = 1 \) and \( \# \MMM
(\phi) = \pm 1\), and when
\( k \gg 0 \), \(n_{z_k}(\psi_\bfy) = n_{z_s}(\psi_\bfy) = 0.\) If we
make the area of the spiral region very small compared to the area of
the other components of  \( \Sigma- \alpha - \beta \),
 the area function induces a
filtration with respect to which \(f(\overline{\bfy} ) = \bfy + \text
{lower order terms}\) ({\it cf.} section 8 of \cite{OS2}.) It follows
that \(f\) is an isomorphism of chain complexes. (In fact, we conjecture that when the spiral is sufficiently tight, the
map \(f\) is precisely given by \(f(\overline{\bfy}) = \bfy \).)
\end{proof}
It is easy to see that  for \( k \ll 0 \) there is an  analogous result
 relating \(\cfp (\E, \Nbar,z_a, \spibar ) \) with
\( \cfp (\E,n,E, \spi_k) \).

\begin{cor}
If \(E' \) is another good \(\epsilon \)-class, then for \(k \gg 0 \),
there is an
isomorphism between \( \cfp (\E,n, E', \spi_k) \) and
\( \cfp (\E,n, E, \spi_k) \). Moreover, if \( n \) and \(n'\) are both
large enough for good \( \epsilon\)-classes to exist, there is an
 isomorphism \( \cfp (\E,n', E, \spi_k)  \cong
 \cfp (\E,n, E, \spi_k) \).
\end{cor}

\begin{proof}
Compare both complexes with \(\cfp (\E, \Nbar,z_s, \spibar ) \).
\end{proof}

In light of these facts, we make the following 

\begin{defn}
The stable complex of \(\E \), written \( \cfs (\E) \), is the complex
\begin{equation*}
\cfhat (\E, \Nbar,z_s, \spibar )\cong \cfhat (\E,n,E, \spi_k) 
\end{equation*} 
for \(k \gg 0 \). The antistable
 complex \( \cfa (\E) \) is \(\cfhat  (\E, \spi_k ) \) for \(k \ll 0
 \). The complexes \(\cfpms  (\E) \), and \(\cfpma (\E) \)
 are defined analogously. We refer to the 
 \( \spi _k \) with \( k \gg 0  \) (resp.
 \( k \ll 0 \)) collectively  as the stable \( \spinc
 \) structure \( \spi _s \) (resp. the antistable \( \spinc \)
 structure \( \spi _a \).) 
\end{defn}

\noindent {\bf Remarks}: Although \(\hfp (\Nbar, \spi) \)  is the
homology of both \(\cfps (E) \) and \(\cfpa
(E) \), the two are usually not
isomorphic as complexes. Note that at the level of homology  
 Lemma~\ref{Lem:SpiralCrush} follows from the exact triangle of  \cite{OS2},
 which gives a long exact sequence 
\begin{equation*}
\begin{CD}
 @>>> \hfp (\Nbar, \spi) @>>> \bigoplus _{i \equiv k \ (n)} \hfp
 (N(0,1), \spi _i ) @>>> \hfp (N(n,1), \spi _k) @>>>
\end{CD}
\end{equation*}
 \( \hfp (N(0,1))  \) is supported in a finite number of \( \spinc \)
 structures,  so when \(n\)
is very large, most of the \(\hfp (N(n,1), \spi _k) \) will be
isomorphic to \(\hfp (\Nbar, \spi) \). On the other hand, we have
already seen that most of the \(\cfp (\E, N(1,n), \spi _k) \) are
isomorphic either to either \( \cfps (\E) \) or \( \cfpa (\E) \). 
This is the first instance of a principle which we will see more of
in the future: {\it   the behavior  of the exact triangle is
 realized by the complexes} \( \cfp (\E, \spi _k) \).

\subsection{The generators of \(  \cfp (E, \spi _k) \)}
 
 We now consider the complexes \(  \cfp (\E,n, E, \spi_k)\) for
 arbitrary values of \(k\). In analogy with the notation for \(\cfs\),
we drop the dependence on
 \(E\) and \(n\) from the notation, and simply denote  these
 complexes by \( \cfp (\E, \spi_k) \).
 (We will show  this is justified in the next section.)
These complexes all have the same generating set \(E\). Indeed, 
the numbers \( n_{z_k}(\phi) \) are the only things which
distinguish the  \( \cfp (\E, \spi_k) \)  from each other.
 Fortunately, these numbers are  easily expressed
in terms of the Alexander grading on \(E\). The result is summarized
in the following handy

\begin{lem}
\label{Lem:nzChange}
Suppose \( \phi \in \pi_2(\bfy, \bfz) \). Then
\begin{align*}
n_{z_{k+1}}(\phi) - n_{z_k}(\phi) = 
\begin{cases}
1 & \text{if} \  A(\bfz) >k \geq A(\bfy) \\
-1 & \text{if} \ A(\bfy) >k \geq A(\bfz) \\
0 & \text{otherwise}
\end{cases}
\end{align*}
\end{lem}

\begin{proof}
We refer  to Figure~\ref{Fig:Basepoints}. Recall that 
the  \( \beta _g \) component of \( \partial \DDD (\phi ) \) is
oriented to point from \( \bfz \) to \( \bfy \). In the first case, it
traverses the segment separating \(z_{k+1} \) from \(z_k \) 
 once with a positive (upward)
orientation, in the second, once with a negative
orientation, and in the last it does not pass over it at all.
\end{proof}

We can now use the stable complex to describe \(  \cfp (\E, \spi _k)
\). Indeed,  the
complexes \( \cfi (\E, \spi _k ) \) are the same for every choice of
\(k\), so we can realize all of the \( \cfp (\E, \spi _k)
\) as quotients of  \( \cfis (\E) \). To be precise, we label  the
generators \(\{[\bfy, j] 
\ts | \ts \bfy \in E, \ts j \in \Z \} \) of 
 \( \cfis (\E) \) according to the usual convention, so that 
 \(\{[\bfy, j]  \ts | \ts \bfy \in E, \ts j\geq 0\} \) are the generators
of \( \cfps (\E) \). Then we have

\begin{lem}
\label{Lem:CFPGen}
\(\cfp (\E, \spi _k) \) is the quotient complex of \( \cfis (\E) \)
with generators 
\begin{equation*}
 \{ [\bfy, j] \ts |  \ts \bfy \in E, \ts  j \geq \min (k- A(\bfy), 0) \}.
\end{equation*}
\end{lem}

\begin{proof} 
We temporarily denote the generators of \( \cfp (\E, \spi _k ) \) by
\(\{[\bfy, j]_k  \ts | \ts \bfy \in E, \ j\geq 0\} \). Inside \( \cfis \),
 \( [\bfy, j]_k \) is identified with \( [\bfy, j - n_{\bfy}] \) for
 some number \(n_{\bfy} \). To find this number,
  we choose  \(\bfz \in E \) with minimal Alexander
  grading. Since \(\cfis \) is translation invariant, we may as well
  identify \([\bfz, 0]_{k} \) with \([\bfz,0] \).
 Any other \( \bfy \in E \) may be connected to \(\bfz\) by a
  domain \( \phi \in \pi_2 (\bfy,\bfz) \) whose \( \beta _g \)
  component runs  from \( w(\bfz) \) to \( w(\bfy) \) inside the
  spiral. Then we have 
\begin{align*}
\gr_{\spi _k}(\bfy) - \gr _{\spi_k}(\bfz) & = \mu (\phi) -
2n_{z_k}(\phi) \\
& = \mu (\phi ) - 2 n_{z_M}(\phi) - 2(n_{z_k}(\phi) - n_{z_M}(\phi)) \\
& = \gr_{\spi_s}(\bfy) - \gr_{\spi _s} (\bfz) - 2(n_{z_k}(\phi)
 - n_{z_M}(\phi)).
\end{align*}
It follows that \(n_{\bfy} = n_{z_k}(\phi) - n_{z_M}(\phi) \).
 Applying Lemma~\ref{Lem:nzChange}, we see that 
\begin{equation*}
 n_{z_k} (\phi ) - n_{z_M}(\phi) = \begin{cases}
 0 & \text{if} \ k\geq A(\bfy) \\
 A(\bfy)-k & \text{if} \ k \leq A(\bfy).
\end{cases}
\end{equation*}
\end{proof}

\begin{cor}
\label{Cor:Seq1}
There is a short exact sequence of chain complexes 
\begin{equation*}
\begin{CD}
 0 @>>> C_{\spi_k}(\E) @>>> \cfp (\E,\spi_k) @>>> \cfps (\E) @>>> 0 
\end{CD}
\end{equation*}
where \( C_{\spi_k}(\E) \) is the subquotient of \( \cfis (\E) \) with
generators 
\begin{equation*}
 \{[\bfy, j] \ | \  \min (k - A(\bfy),0) \leq j < 0 \}.
\end{equation*}
\end{cor}

\begin{proof}
We have inclusions of complexes \( \cfm (\E,\spi_k) \subset \cfms (\E) 
\subset \cfi (E) \). The given sequence is just the short
exact sequence of quotients.  
\end{proof}

As we will describe in section~\ref{Sec:ExactTriangle}, the associated
long exact sequence serves as a model for the exact triangle, with \(
H_*(C_{\spi_k}) \) playing the role of \( \hfp (N(0,1),\spi _k ) \).

\subsection{The Alexander Filtration}

The exact sequence of Corollary~\ref{Cor:Seq1} is our first example of
the  rich filtration structure  on \( \cfp (\E, \spi _k) \). We
now investigate this structure more systematically. 
Our first step is to extend the notion of the Alexander
grading to \( \cfis (\E) \):

\begin{defn}
The Alexander grading of \( [\bfy,j] \in \cfis (\E) \) is 
 \( A([\bfy,j]) = A(\bfy) + 2j \).
\end{defn}

Viewing  \(\cfhat (\E,\spi _k) \) as a subset of
\(\cfis (\E) \) allows us to restrict the Alexander grading to the
generators of
 \( \cfhat (\E, \spi _k) \). We denote the resulting affine grading on
 \(E\) by \( A_{\spi _k} \).  From
Lemma~\ref{Lem:CFPGen}, we see that 
\begin{equation}
\label{Eq:AlexGr}
A_{\spi_k} (\bfy) = \begin{cases}
A(\bfy) & \text{if} \ k \geq A(\bfy) \\
2k -  A(\bfy) & \text{if} \ k \leq A(\bfy)
\end{cases}
\end{equation}
In particular, \( A_{\spi _a} (\bfy) \) is equivalent to  \( - A(\bfy)
\) as an affine grading.

It is often helpful to think of \( \cfis (\E) \)  in terms of a diagram
like that shown in   Figure~\ref{Fig:Filtrations}. We attach each generator
\([\bfy, j] \) of \(\cfis (\E) \) to the dot with coordinates \( (x,y)
= (-A(\bfy), A([\bfy, j])) \). All of the points attached to a given dot
behave the same way with respect to the Alexander grading. 

It is well-known that \( \cfis (\E) \) has a filtration
\begin{equation*}
\ldots \subset C^s_{-1} \subset C^s_0 \subset C^s_1
\subset \ldots 
\end{equation*}
with \( C^s_j =  \{[\bfy, i] \ts | \ts i \leq j \}  \cong \hfms (\E ) \).
This filtration is indicated by the solid lines in 
Figure~\ref{Fig:Filtrations}. 
 This basic fact becomes
more interesting when we realize that the filtration
 \begin{equation*}
 \ldots \subset C^a_{-1} \subset C^a_0 \subset C^a_1
\subset \ldots 
\end{equation*}
  (shown by the dashed lines)  obtained by viewing \( \cfi (\E) \) 
 as \(\cfia (\E) \) is very different from \(C^s_j\).
Thus  \( \cfi \) is actually equipped with a double filtration. 
The following result is obvious from the figure:

\begin{figure}
\includegraphics{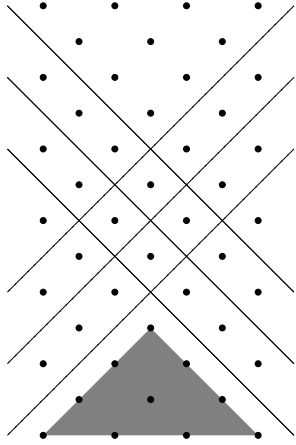}
\caption{\label{Fig:Filtrations} The complex \( \cfi
  (E) \). Differentials coming from the top vertex of the shaded triangle
  must land inside the shaded region.  }
\end{figure}

\begin{prop}
\label{Prop:AlexFilt}
The Alexander grading is a filtration on \( \cfis (\E)\); 
  i.e. if there is a differential from \([\bfy,i] \) to \( [\bfz,j] \) 
 in \( \cfis (\E) \), \(A([\bfy,i] ) \geq A([ \bfz,j]) \). 
\end{prop}

We refer to this filtration as the {\it Alexander filtration}. Note
that we do not gain anything more by considering \( \spinc \)
structures other than \(\spi_s \) and \( \spi _a\) --- the presence of
these two filtrations implies the presence of the filtration
induced by any \( \spi _k \). 

The fact that the Alexander grading is a
filtration can also be derived from the following useful formula for a disk
\( \phi \in \pi_2(\bfy, \bfz) \):
\begin{equation*}
n_{z_a}(\phi) - n_{z_s}(\phi) = A(\bfy) - A(\bfz).
\end{equation*}
\vskip0.05in
\noindent{\bf Remark:}
It is easy to check that the map \(f\) of Lemma~\ref{Lem:SpiralCrush}
respects the Alexander filtration. This implies that the complexes 
 \(  \cfp (\E,n', E, \spi_k)\) are isomorphic for all sufficiently
 large \(n'\) and good \(\epsilon\)--classes, thus justifying our
 omission of \(n'\) and \(E\) from the notation.

\subsection{The reduced stable complex}
As a special case of  Proposition~\ref{Prop:AlexFilt}, the
Alexander grading   induces a filtration on \( \cfs
(\E) \). From this filtration we get a spectral sequence 
\( (\scx_i(\E), d_i) \) which converges to \( \hfhat (N(n,1), \spi_s) \cong
\hfhat (\Nbar, \spibar) \). 
 This  sequence is an object of considerable interest: many of
 the results described below could also be phased in terms of it.
 For some
purposes, however, it is more convenient to have an actual { complex},
rather than a spectral sequence. 
For this reason, we introduce the
following construction, which we refer to as {\it reduction}. Suppose
\( C \) is a filtered complex with filtration \(
 C_1 \subset C_2 \subset \ldots C_m \), and let 
\( C^i = C_{i+1}/C_{i} \) be the filtered quotients. 
The homology groups \(H_*(C^i) \) are the \(E_2 \)
 terms of the spectral sequence associated to the filtration. 
\begin{lem}
\label{Lem:Reduction}
Let \(C\) be a filtered  complex over a field.
 Then up to isomorphism there is a unique
 filtered complex  \( C' \) with the following
properties:
\begin{enumerate}
\item \(C' \) is chain homotopy equivalent to \(C\). 
\item \( C'^i \cong H_*(C^i) \)
\item The spectral sequence of the filtration on \( C' \) has trivial
  first differential. Its higher terms are the same as the higher
  terms of the spectral sequence of the filtration on \(C\). 
\end{enumerate}
We refer to \(C' \) as the reduction of \(C\). 
\end{lem}
 \noindent The proof will be given in
 section~\ref{Subsec:Reduction}. 
%The requirement that \(C\) have filed
% coefficients is a technical one, which we hope can be removed by a
% more careful construction of the reduction. For the moment, however,
% we will stick with field coefficients 
%(over \(\Q\) unless otherwise specified.)
\begin{defn}
The reduced stable complex \( \cfr (N,m, \spi) \) is the reduction of
 \(\cfs (\E) \). 
\end{defn}

The intuitive picture behind this construction is as follows. Recall that each
dot in Figure~\ref{Fig:Filtrations} represents a whole set of
generators with some particular Alexander grading. On \( \cfs (\E) \), we
split \(d\) into two parts: \(d = d' + d'' \), where \(d' \)
preserves the Alexander grading and \( d'' \) strictly reduces
it. Thus \(d' \) involves those differentials which ``stay within'' a
given dot, while \(d'' \) contains those differentials which go from
one dot to another. The condition \(d^2 = 0 \) 
implies that \((d')^2 = 0 \) as well, so each dot is itself a little
chain complex. The idea  is to find a new complex
 chain homotopy equivalent to \( \cfs (\E) \) in which we have
replaced the chain complex inside each dot by its homology.

Our motivation for considering \(\cfr \) instead of \( \cfs \) is that 
the new complex is actually a topological invariant:

\begin{tm}
\label{Thm1}
The filtered complex \( \cfr (N,m,\spi) \) is an invariant of the triple
\((N,m, \spi) \). 
\end{tm}

The proof will be given in 
section~\ref{Sec:Thm1}. For the moment, we remark that the theorem is a
specific instance of the general principle
that the higher terms of spectral sequences 
tend to be  topological invariants.
 (A more familiar example of
 this phenomenon is provided by the Leray-Serre sequence
of a fibration, in which the \(E_1 \) term depends on the
triangulations of the base and fibre, but the \(E_2\) and higher terms are
canonically defined in terms of the { homology} of these spaces.) 
\vskip0.1in

We record some elementary facts about the reduced stable complex below:

\begin{prop}
\label{Prop:StabCxFacts}
The reduced stable complex \( \cfr (N,m, \spi) \) has the
following properties:
\begin{enumerate}
\item Its homology is isomorphic to \( \hfhat (\Nbar, \spi) \).
\item The filtered subquotient \( \cfr ^{(j)} (N, m, \spi) \) is
  isomorphic to \( \scx ^j_2 (\E) \). Its Euler characteristic
 is the coefficient
  of \(t^j\) in the Alexander polynomial \( \Delta  (N, m, \spi ) \). 
%\item The spectral sequence induced by the Alexander filtration on \(
%  \cfr (N,m, \spi) \) has trivial \(d_1 \) and 
%is \( (\scx _i, d_i) \) for \( i \geq 2 \). 
\item \(\cfr (-N, m, \spi ) \) is the dual complex to \( \cfr (N, m, \spi)
\). 
\end{enumerate}
\end{prop}

\begin{proof}
Parts {1)} and 2) follow from the definition of reduction, combined
   with  Lemma~\ref{Lem:SpiralCrush}
 and  Corollary~\ref{Cor:AlexP}, respectively.

To prove part 3), recall that if \( (\Sigma, \alpha, \beta ) \)
is a Heegaard splitting for \( N\), 
\((-\Sigma , \alpha, \beta)\) is  a splitting for \(-N\). 
There is an obvious bijection between the elements of good
\(\epsilon \)-classes  \(E_{\pm N}\)  for these splittings, and the
Alexander grading on \(E_{-N} \) is \(-1\) times the Alexander grading on
\(E_N \). Combining this fact with the isomorphisms
\begin{equation*}
\cfs (\E, \spi ) \cong \cfhat ( \E, \Nbar,z_s, \spibar) \cong 
(\cfhat ( \E, -\Nbar, z_s, \spibar))^* \cong (\cfs (-\E, \spi ))^*
\end{equation*} 
provided by Lemma~\ref{Lem:SpiralCrush} and Proposition 7.3 of
\cite{OS1}, we see that \( \cfs (\E, \spi ) \) is dual to 
\( \cfs (-\E, \spi ) \) as a filtered complex. The result then follows
from the fact that the dual of the reduction is the reduction of the
%this is only true over Q!
dual, which will be obvious from the construction of
section~\ref{Subsec:Reduction}.
\end{proof}

The first part of the proposition is particularly useful when
\((N,m)\) is a knot complement in \(S^3\), so that  \(\hfhat (\Nbar)
\cong \hfhat(S^3) \cong \Z \). 

\begin{defn}
If \(K \) is a knot in \(S^3\), the level of the generator of \(K \), 
written \(s(K) \), is the Alexander grading of the surviving copy of
\(\Z\) in \(\scx _{\infty} ( \Ko, m) \). 
\end{defn}

By Theorem~\ref{Thm1}, \(s(K) \) is an invariant of \(K\). Part 4)
 of the proposition implies that \(s(-K) = -s(K) \). 
We will see in section~\ref{Sec:ExactTriangle} that 
 this invariant  encodes useful information about the exact
triangle for  surgery on \(K\). 

\vskip0.1in
It is not difficult to see that everything we have  done  using
the stable complex works just as well when applied 
 to the antistable complex: 
 the Alexander grading \( A_{\spi _a } = -A \) gives a filtration
on \( \cfa (\E, \spi) \),
there is a   spectral sequence \( (\scx^a _i, d^a_i) \) derived from 
it,
and the reduced  complex \( \cfr ^a \) is a topological
invariant. How is this new object related to \( \cfr \)? Since the
first differential in each spectral sequence preserves the Alexander
grading, it is easy to see that  \( (\scx ^a_2)^j \cong \scx _2 ^{-j}
\). Thus the two complexes have the same filtered subquotients (but
in opposite orders.) Somewhat less obviously, we have 

\begin{prop}
\label{Prop:Symmetry}
\begin{equation*}
\cfr ^a (N, m , \spi ) \cong \cfr (N, m, -\spi  )
\end{equation*}
as filtered complexes.
\end{prop}

\begin{proof}
The basic idea is to combine the identification \( \cfs (\E, \spi)
\cong \cfhat (\E, \Nbar, z_s, \spibar) \) with the conjugation
symmetry \( \hfhat (\Nbar, \spibar) \cong \hfhat (\Nbar, -\spibar)
\). The latter isomorphism is realized by simultaneously reversing the
orientation of \( \Sigma \) and switching the roles of \( \alpha \)
and \( \beta \) in a Heegaard splitting. 

With this in mind, we consider the Heegaard splitting
\( - \E = (-{\Sigma}, \beta, \alpha) \) of \(N(n,1) \). There
is an obvious identification between 
\( \cfhat (\E, E, z_a) \) and \( \cfhat (- \E, E, z_a) \). We will
show that the latter complex is actually of the form \( \cfs ( \E',
-\spi) \) for some new choice of Heegaard data \(\E ' \). Reducing
both sides and applying Theorem~\ref{Thm1} gives  the desired isomorphism.

By Lemma~\ref{Lem:Meridian} we can assume 
 that our original \( m\) has a single geometric intersection with
 \( \alpha _g \). In this case,  \( \alpha _g \) and \(\beta _g \) 
are mirror-symmetric in the spiral region. This is not obvious from
 our usual way of drawing the spiral, but if we view the annulus
 as a cylinder and twist the ends, we can put \( \alpha _g \) and \(
 \beta _g \) into the form of Figure~\ref{Fig:Cylinderspiral},
 in which the symmetry is overt. 

\begin{figure}
\includegraphics{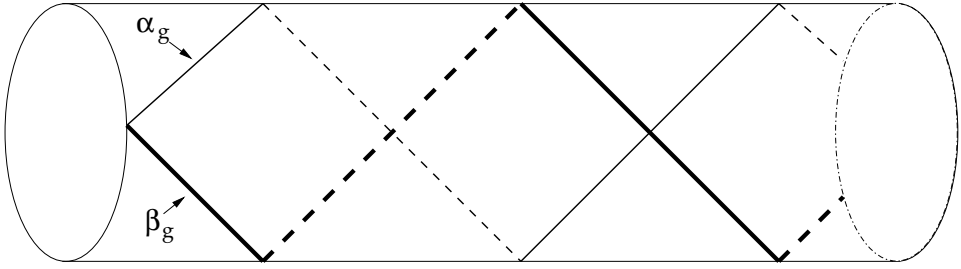}
\caption{\label{Fig:Cylinderspiral} The curves \( \alpha _g \) and \(
  \beta _g \) in the spiral region.}
\end{figure}

Thus our new Heegaard splitting
\( -\E \) has the form we used to define the stable complex: the new
\( \beta _g \) (which used to be \(\alpha_g \)) has been twisted
around \(m\) many times in a positive sense. It is not difficult to
see that \(E\) is a good \( \epsilon \)-class for this new splitting,
and since we reversed the orientation of \( \Sigma \), the basepoint
\(z_a \) is actually on the stable side of the new complex.  The only
thing remaining to check is that the manifold that we get by omitting
\( \alpha _g \) from \( -\E \) is still \(N\). For this, we use

\begin{lem}
Suppose \( (\Sigma, \alpha, \beta) \) is a Heegaard splitting of \(Y\),
and \(m \subset Y \) is represented by a curve on \( \Sigma \) which
intersects \( \beta _g \) geometrically once and misses all of the other \(
\beta _i\)'s. Then \( (\Sigma, \alpha, \beta_1, \ldots, \beta _{g-1}) \)
is a Heegaard splitting of \(N-m \). 
\end{lem}

\begin{proof}
If we push \(m\) off of \(\Sigma \), we see that it punches a single hole in
the handle attached along \( \beta _g \) and misses  the other \(
\beta \) handles.  To get the desired splitting, we retract
 what is left of the \(\beta _g \) handle
back to \(\Sigma \). 
\end{proof}

Applying the lemma to the two splittings \( \E \) and \(- \E \) of 
\(N(n,1) \), we see that omitting \( \alpha _g \) from \( -\E \) 
 gives the same manifold as omitting \( \beta _g \) from \(\E \).
Since latter manifold is \(N\), the claim is proved.
\end{proof}

Thus we do not get any new invariants by considering the antistable
complex. What we do get, however, is an extension of the usual
conjugation symmetry on \( \hfhat \) 
 to the stable complex. For example, when \(N =
\Ko \) is a knot complement, the proposition implies that 
\( \scx _2 ^j(K) \cong \scx_2^{-j} (K) \). 
 This fact can very helpful in
practice, since one of the \(\scx _2^{\pm j} \) may be easier to
compute than the other. 

To make full use of the symmetry between stable and antistable
complexes, it is convenient to introduce yet another affine grading on \(E\):

\begin{defn}
For \( \bfy \in E \), we set
%\begin{equation*}
\(\mutid (\bfy ) =  \gr_{\spi_k}(\bfy) - A_{\spi _k}(\bfy )\).
%\end{equation*}
\end{defn}

The important point is that \( \mutid \) does not depend on which \(
\spi _k \) we use. This is obvious if we think of \( \cfhat (\E,
\spi_k) \) as a subset of \( \cfis (\E) \): the difference
 \( \gr ([\bfy,j]) - A([\bfy,j]) \) is clearly independent of \(j\).
This fact gives us an easy way to compute all the \( \gr _{\spi
  _k}\)'s from any one of them.

Since it depends only on the Alexander and homological gradings, the
definition of \( \tilde{\mu} \) makes sense for for a generator of \(
\cfr \) as well. Suppose that \( N = \Ko \) is a knot complement, and
set
\begin{equation*}
C(K,j,l) = \{ \bfy \in \cfr (K) \ts | \ts A(\bfy) = j, \ts \tilde{\mu}(\bfy)
= l \} 
\end{equation*}
The symmetry between \( \scx _2 ^j \) and \(\scx _2 ^{-j} \) implies
that \( C(K,j,l) \cong C(K,-j,l) \). Thus from the point of view of
\( \gr _{\spi_s}\), \( \cfr  (K) \) is composed of 
symmetrical pieces on which \( \tilde{\mu} \) is constant, but these
pieces are put together in such a way that \(\gr _{\spi_s} \) is not
symmetric on \( \cfr  (K) \). (See the discussion of the \((3,5)\)
torus knot below for an example of such a decomposition.)

We conclude our discussion of the reduced stable complex by describing
\( \cfr (K) \) for two basic types of knots.
\vskip0.1in
\noindent{\bf Example 1: Two-bridge knots} 
If \(K\) is a two-bridge knot with the two-bridge Heegaard
splitting, \( \gr_{\spi_s} = A \). Since  every differential
 reduces the Alexander grading by one, \( \cfr (K) \cong \cfs (\E) \), and
 the rank of  \( \scx ^2_j (K) \) is the absolute value of the
 coefficient of \(t^j \) in \( \Delta _K (t) \).
 In addition, \( s(K) = \sigma (K) / 2
 \), where \( \sigma \) is the ordinary knot signature.
Thus the isomorphism class of \( \cfr (K) \) is completely determined
by classical knot invariants. 
For two-bridge knots,  the symmetry of \( \cfr (K)\) is explicitly realized by
 a natural involution of the Heegaard splitting, as described in
 \cite{2bridge}. 
\vskip0.1in
\noindent {\bf Example 2: Torus knots} If \(K = T(p,p+1)\)
 is the  right-handed
\((p,p+1)\) torus knot, the
coefficients of \(\Delta _K (t) \) are always \(\pm 1 \) or \(0 \). We
show in the appendix  
that the rank of \(\scx (K) _2 ^j \) is the absolute value of
the corresponding  coefficient of \(\Delta _K (t) \). Let 
  \(\bfy _1, \bfy _2, \ldots, \bfy _m \) be the generators of \( \cfr
  (K) \) arranged in 
order of decreasing Alexander grading.  Then there is a differential from \(
\bfy _{2i} \) to \( \bfy _{2i+1} \). The presence of these differentials
enables us to compute \( \gr _{\spi _s} \). It turns out that 
for \(i>j\), \( \gr _{\spi_s} (\bfy _i) >\gr _{\spi_s} (\bfy _i)\), so
 the differentials described above are the only 
differentials in \( \cfr (K) \). In turn, this implies that \( \hfs
(K) \)  is generated by \( \bfy _1 \), so  \(s(T(p,p+1)) =
g(T(p,p+1)) = p(p-1)/2 \).
\vskip0.05in
We illustrate how to how to find \( \gr _{\spi
 _s} \) in the case of 
  the \((3,5) \) torus knot. The Alexander grading on \( \cfr (K) \)
 is shown in Figure~\ref{Fig:(3,5)TorusCx}a. 
There is a differential from \( \bfy _{2i} \) to \( \bfy _{2i+1} \),
 so 
\begin{equation*}
 \gr _s ( \bfy _{2i}) -  \gr _s ( \bfy _{2i+1}) = 1 .
\end{equation*}
 On the other
hand, the symmetry between stable and antistable complexes shows that
in the antistable complex, there must be a differential from \( \bfy
_{2i} \) to \( \bfy _{2i-1} \), so 
\begin{equation*}
 \gr _a ( \bfy _{2i}) -  \gr _a ( \bfy _{2i-1}) = 1 .
\end{equation*}
 Then we compute, for example, that 
\(\mutid (\bfy_6, \bfy_5)  = -1 \), so 
\(\gr _s (\bfy_6, \bfy _5)  = -3 \). The other
gradings can be found by the same method.

\begin{figure}
\includegraphics{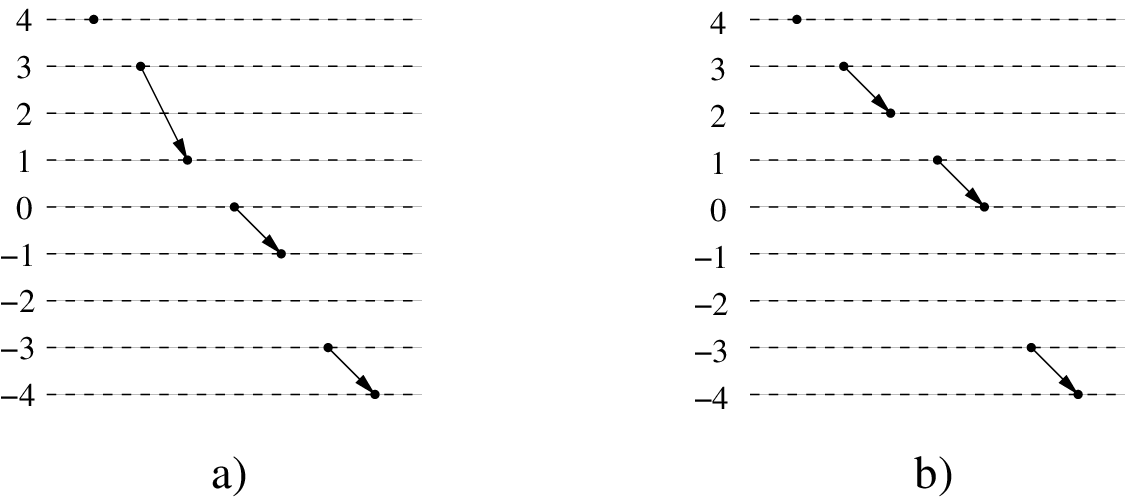}
\caption{\label{Fig:(3,5)TorusCx}} The reduced stable complex of the
\((3,5) \) torus knot, showing the {\it (a)} Alexander  and {\it (b)}
homological  gradings of the generators.
\end{figure}

\subsection{Reduction for \( \cfis (\E) \)}
\label{SubSec:ReducedCFI}
There is no reason that the process of reduction should be limited to
 the stable complex.  In fact, there are reduced versions of
the all the various  \OS Floer complexes. We briefly 
describe these  other reduced complexes and
their relation to the reduced stable complex. 
Since the complex  \( \cfis (\E, \spi) \) contains all of the others,  it
 is a natural place to start. We filter \( \cfis (\E, \spi) \) by the
Alexander grading and denote the associated reduced complex by
\( \cfir (\E, \spi) \). (This is not a finite filtration, but it is
finitely supported in each degree, so it is still possible to take the
reduction.)

\begin{lem}
As a group, 
\(\cfir (\E) \cong \bigoplus_{i,j} \ts [\scx_2^i(\E, \spi), j] \).
\end{lem}

\begin{proof}
Let \( C^l \) be the subquotient of \( \cfis (\E, \spi) \) generated by the
elements of Alexander grading \(l\). It is easy to see from
Figure~\ref{Fig:Filtrations} that as a complex, 
\begin{equation*}
C^l = \oplus_{i+2j = l} \ts [\scx_1^i(\E, \spi), j].
\end{equation*}
\end{proof}

 \( \cfir (\E, \spi) \) retains the same
double filtration structure as \(\cfis (\E, \spi) \), so it has
subcomplexes \( \cfmr (\E, \spi _k) \)
and quotient complexes \( \cfpr (\E, \spi _k) \)
specified by the same sets of dots
as \( \cfm (\E, \spi _k) \) and \( \cfp (\E, \spi _k) \). We would
like to know that these complexes are the same as those obtained by
reducing \( \cfpm (\E, \spi _k) \). This is implied by the following
lemma, which will be proved in section~\ref{Subsec:Reduction}:

\begin{lem}
\label{Lem:ReducedSequence}
Suppose \(B\) is a filtered complex, and that there is a short exact
sequence
\begin{equation*}
\begin{CD}
0 @>>> A @>>> B @>>> C @>>> 0 
\end{CD}
\end{equation*}
which respects the filtration, in the sense that \(B^i = A^i \oplus
C^i \) (as complexes). Then there is a commutative diagram
\begin{equation*}
\begin{CD}
0 @>>> A @>>> B @>>> C @>>> 0 \\
@VVV @VVV @VVV @VVV @VVV \\
0 @>>> A^r @>>> B^r @>>> C^r @>>> 0
\end{CD}
\end{equation*}
in which the vertical arrows are homotopy equivalences. 
\end{lem}

Thus as far as the  homology is concerned, there is no difference
between using the reduced complexes and their unreduced counterparts.
Similarly, we can view  the reduced complex \( \cfr (\E, \spi _k)
\) as a subcomplex of \( \cfpr (\E, \spi _k) \). 

In analogy with Theorem 1, we expect that all of these reduced
complexes should also be
topologically invariant. The next proposition shows that this is true for
\( \cfr (\E, \spi _k) \). For the others, some more careful
accounting in the proof of Theorem 1 might provide a proof. (This is
done in \cite{OS7}.) Since we
do not have any immediate need for this invariance, however,
 we will not pursue the matter here. 

We can use our knowledge of \( \cfr (N, m, \spi ) \) and \( \cfa (N,
m, \spi) \) to understand some of the differentials in the other reduced
complexes. For example, if \(N = \Ko \) is a knot complement, and we 
let  \(S_k (K) \) be the subcomplex of \( \cfr (K) \) generated by those
elements with Alexander grading less than \(k\), then we have 

\begin{prop}
There is a short exact sequence 
\begin{equation*}
\begin{CD}
0 @>>> S_{-k} (K) \oplus S_{k} (K) @>>> \cfr (K, \spi _k) @>>> \scx _2 ^k (K)
@>>>0.
\end{CD}
\end{equation*}
\end{prop}

\begin{figure}
\includegraphics{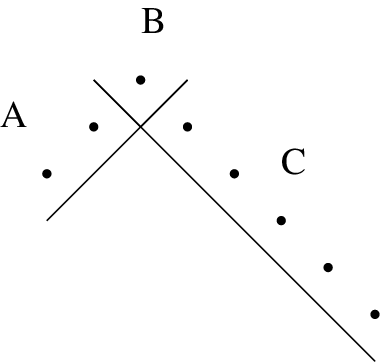}
\caption{\label{Fig:HatCx} A dot picture of \(\hfhat (\spi_k) \) }
\end{figure}

\begin{proof}
It is clear from Figure~\ref{Fig:HatCx} that we have a short exact
sequence 
\begin{equation*}
\begin{CD}
0@>>> A \oplus C @>>> \cfr (K, \spi _k) @>>> B @>>> 0 
\end{CD}
\end{equation*}
If we think of everything as being inside \( \cfir (\E) \), we see
that  \( C \cong S_k (K) \) and \( B \cong \scx_2 ^k (K) \). To identify \(A\)
with \(S_{-k} (K) \), we use the isomorphism \(\cfr (K) \cong \cfa (K)
\). 
\end{proof}

In general, if we know the filtered complexes \( \cfr (N,m, \spi) \)
and \( \cfr (N,m, -\spi) \) we can deduce the structure of 
\( \cfr (N,m, \spi _k) \) for every \( \spi _k\). 
The situation for \( \cfir (\E, \spi ) \) is more complicated. 
Knowledge of \( \cfr (N,m, \spi) \) tells us what the differentials
which go to the edge of the shaded triangle in
Figure~\ref{Fig:Filtrations} are, but does not give us any information
about those  differentials which go to the interior.

%% file: figs/basepoints.pstex_t
\begin{picture}(0,0)%
\includegraphics{figs/basepoints.pstex}%
\end{picture}%
\setlength{\unitlength}{2171sp}%
\begingroup\makeatletter\ifx\SetFigFont\undefined%
\gdef\SetFigFont#1#2#3#4#5{%
  \reset@font\fontsize{#1}{#2pt}%
  \fontfamily{#3}\fontseries{#4}\fontshape{#5}%
  \selectfont}%
\fi\endgroup%
\begin{picture}(11624,4173)(4179,-4942)
\put(12124,-1592){\makebox(0,0)[lb]{\smash{\SetFigFont{6}{7.2}{\rmdefault}{\mddefault}{\updefault}{\color[rgb]{0,0,0}j-2}%
}}}
\put(12709,-1607){\makebox(0,0)[lb]{\smash{\SetFigFont{6}{7.2}{\rmdefault}{\mddefault}{\updefault}{\color[rgb]{0,0,0}j-1}%
}}}
\put(13414,-1607){\makebox(0,0)[lb]{\smash{\SetFigFont{6}{7.2}{\rmdefault}{\mddefault}{\updefault}{\color[rgb]{0,0,0}j}%
}}}
\put(13969,-1652){\makebox(0,0)[lb]{\smash{\SetFigFont{6}{7.2}{\rmdefault}{\mddefault}{\updefault}{\color[rgb]{0,0,0}j+1}%
}}}
\put(14689,-1652){\makebox(0,0)[lb]{\smash{\SetFigFont{6}{7.2}{\rmdefault}{\mddefault}{\updefault}{\color[rgb]{0,0,0}j+2}%
}}}
\put(12169,-3917){\makebox(0,0)[lb]{\smash{\SetFigFont{6}{7.2}{\rmdefault}{\mddefault}{\updefault}{\color[rgb]{0,0,0}j-3}%
}}}
\put(12889,-3902){\makebox(0,0)[lb]{\smash{\SetFigFont{6}{7.2}{\rmdefault}{\mddefault}{\updefault}{\color[rgb]{0,0,0}j-2}%
}}}
\put(13369,-3917){\makebox(0,0)[lb]{\smash{\SetFigFont{6}{7.2}{\rmdefault}{\mddefault}{\updefault}{\color[rgb]{0,0,0}j-1}%
}}}
\put(14629,-3932){\makebox(0,0)[lb]{\smash{\SetFigFont{6}{7.2}{\rmdefault}{\mddefault}{\updefault}{\color[rgb]{0,0,0}j+1}%
}}}
\put(13954,-3932){\makebox(0,0)[lb]{\smash{\SetFigFont{6}{7.2}{\rmdefault}{\mddefault}{\updefault}{\color[rgb]{0,0,0}j}%
}}}
\put(12826,-961){\makebox(0,0)[lb]{\smash{\SetFigFont{10}{12.0}{\rmdefault}{\mddefault}{\updefault}{\color[rgb]{0,0,0}\(\beta _g\)}%
}}}
\put(13201,-4861){\makebox(0,0)[lb]{\smash{\SetFigFont{11}{13.2}{\rmdefault}{\mddefault}{\updefault}{\color[rgb]{0,0,0}b)}%
}}}
\put(5701,-4861){\makebox(0,0)[lb]{\smash{\SetFigFont{11}{13.2}{\rmdefault}{\mddefault}{\updefault}{\color[rgb]{0,0,0}a)}%
}}}
\end{picture}

%% file: theorem1.tex
\section{Invariance of \( \cfr \)}
\label{Sec:Thm1}

 The first goal of this section is to  explain
the process of reduction, show that it is well defined, and prove that a
filtration preserving chain map induces a chain map of the reduced
complexes. Once this has been taken care of, we turn to the proof of
Theorem~\ref{Thm1}.
We need to show that \( \cfr \) is does not depend on the
various choices we have made in defining it, such as 
the Heegaard splitting of \(N\),
the geometric representative of the meridian, and the position of the
basepoint \(z\). In each case, we examine the proof of invariance of 
\( \hfhat \) given in \cite{OS1} and check that it can be adapted to
show that  \( \cfr \) is invariant as well. 

As an example of how this process works, suppose we transform the
Heegaard diagram \( \E \) into a new diagram \( \E' \) by  a handleslide.
\Ozsvath and \Szabo exhibit a chain map \(F \: \cfhat (\E) \to
\cfhat (\E') \) and show that the induced map on homology is an isomorphism.
From our perspective, the key point is that since  \(F\) is defined
using holomorphic triangles, it preserves the Alexander
filtration. Thus there are induced  maps
 \( F^i \: \scx_1^i (\E) \to \scx _1 ^i
(\E') \). The argument used to show that \( \hfhat \) is invariant
under handleslide  actually implies that the maps 
  \(F^i_* \: \scx_2 ^i (\E) \to \scx _2 ^i (\E') \) are
isomorphisms. It is a standard fact (see {\it e.g.}  Theorem 3.1 of
\cite{UsersGuide}) that this  implies  that all the higher
terms in the spectral sequence are the same as well. It also enables us
to show that the associated reduced complexes are isomorphic.

\subsection{Reducing a filtered complex}
\label{Subsec:Reduction}

Throughout this section, we work with field coefficients.
% (There should be some
%better way to do all of this, but I can't see what it is going to be.) 
Recall the setup of Lemma~\ref{Lem:Reduction}: we have a complex \((C,d)
\) with a filtration \(C_1 \subset C_2 \ldots \subset C_m \)
and filtered quotients \( (C^i,d^i) = C_{i+1}/C_{i} \). We wish to find an
equivalent complex \( (C',d') \) such that \(d'^i \equiv 0\).
 The basic tool is the  well-known ``cancellation lemma'':

\begin{lem}
\label{Lem:Cancellation}
Suppose \((C,d)\) is a chain complex freely generated by elements \(y_i
\), and write \( d(y_i,y_j) \) for the \(y_j\) coordinate of \( d(y_i)
\). Then if \(d(y_k,y_l) = 1 \), we can define a new complex
 \((\overline{C},\overline{d})\) with generators \( \{
y_i \ts | \ts i \neq k,l\} \) and differential 
\begin{equation*}
\overline{d} (y_i) =d(y_i) + d(y_i,y_l)d(y_k)
\end{equation*}
which is chain homotopy equivalent to the first one. 
\end{lem} 

\noindent {\bf Remarks:} A proof of this fact may be found in
\cite{Floer1}.  The chain homotopy
equivalence \( \pi \: C \to \overline{C} \) is just the projection,
while the equivalence \(\iota \: \overline{C} \to C \) is given by 
\(\iota (y_i) = y_i -d(y_i,y_l)y_k \).

\begin{proof}
(Of Lemma~\ref{Lem:Reduction}) If all the \(d^i\) are \(0\), the
  result obviously holds. If not,  we can
  find \(x, y \in C^i \) with \(d^i(x,y) \neq 0 \). Since we
  are using field coefficients, we can assume that \(x \) and
  \(y\) are generators of \(C \) and scale \(y \) so that
\( d^i(x,y) = 1 \). This implies that \(d(x,y) = 1 \) as well,
  so we can apply the cancellation lemma to obtain a
  new  complex \( (\overline{C}, \overline{d}) \) chain homotopy
  equivalent to our original \(C\).

We claim that the filtration on 
  \(C\) induces a filtration on \( \overline{C} \). Indeed, suppose
  that 
\begin{equation*}
 \overline{d} (u, v)  = d(u,v) + d(u,y)d(x,v) \neq 0. 
\end{equation*} 
Then either \( d(u,v) \neq 0 \), which implies that \(u \geq v
\), or \(d(u,y)d(x,v) \neq 0 \), which implies \( u \geq y
\) and \(x \geq v \). By hypothesis, \(x \) and \(y \)
have the same filtration, so \( u \geq v \) in this case as well. 
Moreover, the chain homotopy equivalences \(\iota \) and \( \pi\)  are
 filtered maps. This is obvious for \(\pi\), and is true for \( \iota
 \) because \( d(u,y) \neq 0 \) implies that \(u \geq y\geq x \). 

The maps \( \iota_*^j\: H_*(\overline{C}^j) \to H_*(C^j)
  \) and \( \pi_*^j\: H_*({C}^j) \to H_*(\overline{C}^j)
  \)
are isomorphisms. Indeed, when
\(j \neq i \),  it is easy to see
that \( d^j = \overline{d}^j
\). On the other hand, \(\overline{C}^i \) is just the complex obtained from
\(C^i \) by cancelling \(x \) and \(y\). It follows that \( \iota
  \) and \(\pi\) induce an isomorphisms on the \(E_2 \) terms of the spectral
  sequences for  \(C\) and \( \overline{C} \). It is then a standard
  fact (see {\it e.g.} Theorem 3.1 of
\cite{UsersGuide}) that they induce isomorphisms on all higher terms as well.

We now repeat this process. Since \(C\) was
 finitely generated, the cancellation  must terminate after a finite number
 of steps. The result is the desired complex \((C', d' ) \). To show
 that \(C'\) is unique up to isomorphism, suppose we have another
 reduction \((C'',d'') \). Then by transitivity, \(C'\) is chain
 homotopy equivalent to \(C''\), and the induced map on the \(E_2\)
 terms of the spectral sequences is an isomorphism. It
 follows that the map \( C' \to C'' \) is an isomorphism at the level
 of groups. But a chain map which is an isomorphism on groups
 is  an isomorphism of chain complexes. 
\end{proof}

The following lemma says that we can reduce maps as well as complexes:

\begin{lem}
Suppose \(f\: C \to D \) is a filtration preserving map. Then there is
an induced map \(f'\: C' \to D' \). Moreover, 
\(f \) and \(f' \) induce the same maps on spectral sequences.
\end{lem}

\begin{proof}
Choose chain homotopy equivalences \( \iota \: C' \to C \) and \(\pi
 \:  D \to D'
\), and set \(f' = \pi  \circ f \circ \iota  \). Since \(\iota \) and \(\pi \)
induce isomorphisms on spectral sequences, the map on spectral
sequences induced by \(f' \) is the same as the map induced by \(f\). 
\end{proof}

We would also like to know that the reduced map \(f'\) is
unique.  Since we have only defined the reduced complex up to
isomorphism, the best we can hope to show is that any two reductions
\(f' \) and \(f''\) differ by composition with isomorphisms at either
end.
% {\it i.e.} \(f'' = \iota _1 \circ f' \circ \iota _2 \).
 This is
indeed the case. Supppose, for example, that  \( \iota_1 \: C' \to C \) is
another chain homotopy equivalence, and let \( \pi _1 \) be its
homotopy inverse. Then \(\iota \) is chain homotopic to
 \(  \iota_1 \circ \pi_1 \circ \iota  \), and \(\pi_1 \circ \iota 
 \) is an isomorphism from \( C' \) to itself. Finally, 
\(\pi \circ f \circ \iota \) is equal (not just homotopic!) to 
\( \pi \circ f \circ \iota _1 \circ (\pi_1 \circ \iota ) \), 
since both sides induce
the same map on \(E_2 \) terms of spectral sequences. A similar
argument handles the case of a different \(j\).

It is now easy  to prove Lemma~\ref{Lem:ReducedSequence}:

\begin{proof}
Suppose \(x,y \in B^i\) are a pair of cancelling elements, as in the
proof of Lemma~\ref{Lem:Reduction}. Since \( B^i \cong A^i \oplus C^i
\), \(x\) and \(y\) are either both in \(A^i\) or both in
\(C^i\). For the moment, we suppose the former. Then for \(u \in A
\), it is easy to see that \( \overline{d}|_A \) is the same as the
differential obtained by cancelling \(x \) and \(y\) in \(A\). On the
other hand, if \(z,w \in B \backslash A\),  \(\overline{d} (z,w)
= d(z,w) \). Thus we have a commutative diagram
\begin{equation*}
\begin{CD}
0 @>>> A @>>> B @>>> C @>>> 0 \\
@VVV @VVV @VVV @VVV @VVV \\
0 @>>> \overline{A} @>>> \overline{B} @>>> C @>>> 0
\end{CD}
\end{equation*}
in which the vertical arrows are chain homotopy equivalences. 
If \(x,y\) are in \(C\), there is similar diagram involving \(A\),
\(\overline{B} \) and \( \overline{C} \). 
We now repeat, stacking the  diagrams on top of each other as
we go. Compressing the resulting large diagram gives the statement of
the lemma.
\end{proof}

\subsection{Filtered chain homotopies} 
Let  \(C \) and \(D \) be chain complexes. A map \( i \:
C \to D \) is a chain homotopy equivalence if there exists a homotopy
inverse \(j \: D \to C \), and chain homotopies \( H_1 \: C \to C \),
\(H_2 \: D \to D \) such that 
\begin{align}
\label{Eq:CHT}
1_C - ji & = dH_1 - H_1d \\
1_D - ij & = dH_2 - H_2d
\end{align}
Now suppose that that \(C\) and \(D\) are both filtered complexes. We
say that \(i \) is a {\it filtered chain homotopy equivalence } if
\(i,j,H_1\), and \(H_2 \) all respect the filtration.

\begin{lem}
\label{Lem:CHT}
Let \(i : C \to D \) be a filtered chain homotopy equivalence. Then
all the induced maps \( i^k : C^k \to D^k \) are chain homotopy
equivalences as well.
\end{lem}

\begin{proof}
Since \(i, j, H_1\) and \(H_2 \) all respect the filtration, they
induce maps \(i^k,j^k,H_1^k\), and \(H_2 ^k \) on the filtered
quotients. It is easy to see that the analogues of equations
\ref{Eq:CHT}
and 3 hold for these maps as well. 
\end{proof}

We can now explain the strategy for proving Theorem~\ref{Thm1} in a
bit more detail.  First, we list the choices we made in
defining \( \cfs (\E) \):
\begin{enumerate}
\item A generic path of almost-complex structures \(J_s\) on \( s^g \Sigma
  \) and a generic complex structure \( \jjj \) on \( \Sigma \). 
\item A Heegaard splitting of \(N\).
\item The  geometric representative \( \emm \) of \(m\). 
\item The position of the basepoint \(z\).
\end{enumerate}

If we change any one of these things, we get a new stable complex,
which we  denote by 
\(\cfs (\E')\). Now in many cases, 
the proof in \cite{OS1}  that \( \hfhat \) is a
topological invariant gives us a chain homotopy equivalence \(i \:
\cfhat (\E) \to \cfhat (\E') \). In the sections below, we  check
that \(i\) is actually a filtered chain homotopy equivalence. Assuming
that this is the case, the remainder of the proof is straightforward.
Indeed, Lemma~\ref{Lem:CHT} tells us that \(i \) induces isomorphisms
\(i^k \: \scx _2 ^k (\E) \to \scx _2 ^k (\E') \). It follows that the
induced map \( i' \: \cfr (\E) \to \cfr (\E') \) is an isomorphism at
the level of groups, which implies that \( \cfr (\E) \cong \cfr (\E')
\).

As a model case, let us check that \( \cfr \) does not depend on
\(J_s \) or  \(\jjj \). This is proved for
\(\hfhat \) in Lemma 4.9 of \cite{OS1}, using a  standard argument in 
Floer theory. Given two paths \(J_s(0) \) and \(J_s(1)
\), one connects them by a generic homotopy of paths \(J_s(t)
\) and defines a chain homotopy equivalence \(\Phi _{J_{s,t}}\:
 \cfhat (\E) \to \cfhat (\E') \) by 
\begin{equation*}
\Phi _{J_{s,t}} (\bfy ) = \sum _{\bfz} \sum _{\phi \in D} \#\MMM_{J_{s,t}} 
(\phi) \cdot \bfz
\end{equation*}
where \( D= \{\phi \in \pi _2 (\bfy, \bfz) \ts | \ts \mu
  (\phi )  = n_{z_s} (\phi ) = 0\} \) and 
 \( \MMM_{J_{s,t}}(\phi) \) is the moduli space of time dependent
holomorphic strips associated to the path \(J_s(t) \). The domain of
  such a
  \( \phi \) looks just like the domain of a differential (except that
  it has Maslov index 0), so 
\( n_{z_a} (\phi) -  n_{z_s} (\phi) = A(\bfy) - A(\bfz) \). 
Since  \( n_{z_s}(\phi) = 0\), we must have 
\(A(\bfy) \geq A(\bfz) \) whenever \( \bfz \) has a nonzero
coefficient in the sum. Thus \( \Phi _{J_{s,t}} \) is a filtration 
preserving map. Its homotopy inverse  \( \Phi _{J_{s,1-t}} \) is
 also filtration preserving  (by the same argument). 
The chain homotopy between \(\Phi _{J_{s,t}} \circ  \Phi _{J_{s,1-t}}
  \) and the identity is given by
\begin{equation*}
H_{J_{s,t,\tau}} = \sum _{\bfz} \sum _{\phi \in D'} \#
\MMM_{J_{s,t,\tau}} (\phi) \cdot \bfz 
\end{equation*}
where \(J_{s,t,\tau} \) is a suitable two parameter family of almost
 complex structures and \(D' = \{\phi \in \pi _2 (\bfy, \bfz) \ts | \ts \mu
  (\phi ) = -1, \ts n_{z_s} (\phi ) = 0\} \). It is easy to see that 
\( H_{J_{s,t,\tau}} \) is also filtration preserving. 
 Thus Lemma~\ref{Lem:CHT} applies, and we
conclude that \( \cfr (\E) \cong \cfr (\E') \). The argument used in
 \cite{OS1} to show that \(\hfhat \) is independent of the generic
 almost complex structure \(\jjj \) on \(\Sigma \) carries over
 without change to \( \cfr \). 

In the following sections, we use arguments like this one
to  show that \(\cfr \) is invariant under isotopies, handleslides,
stabilizations, and changes of basepoint. Finally, in
section~\ref{Subsec:EndProof}, we put  these results together to prove
the theorem. 

\subsection{Isotopy invariance}

\begin{prop}
\label{Prop:MeridianIsotopy}
\(\cfr \) is invariant under isotopies of \(\emm\), \(\alpha_i\), and
\( \beta_j \) supported away from \( z_s \). 
\end{prop}

\begin{proof}
If the isotopy does not create or destroy intersections,
 it is equivalent to an isotopy of the conformal structure
on \( \Sigma \). We observed above that \( \cfr \) is invariant
under such isotopies. Thus it suffices to consider the case in which a
pair of intersection points between \(\emm\) or one of the \( \beta _j\)'s
 and one of the \( \alpha _i \)'s is
created or destroyed. Since an isotopy that moves 
\(\beta _j \) over \( \alpha _i
 \) is indistinguishable from one that moves \( \alpha _i \) over \( \beta
 _j \), we  assume without loss of generality that the curve that
 is moving is one of the \( \alpha _i \), and that \(\emm\) and all the
 \( \beta _j\)'s are fixed.  

 The isotopy invariance of \(\hfhat \) in this case is proved in Theorem
4.10 of \cite{OS1}. The argument is as follows. We represent the
isotopy of \( \alpha _i \) by a family of Hamiltonian diffeomorphisms
\( \Psi _t \) of \(\Sigma \),
 and consider the map \( F _{\Psi_t} \: \cfr (\E) \to \cfr (\E') \)
defined by 
\begin{equation*}
F  _{\Psi_t} (\bfy)
 = \sum _{\bfz \in E'} \sum_ {\phi \in D} \# \MMM ^{\Psi_t} (\phi) \cdot \bfz
\end{equation*}
where \( D =  \{\phi \in \pi_2^{\Psi_t} (\bfy, \bfz) \ts | \ts \mu (\phi) =
n_{z_s}(\phi) = 0 \} \). The moduli space \(\MMM^{\Psi_t} (\phi) \) counts
holomorphic discs subject to an appropriate time dependent boundary
condition determined by \( \Psi _t \). 
The usual Floer homology arguments show that
this is a chain map, and that \( F_{\Psi_{1-t}} \circ F_{\Psi_{t}} \:
\cfr (E) \to \cfr (E) \) is chain homotopic to the
identity. The chain homotopy is defined by choosing a homotopy \( \Phi
_{t,\tau}\) between \( \Psi
_{1-t} * \Psi _t \) and the identity and setting 
\begin{equation*}
 H(\bfy) = \sum _{\bfz \in E} \sum_{\phi \in D'} \# \MMM
^{\Phi_{t,\tau}(\phi)}(\phi) \cdot \bfz 
\end{equation*}
where \( D' = \{ \phi \in \pi _2 ^{\Phi_{t, \tau}}(\bfy, \bfz) \ts |
\ts n_{z_s} = 0, \ts \mu (\phi) = -1\} \). Notice that this proof
never uses the hypothesis that only a single pair of points are
created or destroyed. Thus it applies just as well when \(
\alpha _i \) isotopes over a meridian, which creates or destroys many
pairs of intersection points in the Heegaard splitting for \( N(n,1)
\). 

To prove invariance, we need to check 
 that \(F \) and \(H \) preserve the Alexander
filtration. We argue much as we did in proving invariance with respect
 to \(J_s\). The only real difference
 is that the \( \DDD (\phi) \) is a bit more difficult to control,
 since the  \( \alpha _i \) is only known to lie in \( \Psi_t (\alpha
 _i)\). By hypothesis, however, the 
 isotopy is supported away from \( z _s
\), so the multiplicity \( n_{z_s} (\phi) \) is  well
defined for \(\phi \in D \) or \(D'\).
It still makes sense to talk about the \( \beta _g \) component of
\( \partial \DDD (\phi) \), and  the usual argument still  shows that 
the \( \beta _g \) component 
is supported inside the spiral. Thus   \( n_{z_a} (\phi) = A(\bfy) -
A(\bfz) \) still holds.
 Since any \( \phi \) with \( \# \MMM (\phi) \neq 0 \) has 
\( n_{z_a}(\phi) \geq 0 \), \(F\) and \(H\) respect the Alexander
 filtration. The proposition now follows from Lemma~\ref{Lem:CHT}. 
\end{proof}

\subsection{Handleslide invariance}
The proof that \( \cfr \) is invariant under handleslide is
considerably more complicated than the preceding arguments. There
are three separate cases to consider. In the first, 
one of the \( \beta \)  handles in the Heegaard
 splitting of \(N\) slides over another; in the second, one of the \(
 \alpha \) handles slides over another; in the last,
 \(\mathfrak{m}\) slides over one of the \( \beta \) handles in the Heegaard
 splitting. 

We first outline the proof  that \( \hfhat \) is invariant under
 handleslide
 given by \Ozsvath and \Szabo in 
\cite{OS1}. We adopt their notation, which is as follows:
the two original sets of
attaching handles will be denoted by \( \alpha \) and \( \beta\).
 \( \delta \) is a set of attaching handles
obtained by  small isotopies of the \( \beta \) handles, so that each \(
\beta _ i \) intersects \( \delta _i \) at precisely two points.
We let \( \gamma _1 \) be the handle obtained by sliding
\( \beta _1 \) over \( \beta _2 \). The other \(\gamma _i\)'s are
 obtained by  small isotopies of the \( \beta _i\)'s, so that
\( \beta _i \cap \gamma _i \) and \( \gamma _i  \cap  \delta _i \)
 each consist of two points. 
\vskip0.05in
 The main steps of the
proof in \cite{OS1} are as follows.
\begin{enumerate}
\item Show that \( \hfhat (\TTT_\beta, \TTT _\gamma) \cong \hfhat(\TTT
  _\gamma, \TTT _ \delta) \cong \hfhat(\TTT
  _\beta, \TTT _ \delta) \cong H_*(T^g) \). Label the top generators
  of these three groups  \( \theta _1, \theta _2\), and \(  \theta _3 \)
  respectively.
\item Define maps 
\begin{align*}
 \hat{f}_{\alpha, \beta, \gamma} (\cdot \otimes \theta _1)\: \  & \cfhat
 (\TTT _\alpha, \TTT _\beta ) \to \cfhat ( \TTT _\alpha , \TTT
 _\gamma) \\
\hat{f}_{\alpha, \gamma, \delta} (\cdot \otimes \theta _2)\: \  & \cfhat
 (\TTT _\alpha, \TTT _\gamma ) \to \cfhat ( \TTT _\alpha , \TTT
 _\delta) \\
\hat{f}_{\alpha, \beta, \delta} (\cdot \otimes \theta _3)\: \  & \cfhat
 (\TTT _\alpha, \TTT _\beta ) \to \cfhat ( \TTT _\alpha , \TTT
 _\delta) 
\end{align*}
by counting holomorphic triangles. For example
\begin{equation*}
\hat{f}_{\alpha, \beta, \gamma} (\bfy \otimes \theta _1) =
\sum _ {} c(\bfy, \theta _1, \bfz) \cdot \bfz
\end{equation*}
where
\begin{equation*}
 c(\bfy, \theta _1, \bfz) = \sum_{ \{\psi \in \pi _2 (\bfy,
 \theta _1, \bfz) \ts | \ts \mu (\psi) = n_{z_s} (\psi) =  0\}} \#
 \MMM (\psi). 
\end{equation*}
Check that these maps are chain maps, and label the induced maps on
homology by \(\hat{F}_{\alpha, \beta, \gamma}(\cdot \otimes \theta _1) \),
\(\hat{F}_{\alpha, \gamma, \delta}(\cdot \otimes \theta _2) \), and
\(\hat{F}_{\alpha, \beta, \delta}(\cdot \otimes \theta _3) \)
\item Prove the associativity relation
\begin{equation*}
\hat{F}_{\alpha, \gamma, \delta} ( \hat{F} _{\alpha, \beta,
  \gamma}(\bfy \otimes \theta _1)\otimes \theta _2)
 = \hat{F} _{\alpha, \beta, \delta} (\bfy \otimes \theta _3).
\end{equation*}
\item Check at the chain level that \( \hat{f} _{\alpha, \beta,
  \delta} \) is the identity map, so that \(\hat{F} _{\alpha, \beta,
  \gamma} \) and \( \hat{F}_{\alpha, \gamma, \delta} \) are inverse
  isomorphisms between \( \hfhat (T_\alpha, T_\beta ) \) and \( \hfhat
  (T_\alpha, T_\gamma )\). 
\end{enumerate}

This argument can be extended to show that \( \cfr \) is invariant
 under handleslide as well. To do so, we must
 check that all the relevant maps respect the Alexander
 filtration. Our starting point is  the following fact about
 holomorphic polygons:

\begin{lem}
\label{Lem:HolPolygon}
 Let \( \psi \) be a holomorphic polygon whose sides lie on
 the Lagrangians \(\TTT_{\alpha}, \TTT _{\beta}\), {\it etc}. Suppose
 futher that \(\bfy\) and \(\bfz\) are two corners of \(\psi\) such
 that as we go from \( \bfz \) to \( \bfy \) along \( \partial (\psi)
 \), all the sides contain a
 component parallel to \( \beta _g \), but as we go from \( \bfy \) to
 \(\bfz \), none of them do. Then if \(n_{z_s}(\psi) = 0\),
\( A(\bfy) \geq A(\bfz) \).
\end{lem}

\begin{figure}
\includegraphics{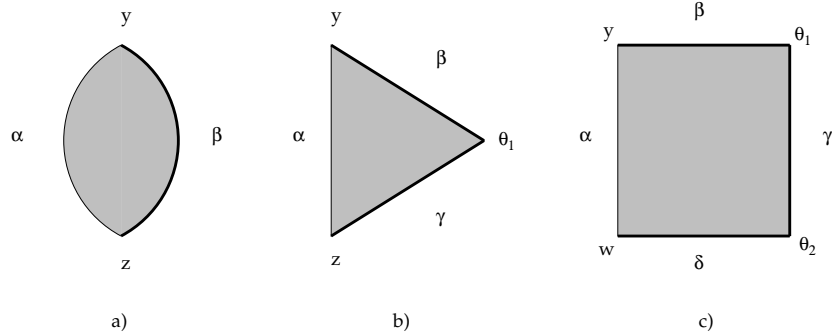}
\caption{\label{Fig:FiltProof} Holomorphic disks, triangles, and
  rectangles used to define {\it a)} \( d(\bfy, \bfz) \), {\it b)} 
 \(F_{\alpha, \beta, \gamma} (\bfy \otimes \theta _1 )\),  and {\it c)}
 \(\MMM (\bfy, \theta _1,
  \theta _2, \bfz) \). The heavy edges have sides parallel to \( \beta
  _g \).} 
\end{figure}

This statement generalizes the  fact that differentials do not increase the
Alexander grading. Some relevant cases  are illustrated in
Figure~\ref{Fig:FiltProof}. 

\begin{proof}
The edges of \(\psi\) between \( \bfz \) and \( \bfy \) give us a path
which runs from \(w(\bfz)\) to \( w(\bfy ) \) along the various
parallel copies of \(\beta _g\). The first thing we do is push all of these
 copies on top of each other, to  get a path from 
\(w(\bfz)\) to \( w(\bfy ) \) along
 \(\beta _g \). We claim that this path is supported inside the spiral
 region. Indeed, if this was not true, the path would have
 support along the whole exterior part of the spiral. It
 is easy to see that \( \psi \) cannot be a positive region if this is
 the case.

Since the other edges of \( \partial (\psi) \) do not contain any
 components parallel to \( \beta _g \), we can compute
 from Figure~\ref{Fig:Basepoints} that \(n_{z_a}(\psi)  -n_{z_s}(\psi) =
 A(\bfy) - A(\bfz) \). This proves the lemma. 
\end{proof}

Now suppose we are slide one \( \beta \) handle in the Heegaard
splitting for \(N\) over another. For the moment, let us
 assume that the slide is
supported away from \(\beta _g\) (and in particular, away from \(m\).)
 Thus there is a corresponding handleslide on \( N(n,1) \). We study
 the maps on \( \cfs \) induced by this handleslide. 

\begin{lem}
\label{Lem:TriangleFilt}
The maps
\( \hat{f} _{\alpha, \beta, \gamma}, \hat{f} _{\alpha, \gamma, \delta}
\), and \( \hat{f} _{\alpha, \beta, \delta} \)  preserve the
Alexander filtration. 
\end{lem}

\begin{proof}

 The form of the triangles which  \(c(\bfy, \theta _1, \bfz)\) counts
 is shown in Figure~\ref{Fig:FiltProof}b. The lemma implies
 that \( A(\bfy) \geq A(\bfz) \) whenever
\( c(\bfy, \theta _1, \bfz)  \neq 0\), which is what we wanted to prove.
 The same argument shows that 
\( F_{\alpha, \gamma, \delta} \) and \( F_{\alpha, \beta, \delta} \)
 are filtration preserving as well. 
\end{proof}

It follows that there are induced maps 
\begin{align*}
\hat{F}^i_{\alpha, \beta, \gamma} (\cdot \otimes \theta _1)
 &\: \scx_2^i (\TTT _\alpha, \TTT _
\beta) \to \scx_2^i (\TTT _\alpha, \TTT _\gamma) \\
\hat{F}^i_{\alpha, \gamma, \delta} (\cdot \otimes \theta _2)
&\: \scx_2^i (\TTT _\alpha, \TTT _
\gamma) \to \scx_2^i (\TTT _\alpha, \TTT _\delta) \\
\hat{F}^i_{\alpha, \beta, \delta} (\cdot \otimes \theta _3)
 &\: \scx_2^i (\TTT _\alpha, \TTT _
\beta) \to \scx_2^i (\TTT _\alpha, \TTT _\delta) 
\end{align*}
Our next step is to check that associativity still holds:
\begin{lem}
\label{Lem:Associate}
\begin{equation*}
\hat{F}^i_{\alpha, \gamma, \delta} (\hat{F}^i_{\alpha, \beta,
  \gamma}(\bfy \otimes \theta _1) \otimes \theta _2) = 
\hat{F}^i_{\alpha, \beta, \delta}(\bfy \otimes \hat{F}_{\beta, \gamma,
  \delta}( \theta _1 \otimes \theta _2))
\end{equation*}
\end{lem}

\begin{proof}
We recall the proof of associativity in Lemma 5.11 of \cite{OS1}.
 By examining the ends of an appropriate moduli space \( \MMM
(\bfy, \theta _1, \theta _2, {\bf w}) \) of psuedoholomorphic rectangles,
one sees that the statement of the lemma is false at the chain level,
but that the error is given by 
\begin{equation*}
\sum _{{\bf w'}\in \TTT_{\alpha} \cap \TTT_{\delta}} \# \MMM(\bfy,
 \theta_1, \theta _2, {\bf w'}) \cdot
d ({\bf w'}) 
\end{equation*}
Since \( d({\bf w'})\) is exact, 
the statement is true at the level of homology.

The same argument will  prove our version of the
statement if we can show that   for \( \bfy \in \scx^i_1 (\TTT_\alpha,
\TTT_{\beta}) \), the error term is exact in \(\scx^i_1 (\TTT_\alpha, 
\TTT_\delta)\). To prove this, it suffices to check
\(\# \MMM(\bfy, \theta_1, \theta _2, {\bf w'}) = 0 \) whenever \(
A ({\bf w'}) > A (\bfy) \).
Referring to Figure~\ref{Fig:FiltProof}c, we see that this follows
from  Lemma~\ref{Lem:HolPolygon}.
\end{proof}

\begin{prop}
\label{Prop:HandleSlide}
\( \cfr (N,m,\spi) \) is invariant under \( \beta \) handleslides
in \(N\) whose support is disjoint from \(\emm\). 
\end{prop}

\begin{proof}
If necessary, we slide and isotope \( \beta _g \) so that it is
disjoint from the support of the handleside in question. Since
\( \cfs (\E ) \cong \cfs (\E, \Nbar ) \),
 this has no  effect on \( \cfs (\E) \). 
Thus we may assume that the handleslide on our Heegaard splitting for
\(N\) extend to a handleslide on our splitting for \( N(n,1) \).

We have a filtered chain map  \( \hat{f}_{\alpha, \beta, \gamma}
 : \cfs (\TTT _\alpha,
  \TTT _{\beta}) \to \cfs (\TTT _\alpha, \TTT _{\gamma})\). To prove
  that the associated reduced map is an isomorphism, it suffices to
 show that the maps \(  \hat{F}^i_{\alpha, \beta, \gamma} : \scx ^i _2 (
\TTT _\alpha, \TTT _{\beta}) \to \scx ^i _2 (
\TTT _\alpha, \TTT _{\gamma}) \) are isomorphisms. This follows from
the associativity relation proved in the previous lemma, the identity
\( \hat{F}_{\beta, \gamma , \delta } (\theta _1 \otimes \theta _2) = \theta
 _3 \), and the fact (Lemma 5.12 of \cite{OS1}) that 
\( \hat{f}_{\alpha, \beta, \delta} : \cfhat (\TTT _\alpha , \TTT _
 \beta) \to \cfhat (\TTT _\alpha , \TTT _\delta ) \) is a filtration
 preserving  isomophism.
\end{proof}

Although this proof is more complex than our previous arguments
involving filtered chain homotopies, it is  similar in spirit. We still have
a filtered chain map \( \cfs (\E)  \to  \cfs (\E') \). The role of
the chain homotopy is now played by the associativity relation,
 and the fact that the chain
homotopy is filtered is  replaced by the fact that the relevant holomorphic
rectangles  respect the Alexander filtration. 
\vskip0.05in

Next,  suppose  that we are sliding an \( \alpha \) handle in the Heegaard
splitting of \(N\). We modify our labeling so that \( \gamma \) is the
set of \( \alpha \) handles after the handleslide, and \( \delta \) is
a parallel copy of \( \alpha \). There are chain maps 
\begin{align*}
 \hat{g}_{\alpha, \gamma, \beta} (\cdot \otimes \theta _1)\: \  & \cfhat
 (\TTT _\alpha, \TTT _\beta ) \to \cfhat ( \TTT _\gamma , \TTT
 _\beta) \\
\hat{g}_{ \gamma, \delta, \beta} (\cdot \otimes \theta _2)\: \  & \cfhat
 (\TTT _\gamma, \TTT _\beta ) \to \cfhat ( \TTT _\delta , \TTT
 _\beta) \\
\hat{g}_{\alpha, \delta, \beta } (\cdot \otimes \theta _3)\: \  & \cfhat
 (\TTT _\alpha, \TTT _\beta ) \to \cfhat ( \TTT _\delta , \TTT
 _\beta) 
\end{align*}
analogous to the \( \hat{f}\)'s. For example, \( 
\hat{g}_{\alpha, \gamma, \beta} \) is defined by 
\begin{equation*}
\hat{g}_{\alpha, \gamma, \beta} (\bfy \otimes \theta _1) = 
\sum c(\bfy, \bfz, \theta _1) \cdot \bfz
\end{equation*}
where
\begin{equation*}
c(\bfy, \bfz, \theta _1) = \sum _{ \{ \psi \in \pi _2 (\bfy, \bfz, \theta
  _1) \ts | \ts \mu (\psi) = n_{z_s} ( \psi ) = 0 \}} \# \MMM (\psi).
\end{equation*}
(The easiest way to see  this is the right definition is to
reverse  the orientation of \( \Sigma \) and switch the
roles of the \( \alpha \) and \( \beta \) handles. This leaves \( \cfs
(\E ) \) unchanged, but now we are sliding ``\( \beta \)'' handles, so
we can use the definitions of \cite{OS1}. Translating back to our
original setting, we arrive at the definition above.)

\begin{prop}
\label{Prop:HandleSlide2}
\(\cfr \) is invariant under the action of sliding one \( \alpha \) handle
over another.
\end{prop} 

\begin{figure}
\includegraphics{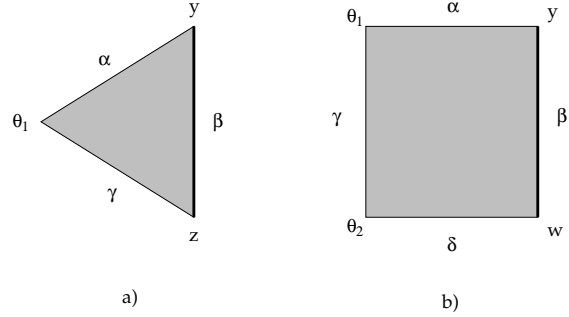}
\caption{\label{Fig:FiltProof2} Holomorphic triangles and
  rectangles used to define {\it a)}   \( G_{ \alpha
  , \gamma, \beta} (\bfy \otimes \theta _1 )\) and {\it b) }
\( \MMM (\bfy, {\bf w}, \theta _2, \theta _1) \). The heavy edges have
  a side parallel to \( \beta _g \).} 
\end{figure}

\begin{proof}

We argue as in the proof of Proposition~\ref{Prop:HandleSlide}. First,
we need to check that the \( \hat{g}\)'s respect the Alexander
filtration. The form of the triangles we count to define \( \hat{g}_{\alpha,
  \gamma, \beta }\) is shown in Figure~\ref{Fig:FiltProof2}a. Since \(
  \alpha \) and \( \gamma \) do not have any sides parallel to the
  spiral handle \( \beta _g \), Lemma~\ref{Lem:HolPolygon} applies. It
  follows that \( c(\bfy, \bfz, \theta _1 ) = 0 \) unless \( A(\bfy)
  \geq A(\bfz) \), so \( \hat{g}_{\alpha,\gamma, \beta} \) respects the
    Alexander filtration. Similar arguments apply to the other \(
    \hat{g}\)'s. 

Next, we need to check that associativity holds, {\it i.e.} that
\begin{equation*}
\hat{G}^i_{\gamma, \delta, \beta}(\hat{G}^i_{\alpha, \gamma, \beta}
(\bfy \otimes
\theta _1) \otimes \theta _2) = \hat{G}^i_{\alpha, \delta, \beta}( \bfy
\otimes \hat{G}_{\alpha, \gamma, \delta}(\theta _1 \otimes \theta _2))
\end{equation*}
For the unfiltered maps (without the \(i\)'s), this
 is proved in the usual way, by studying degenerations of the
space of holomorphic rectangles of the sort illustrated in 
Figure~\ref{Fig:FiltProof2}b. 
As in Lemma~\ref{Lem:Associate}, the statement is
 false at the chain level, but the error is 
\begin{equation*}
\sum _{{\bf w}' \in \TTT _ \delta \cap \TTT _ \beta} \# \MMM (\bfy,
 {\bf w}', \theta_2, \theta _1) d({\bf w}')
\end{equation*}
Referring to the rectangle of Figure~\ref{Fig:FiltProof2}b, we see
that Lemma~\ref{Lem:HolPolygon} applies in this case as well, so
\( \# \MMM (\bfy, {\bf w}', \theta _2, \theta _1) = 0 \) unless
\(A( \bfy) \geq A({\bf w}')  \). This proves the desired relation.

It is straightforward to check that the remaining elements of the
proof of handleslide invariance --- the fact that \(
\hat{G}_{\alpha, \gamma, \delta}(\theta _1 \otimes \theta _2) = \theta _3 \)
and that \( \hat{g}_{\alpha, \delta, \beta}( \cdot \otimes \theta _3) \) is
an isomorphism at the chain level --- extend to the \( \alpha \)
handleslide case. The result now follows by the same arguments
as Proposition~\ref{Prop:HandleSlide}.
\end{proof}

We can also use handleslides to modify \(\emm\). Indeed, it is not
difficult to see that the curve \(m'\) obtained by sliding \(\emm\) over
any one of the \( \beta \) handles in the Heegaard splitting of \(N\)
is another geometric representative of the meridian. We would like to
show that \( \cfr \) is invariant under the operation of replacing
\(\emm\) by \(\emm'\). 
 {\it A priori}, this situation seems to
be more complicated than the previous ones,
 since the corresponding move in the closed manifold \(N(n,1) \) is
 not a handleslide.
 Fortunately, Lemma~\ref{Lem:SpiralCrush}
 tells us that \( \cfs (\E, \spi) \cong \cfs
(\E, \Nbar, \spibar) \) as a filtered complex. Thus we can  study the
situation in \( \Nbar \), where sliding \(\emm\) is an honest
handleslide, albeit one involving the special handle. 

\begin{prop}
\label{Prop:MeridianSlide} 
\(\cfr \) is invariant under the operation of sliding \(\emm\) over
one of the \(\beta _i\)'s, so long as the slide is supported away from
\( z_s \). 
\end{prop}

\begin{proof}
Let \( \Tbar _\beta, \Tbar _\gamma \) and \( \Tbar _ \delta \) be the
tori in \( s^g \Sigma \) induced by the \( \beta\), \( \gamma \) and
\( \delta \) handles in the Heegaard splittings for \(\Nbar\). (We
revert to our original notation in which \( \gamma \) and \( \delta \) are
derived from \( \beta \).)
As usual, there are maps 
\begin{align*}
 \hat{f}_{\alpha, \beta, \gamma} (\cdot \otimes \theta _1)\: \  & \cfhat
 (\TTT _\alpha, \Tbar _\beta ) \to \cfhat ( \TTT _\alpha , \Tbar
 _\gamma) \\
\hat{f}_{\alpha, \gamma, \delta} (\cdot \otimes \theta _2)\: \  & \cfhat
 (\TTT _\alpha, \Tbar _\gamma ) \to \cfhat ( \TTT _\alpha , \Tbar
 _\delta)  \\
\hat{f}_{\alpha, \beta, \delta} (\cdot \otimes \theta _3)\: \  & \cfhat
 (\TTT _\alpha, \Tbar _\beta ) \to \cfhat ( \TTT _\alpha , \Tbar
 _\delta)
\end{align*}
which induce isomorphisms on \( \hfhat \). 

Although there is no longer a spiral handle, it still makes sense to
talk about the Alexander grading on \( \cfs (\E, \Nbar) \)
--- it is defined by the relation \(A(\ybar) - A(\ybar) = n_{z_s} (\phi)
- n_{z_a} (\phi) \) for \( \phi \in \pi_2(\ybar, \zbar) \). 

We would like to show that
the \( \hat{f}\)'s respect  the Alexander filtration. To do this,
observe that for every \( \ybar \in \TTT _ \alpha \cap \Tbar _\beta \) 
there is a  ``close'' point \( \ybar ' \in 
\TTT _ \alpha \cap \Tbar _\gamma \). Suppose 
\( \ybar \in \TTT _ \alpha \cap \Tbar _\beta \) and \( \zbar \in 
\TTT _ \alpha \cap \Tbar _\gamma \). Then to any 
\( \psi \in \pi _2 (\ybar, \theta _1, \zbar )\), there is  associated a disk
\( \psi ' \in \pi _2 (\ybar', \zbar) \). To define \( \psi ' \), we
superimpose the \( \beta \) and \( \gamma \) handles, so they are
indistinguishable except in the neighborhood of the handleslide shown
in Figure~\ref{Fig:MeridianSlide}. Then \( \partial \DDD (\psi') \)
looks just like  \( \partial \DDD (\psi) \), except that every time
the segment \( B \) appears in \( \partial \DDD (\psi) \), we
replace it with the segment labeled \(C \) 
together with an entire copy of \( \gamma _2
\). Since \( C \cup \gamma _2 \) is homologous to \( B \), the
multiplicities of \( \DDD (\psi ) \) and \( \DDD (\psi ' ) \) agree
outside of the shaded region in the figure. 

By hypothesis the handleslide was supported away from \( z\), so 
\( n_{z_a} (\psi') = n_{z_a}(\psi) \). Thus for \( \DDD (\psi ) \) to
be positive, we must have \( A(\bfy ') - A( \bfz) = n_{z_a} (\psi ')
\geq 0 \).  
Thus \( \hatf_{\alpha, \beta, \gamma} \) is a filtered map. Similar
arguments (always comparing triangles with differentials in 
\( \cfhat (\TTT_{\alpha}\cap \Tbar_{\gamma}) \)) show that the other \(
\hat{f} \)'s are filtered and that the associativity relation
holds. The proof then proceeds along the lines of
Proposition~\ref{Prop:HandleSlide}. 

\begin{figure}
\includegraphics{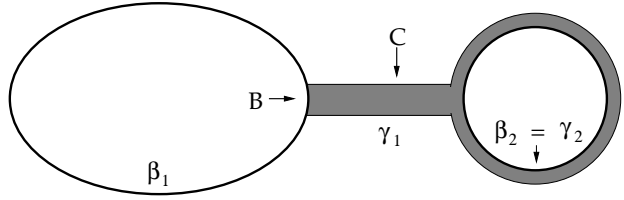}
\caption{\label{Fig:MeridianSlide} Comparison of the domains of \(
  \psi \) and \(\psi '\). }
\end{figure}
\end{proof}

\subsection{Stabilization invariance}
\begin{lem}
\label{Lem:Stabilization}
\( \cfr \) is invariant under stabilization of the Heegaard splitting
of \(N\).
\end{lem}

\begin{proof}
Let \( \cfr (\E') \) be the reduced complex for the stabilized Heegaard
splitting.
We claim that \( \cfr (\E') \) is independent of the point \( \sigma
\in \Sigma - \alpha - \beta \) at which we stabilized. Indeed,
stabilizing at another point \( \sigma ' \) is the same as
stabilizing at \( \sigma \) followed by a sequence of handleslides. 
By Propositions~\ref{Prop:HandleSlide} and \ref{Prop:HandleSlide2},
 we know that \( \cfr \) does not
change under handleslide. 
%(Since \( \Sigma - m \) is connected, we don't
%even need to worry about sliding across \(m\).)
 
Thus we may as well assume that \( \sigma \)
is in the same component of \( \Sigma - \alpha - \bbar \) as the
point \(z_s\). 
In this case, we are in the situation of Theorem 6.1 of \cite{OS1}, which
proves that \( \hfhat \) is invariant under stabilization by showing
that \(\cfs (\E) \) and \( \cfs (\E') \) are identical.  Combined with
the obvious fact that the Alexander grading on \( \cfs(\E') \) is the
same as the Alexander grading on \( \cfs (\E) \), this proves the
proposition. 

\end{proof}

\subsection{Invariance under change of basepoint}

The situation with regard to our choice of basepoint is rather
different than for \( \hfhat\). In the latter case, any
basepoint will do, but in order to define \( \cfs \), we required
\(z_s \) to be adjacent to \(\emm\) and on the stable side (as defined by
our choice of longitude   back in
section~\ref{Ssec:GenAlex}.)
Thus the only real way that we can vary \( z_s \) is to move it over
one of the \( w_i\)'s.

\begin{lem}
\label{Lem:Basepoint}
\( \cfr (N) \) does not change if we move \( z_s \) 
over one of the segments \(w^i\).
\end{lem}

\begin{proof}
This is immediate from the proof that  \( \hfhat
\) is invariant under change of basepoint  given Theorem 5.15 of
\cite{OS1}. Indeed, the argument there is based
on the observation that we can obtain the effect of moving \( z\) over
\( \alpha _i \) by a series of isotopies and handleslides of
\(\alpha_i \) supported away from \(z \). But we  know from
Propositions \ref{Prop:MeridianIsotopy} and 
\ref{Prop:HandleSlide2} that \( \cfr \) is invariant under these moves.
\end{proof}

Note that it is usually not true that the complex \( \cfr(N,m,\spi)
\) is the same as the complex \( \cfr^a(N,m,\spi) \) that we get by moving the
basepoint to the other side of \(m\); this is only the case when
\(\spi = \spibar\). This is the reason we required
\( z_s \) to be on the side corresponding to positive \( \spinc \)
structures. 

\subsection{Conclusion of proof}
\label{Subsec:EndProof}
To wrap things up, we need a pair of geometric lemmas.

\begin{lem}
Let \(N\) be a knot complement in a general three-manifold. Any two Heegaard
splittings of \(N\) differ by a sequence of isotopies, handleslides,
stabilizations, and destabilizations.
\end{lem}

\noindent This is proved in the same way as the closed case.

\begin{lem}
For a fixed Heegaard splitting of \(N\), 
any two geometric representatives of \(m\) are connected
by a series of isotopies and handleslides over the \( \beta _i\)'s.
\end{lem}

\begin{proof}
Recall that \( \partial N \) is obtained by surgering \(
\Sigma \) along \( \beta _1, \beta _2, \ldots \beta _{g-1} \). Thus we can
think of \( \partial N\)
 as a torus with \( 2g-2 \) small marked disks (the surgery disks) on it. Any
 two geometric representatives of \(\emm\) are connected by an isotopy on
 this torus, and it is easy to see that this isotopy extends to \(
 \Sigma \) except when \(\emm\) crosses over one of the marked disks. At
 these points, \(\emm\) slides over the corresponding \( \beta _i \) on
 \(\Sigma \). 
\end{proof}

\begin{figure}
\includegraphics{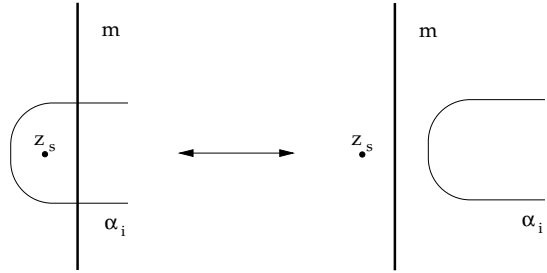}
\caption{\label{Fig:BaseMove} To do this isotopy, first move \(z_s \)
  over \( \alpha _i \).}
\end{figure}

We now prove of Theorem~\ref{Thm1}. 
Suppose we have two Heegaard splittings \((\Sigma, \alpha, \beta) \)
and \( (\Sigma', \alpha', \beta') \) with meridians \(\emm\) and \(\emm'
\) and basepoints \(z_s \) and \(z_s'\).
 We start by transforming the first Heegaard splitting into the second by a
sequence of isotopies, handleslides, stabilizations, and
destabilizations. We have already checked that these moves preserve \(
\cfr \). As we make these moves, however, we may also have to change our
choice of \(z_s\) and \(\emm\), since we may have isotoped an \( \alpha \)
handle over \(z_s\) or a \( \beta \) handle over \(\emm\). In the former
case, we use Lemma~\ref{Lem:Basepoint} to see that we can move \(z_s\)
over the \( \alpha \) handle without changing \( \cfr \). In the
latter, we isotope and slide \(\emm\) to keep it from intersecting any of
the \( \beta \) handles. Propositions~\ref{Prop:MeridianIsotopy} and
\ref{Prop:MeridianSlide} ensure that \( \cfr \) is invariant under
these moves as well. 
As we isotope \(\emm\), we move \(z_s \) along with it, so that it is
always adjacent to \(\emm\). The only situation in which it is not
possible to do this is shown in Figure~\ref{Fig:BaseMove}. In this case,
 the problem is rectified simply by moving \(z_s\) over 
\(\alpha _i \) before doing the isotopy. 

Once we have turned \((\Sigma, \alpha, \beta) \) , into 
\( (\Sigma', \alpha', \beta') \) we
isotope and slide the resulting \(\emm\) until it is the same as
\(\emm'\). (In doing so, we may need to move \(z_s\) over some more \(
\alpha \) handles.) We have thus transformed our first diagram into
a diagram which is identical to our second diagram except for the
position of the basepoint. Since \(z_s \) is by definition on the side
of \(\emm\) corresponding to the positive \( \spinc \) structures, the
resulting \(z_s \) is on the same side of \(\emm\) as \(z_s'\). Thus we can
 move  \(z_s\) over \(
\alpha \) handles until it is in same position as \(z_s'\). This
concludes the proof. 

%% file: quasiperfect.tex
\section{Perfect Knots}
\label{Sec:Quasiperfect}

Up to this point, we have considered the \OS Floer homology from a
purely theoretical standpoint. We now describe a case in which we can
actually compute the groups in question.

\begin{defn}
Let \(K\) be a knot in \(S^3\). We say that \(K\) is perfect if the
homological and Alexander gradings on \( \cfr (K) \) are the same.
\end{defn}

If \(K\) is perfect, the rank of \( \cfr ^{(j)} (K) \) is equal to the
absolute value of the coefficient of \(t^j\) in 
\( \Delta_K (t) \). Perfection is thus a very strong
condition on 
\(K\). Despite this fact, it is satisfied by a surprisingly large
number of knots.
In \cite{2bridge}, we showed that two-bridge knots are
perfect and described how to compute their \OS Floer homology. In
fact, the
arguments of that paper apply equally well to any perfect knot. This
fact was alluded to in \cite{2bridge}, but with a more restrictive
definition of ``perfect,'' which required the Alexander and
homological gradings to agree on \( \cfs (E) \). We will show in
section \ref{Sec:Alternating} that the new definition is a substantial
improvement: there are many knots whose stable complexes are not
perfect, but whose reduced stable complexes are. 

In this section, we repeat the calculation of \cite{2bridge}, both for
the sake of completeness and  because our improved understanding
of the Alexander filtration enables us to give a (hopefully) cleaner
exposition. We conclude by showing that the the class of perfect knots
is closed under connected sum. 

\subsection{Preliminaries}
 We fix  some basic conventions regarding
perfect knots. As usual, we choose a Heegaard splitting \( (\Sigma,
\alpha, \beta ) \) for \(\Ko \), twist up, and choose some good \(
\epsilon\)-class \(E\). From the discussion in
section~\ref{SubSec:ReducedCFI}, we
know that \( \cfp (\E, \spi_k) \) is chain homotopy equivalent to a
reduced complex, which we denote by \( \cfpr (K, \spi _k) 
\). (Although we have dropped the reference to \(\E\), we do not need
to know that \( \cfpr (K, \spi _k) \) is canonical in what follows ---
just that it exists.) 

\begin{lem}
If \(K\) is perfect, the
 Alexander and homological gradings on \( \cfpr (K, \spi _k ) \)
are the same. 
\end{lem}

\begin{proof}
The generators of \( \cfpr (K, \spi _k) \) are of the form \( [\bfy,
  i] \), where \(\bfy \) is a generator of \( \cfr (K) \). 
By definition, \(A([\bfy, i])= A(\bfy) + 2i \) and \( \gr ([\bfy, i])
= \gr (\bfy) + 2i\). Since \(K\) is perfect, \( A(\bfy) = \gr (\bfy) \).
\end{proof}

It is often convenient to think of the complexes \( \cfpr (K, \spi
_k) \) as  ``dot stacks'' like the one shown in
Figure~\ref{Fig:Filtrations}. Recall that each dot represents a
collection of generators, all with the same Alexander gradings. By
the lemma,  the vertical axis of such a figure
represents not only the Alexander grading on \( \cfpr (K, \spi _k )
\), but also the homological grading. Since the two are the same, we
often just refer to the grading of a generator in \( \cfpr (K, \spi
_k ) \), without specifying which one we are talking about. 

\begin{lem}
\label{Lem:Hfps}
 \( \hfps (K ) \cong \Zu  \) where \( 1 \in \Zu \) has grading \(
 s(K) \). Similarly \( \hfms (K) \cong \Z[u] \), where \(1 \in \Z[u]
 \) has grading \( s(K) - 2\). 
\end{lem}

\begin{proof}
Since \(K\) is a knot in \(S^3\), \( \hfs (K) \cong \Z \). The
generator has grading \(s(K) \) by definition. The obvious filtration
on \( \cfps (K) \) gives a spectral sequence with \(E_2 \) term \(
\hfs (E) \otimes \Zu \) converging to \( \cfps (K) \). Since 
\( \hfs (K) \cong \Z \), the sequence collapses at the \(E_2 \) term. 
The result for \( \hfms (K) \) follows from the analogous spectral sequence.
\end{proof}

In what follows, we will often consider the operation of {\it
  truncating} a chain complex \(C\) above  level \(k\). The
  truncated complex is just the quotient complex \( C_{\geq k} \)
  consisting of all elements of \(C\) with homological grading \( \geq
  k \). It is clear that \( H_*(C_{\geq k} ) \) is trivial for \( * <
  k \), is isomorphic to \( H_*(C) \) for \( *>k \), and has \( H_k
  (C) \) as a subgroup for \(* = k \). Similar results hold for the
  subcomplex \( C_{\leq k } \) obtained by truncating \(C\) below
  level \(k\). 

\subsection{Computing \(H_*(C_{\spi _k}(K))\)}

\begin{figure}
\includegraphics{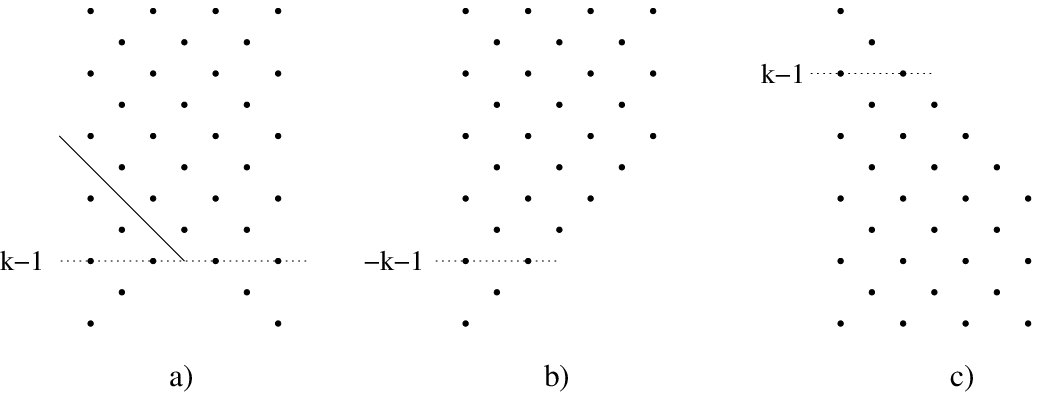}
\caption{\label{Fig:C_k} {\it a)}
 the subcomplex \(C_{\spi _k} (K) \), compared
with {\it b)} \(\cfpa (K)\) and {\it c)} \(\cfms (K)\).}
\end{figure}

Recall from Corollary~\ref{Cor:Seq1} that we have a short exact
sequence 
\begin{equation*}
\begin{CD}
 0 @>>> C_{\spi_k} @>>> \cfpr (K,\spi_k) @>>> \cfpr (K, \spi _s) @>>> 0 
\end{CD}
\end{equation*}
where \( C_{\spi_k} \) is the subcomplex depicted in
Figure~\ref{Fig:C_k}.

We will compute \( H_*(C_{\spi_k}(K)) \)
and then use the long exact sequence on homology to 
find \(\cfpr (K,\spi_k) \). Since \( \hfp
(K,\spi_k)  \cong \hfp (K,\spi_{-k})\), we may as well assume \( k
\geq 0 \). 

We denote the \( \Z[u] \)  module \( \Z[u]/(u^n) \) by \( \tmod _n \).

\begin{prop}
For  \(k \geq 0 \), 
\begin{equation*}
H_*(C_{\spi_k}(K)) \cong \tmod _{n_k} \oplus \Z ^{m_k}
\end{equation*}
where
\begin{align*}
n_k & =  \max (\lceil (|s(K)| - k)/2 \rceil, 0) \\
m_k & =  \chi (C_{\spi_k}(K)) -(-1)^{s(K) - k} n_k\\ & =  \sum _{i>k} (i-k)a_i
- (-1)^{s(K) - k}n_k
\end{align*}
where \(a^i \) is the coefficient of \(t^i \) in \( \Delta _K (t) \).
The grading of the \( \Z ^{m_k}\) summand is always  \(k-1\).
The grading of \(1 \in \tmod_{n_k} \) is \(s(K) - 2
\)  if \( s(K) > 0 \) and either \(k-1\) (if \(k \equiv s(K) \mod 2\) ) or
\(k-2 \) (if  \(k \not \equiv s(K) \mod 2\)) when \(s(K) < 0 \). 
\end{prop}

Less formally, \(H_*(C_{\spi_k}(K))\) is composed of two parts. The first
part is concentrated in grading \(k-1\) and has trivial \(u\)
action, and the second is a torsion module generated by a single
element. If \(s(K) > 0 \), this module starts at level \(s(K) \) and
stretches down to level \(k-1\) (or just above it), while if \(s(K) <
0 \) it starts at or just below \(k-1\) and extends downward. If
\(n_k \geq s(K) \), the torsion module is empty. 
Note that \( \chi ( C_{\spi _k}(K)) = \chi ( \hfp (K(0,1), \spi _k))
\). In fact,  in \cite{OS2}, \Ozsvath and \Szabo compute 
the right-hand Euler characteristic by showing it is equal to the
left-hand one. 

\begin{proof}
Observe from Figure~\ref{Fig:C_k} that if we truncate \(C_{\spi_k}(K)\) above
level \(k-1\), we get same complex that we do by truncating \( \cfms
(K) \) above level \(k-1\). Similarly, if we truncate \( C_{\spi_k}(K)
\) below
level \(k-1\), we get a complex which is isomorphic to \( \cfma
(K) \) truncated at level \(-k-1 \). It follows that 
\begin{align*}
H_*(C_{\spi _k}(K)) \cong 
\begin{cases}
\hfms_* (K) & \text{if} \  * > k-1 \\
\hfpa_{*-2k} (K) & \text{if} \  * < k-1.
\end{cases}
\end{align*}
For \( * \neq k-1 \), the claim now follows from Lemma~\ref{Lem:Hfps}
and the isomorphism \( \hfpa (K) \cong \hfps (K) \). Note that our
assumption that \( k \geq 0 \) ensures that one of \( \hfpa (K) \) or \(
\hfms (K) \) must be trivial over the range in question. 

Since we know the rank of \( H_*(C_{\spi _k}(K)) \) for \( * \neq k \),
we can use \( \chi
(C_{\spi _k}(K)) \) to compute the rank of  \( H_{k-1}(C_{\spi _k}(K))
\).  To see  this group is torsion free, we note that
the argument given above also applies if we use coefficients in \( \Z
/ p \). Since \(  H_*(C_{\spi _k}(K); \Z/p) \cong H_*(C_{\spi
  _k}(K)) \otimes \Z/p \) for \(* \neq k-1 \), the universal coefficient
theorem implies that \( H_{k-1}(C_{\spi _k}(K)) \) is free. 

To finish the proof, we need only check that the \( \Z[u] \) module
structure is as described.
 If \( \tmod _{n_k} \) and \( \Z ^{m_k} \) have different
parities, this is obvious from the grading. If they have the same parity,
 we must show, {\it e.g} that if \(s(K) > 0 \) and \(n_k > 1 \),
then \(H_{k-1}(C_{\spi_k}(K))\) contains an element in the image of
\(u\). But this is obvious, since \( H_{k-1} (C_{\spi_k}(K)) \) contains 
\( (\hfms)_{k-1}(K) \) as a submodule, and the (unique) generator of the
latter group is in the image of \(u\). The case when \(s(K) < 0 \) is
similar. 
\end{proof}

\subsection{Computing \( \hfp (K, \spi _k) \) }
We have a long exact sequence of \( \Z[u] \) modules
\begin{equation*}
\begin{CD}
@>>> H_*(C_{\spi_k}(K)) @>>> \hfp _* (K, \spi _k) @>g>> \hfps _*(K) @>>>
\end{CD}
\end{equation*}
Since \( \cfpr(K, \spi _k) \) is identical to \( \cfpr (K, \spi
_s) \) when the grading is large enough, \(\hfp (K, \spi _k) \) has a  \(
\Zu \) summand \( A\) which maps onto  \( \cfpr (K, \spi _s) \cong \Zu \).
Thus the long exact sequence reduces to a short exact sequence
\begin{equation*}
\begin{CD}
0@>>> H_*(C_{\spi_k}(K)) @>>> \hfp _* (K, \spi _k) @>g>> \Zu @>>> 0
\end{CD}
\end{equation*}
from which we can easily compute \( \hfp _* (K, \spi _k) \).

\begin{thrm}
\label{Thm:PerfectKnots}
For \(k \geq 0 \), 
\begin{align*}
\hfp (K, \spi _k) \cong
\begin{cases}
\Z^{m_k} \oplus \Zu & \text{if} \  s(K) \geq 0 \\
H_*(C_{\spi_k}(K)) \oplus \Zu & \text{if} \ s(K) \leq 0 
\end{cases}
\end{align*}
If \(s(K) \geq 0 \), the grading of \( 1 \in \Zu \) is \(s(K) - 2n_k
\), while if \( s(K) \leq 0 \) it is \( s(K) \). 
\end{thrm}

\begin{proof}
When viewed as abelian group, \( \hfp (K, \spi _k) \) decomposes as \(
 H_*(C_{\spi_k}(K)) \oplus \hfps (K) \), but this does not necessarily hold
 when we view it as a \( \Z[u] \) module. What we do know is that 
\( H_*(C_{\spi_k},(K)) \cong A \oplus B \), where \(A \cong \Zu \) and \(B
 \) is a finitely generated torsion module. \(A\) maps onto \(\hfps
 (K)\), but the map may have some kernel. 
To determine if this is the case, we
rely on the grading. By Lemma~\ref{Lem:Hfps}, the grading of \( 1 \in
 \hfps (K) \cong \Zu \) is \(s(K) \). We now analyze the cases \( s(K)
 \geq 0 \) and \( s(K) \leq 0 \) separately.

Suppose first that \(s(K) \leq 0
 \). In this case, the lowest
 grading attained by any element of \(  H_*(C_{\spi_k}(K)) \) is \(k-1 \)
 if \(n_k\) = 0, and either  \(k-2n_k \) or \(k+1 - 2n_k \) if \(n_k
 > 0 \). In particular, this lowest grading is always greater than or
 equal to \(s(K) -1 \). It follows that the grading of \( 1 \in A \)
 cannot be less than \(s(K) \). (It is congruent to \( s(K) \mod 2\).)
 Thus it must be equal to \(s(K) \), and the map \(A \to  \hfps(K) \) is an
 isomorphism. It follows that the short exact sequence splits as a
 sequence of \( \Zu \) modules. This proves the theorem when \(s(K)
 \leq 0 \). 

When \(s(K) \geq 0 \), we consider two subcases. If \(k \geq s(K) \),
\(n_k = 0 \), and all the elements of \( H_*(C_{\spi_k}(K)) \) have
grading \(k-1\). Since this is greater than or equal to \(s(K) - 1\),
 the same argument shows that the exact sequence splits. On the other
 hand, when \(k < s(K) \), we claim that the grading of \(1 \in A \)
 is either \(k\) or \(k-1\) (depending on the parity of \(k\).)
To see this, we refer back to Figure~\ref{Fig:C_k} and observe that
if we truncate \( \cfpr (K, \spi _k ) \) above level \(k\), we get 
\( \cfi (K) _{\geq k} \). Since \( \hfi (K) \cong \Z[u, u^{-1}] \), 
\( \hfp (K, \spi _k ) \) contains a \(\Zu \) summand with an element
of grading either \(k\) or \(k-1\). The grading of
every element in \( \hfp (K, \spi _k ) \) is at least \(k-1\), so the
claim is true. Thus the \( \tmod _{n_k} \) summand of \(
H_*(C_{\spi_k}(K)) \) is the kernel of the map \( A \to \Zu \), and 
\( \hfp (K, \spi _k) \cong \Z^{m_k} \oplus \Zu \). 
\end{proof}

\vskip0.1in
\noindent{\bf Remark:} The method of computation employed here is not
strictly limited to perfect knots. There exist knots which
 are not perfect, but which still have the
property that \(C_k \) looks like \( \cfms (K) \) when truncated above
some level \(l \) , and like \( \cfpa (K) \) when truncated below
\(l\). The argument above can be easily modified to compute \( \hfp
(K, \spi _k ) \) in such cases. Examples of such knots include the
torus knots and the knot \(10_{139} \).

\subsection{\(s(K)\) and \(\sigma (K)\)}
\label{Subsec:Signature}
Theorem~\ref{Thm:PerfectKnots} can be interpreted as saying that for a
perfect knot \(K\), the groups \( \hfp (K, \spi_k ) \) are completely
determined by the Alexander polynomial of \(K\) and the invariant
\(s(K) \). For all perfect knots we are aware of, \(s(K) \) is
determined by the classical knot signature: \( s(K) = \sigma (K) /2
\). In fact,
 \Ozsvath and \Szabo have proved in \cite{OS8}
 that this equality holds for all  alternating knots. 

 It would be nice to have a geometrical explanation for this
 phenomenon. One way to look at the matter is to say that for
 small perfect knots, \(s(K) \) and \( \sigma (K) \) are determined by
 the same set of skein theoretic axioms. If \(L\) is an oriented link in
 \(S^3\),  we can
 inductively construct a {\it skein resolution tree} \(T(L)\) for \(L\) as
 follows. If \(L\) is the unlink, \(T(L) = {L}\). If \(L\) is not the
 unlink, pick a crossing of \(L\), and let \( L_0 \) and \(L'\) be the
 links obtained by resolving the crossing in an oriented fashion and
 changing the crossing. (Of course, such trees are not unique.) Given
 a skein tree for a knot \(K\), we can use the skein relation to
 inductively compute \( \Delta _K(t) \). It is also possible to
 compute \( \sigma (K) \) from a skein tree. (See, {\it e.g}
 \cite{Kauffman}, p. 207.) 

Now suppose that \(L\), \(L'\) and \(L_0\) are all perfect, and that
two of the three have \( s = \sigma /2\). Then using the skein exact
triangle, it is possible to show that the same relation holds for the
third. Thus if \(K\) has a skein resolution tree all of whose elements
are perfect, it is necessarily the case that \(s(K) = \sigma (K) /2
\). All the perfect knots we are aware of have such {\it perfect
  resolutions}. It is an interesting question whether all perfect
knots have perfect resolutions, and, more generally, whether all
perfect knots satisfy \( s(K) = \sigma (K) / 2 \).

\subsection{Connected Sums}
Let \( K_i \subset \Nbar _i \) (\(i=1,2\)) be knots, and 
 consider their connected sum
 \(K_1 \# K_2 \subset \Nbar _1 \# \Nbar _2 \).
The complement of \(K_1 \# K_2 \) is obtained from the \(N_i\) by
taking a boundary connected sum along the \(m_i\).
To be precise, we identify a tubular neighborhood
 of \( m_1 \subset \partial N_1 \) with a tubular neighborhood of 
\(  m_2 \subset \partial N_2 \) to produce a new manifold
\(N_1 \#_{m} N_2 \) which is the complement of \(K_1 \# K_2 \).
% (Why is this independent of the identification?) 
The reduced stable complex behaves nicely with
 respect to this operation: 

\begin{prop}
\label{Prop:ConnectSum}
\begin{equation*}
 \cfr (N_1\#_m N_2, m) \cong \cfr (N_1, m_1) \otimes \cfr (N_2, m_2). 
\end{equation*}
\end{prop}

\noindent Specializing to the case of knots in \(S^3\), we get

\begin{cor}
The connected sum of two perfect knots is perfect. 
\end{cor}

\begin{cor}
If \(K_1 \) and \( K_2 \) are knots in \(S^3\), then \(s(K_1 \# K_2) =
s(K_1) + s(K_2) \). 
\end{cor}

\begin{figure}
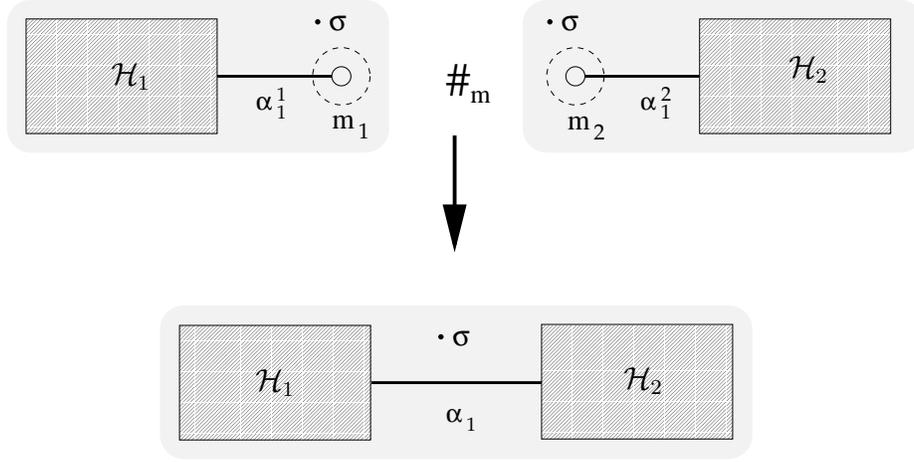
\caption{\label{Fig:KnotSum} Heegaard splittings of \( N_1 \), \(N_2
  \), and \( \Nsum \). }
\end{figure}

The remainder of this section is devoted to the proof of the proposition. 
First, we explain how to get a Heegaard splitting of \( \Nsum \)
from splittings of \(N_1 \) and \(N_2 \). By Lemma~\ref{Lem:Meridian},
we can choose a Heegaard splitting \(\E_1 = (\Sigma^1, \alpha^1,
\beta^1) \)  of
\(N_1 \) with the property that \(m_1 \) has a geometric representative
which intersects \( \alpha _1^1 \) geometrically once and misses all the
other \( \alpha^1_i \)'s. Such a splitting is shown schematically in
Figure~\ref{Fig:KnotSum}. If we choose a  splitting 
\( \E _2 = (\Sigma^2, \alpha^2, \beta^2) \) of \(N_2 \) with the same 
property,
then we can ``sum the Heegaard splittings'' along \( \alpha _1^1 \) and
\( \alpha _1 ^2 \) to obtain a splitting \( \E = 
(\Sigma, \alpha, \beta) \)
 of \( \Nsum \), as indicated in Figure~\ref{Fig:KnotSum}. To be precise,
 \( \Sigma  \) is obtained by cutting \(\Sigma ^i\) open along
 \(m_i \) and then gluing \( \Sigma ^1 \) to \( \Sigma ^2 \) along the
 two resulting circles. It has genus \( g_1 + g_2 -1 \). On \( \Sigma
  \), we join  \( \alpha _1 ^1 \) and \( \alpha _1^2 \) (which have
 been cut open) into a single large attaching circle \( \alpha_1
 \). All the other attaching circles are disjoint from the \(m_i \)
 and remain unchanged.

\begin{figure}
\includegraphics{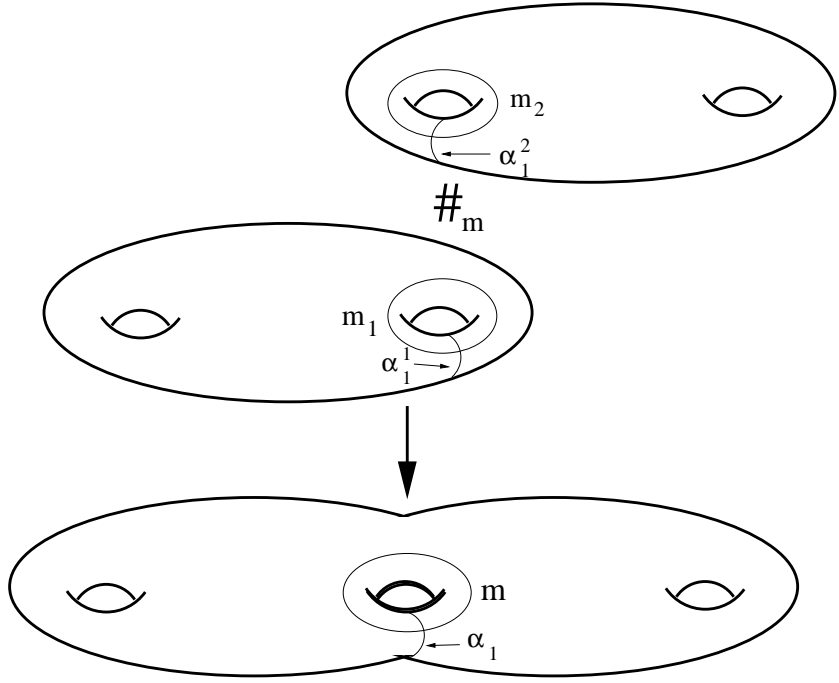}
\caption{\label{Fig:HandleSum} Another view of \( N_1 \#_{m} N_2 \).}
\end{figure}

To see why this procedure works, consider the \( \alpha \)
handlebodies inside  \(N_1\) and \(N_2\). The meridians are circles on
the boundaries of these handlebodies, and by hypothesis there is a
system of compressing disks so that only one disk intersects each
meridian. To take the boundary connect sum along the meridians, we
sandwich the two handlebodies together, as in
Figure~\ref{Fig:HandleSum}. The result is  a new handlebody,
with \( \alpha_1\)  in place of \( \alpha _1^1 \) and \(
\alpha _1^2\). The other attaching circles are  unaffected. 

From this description, it is easy to see that \( \overline{ N_1
  \#_{m} N_2} = \Nbar_1 \# \Nbar _2 \). Indeed, if we attach a
  final handle along \(\emm\), we can cancel it with \( \alpha _1 \)
  to obtain a diagram  which is the connected sum of diagrams for
  \(N_1\) and \(N_2 \). As a consequence,  a \( \spinc \)
  structure on \( N_1 \#_{m} N_2 \) is naturally identified
  with a pair of \( \spinc \) structures on \(N_1\) and \(N_2\). 

Fix \( \spinc \) structures \(\spi _i \) on \(N_i\), and 
let \( E_i \) be a good \(\epsilon\)-class for \((\E_i, \spi_i)\).
 Since \(\emm_1 \)
has only a single intersection with \( \alpha _1^1 \), the elements of
\(E_1 \) are in \(1-1\) correspondence with elements of \( \TTT
_{\hat{\alpha}_1^1} \cap \TTT_{\hat{\beta}^1} \), and similarly for \(E_2
\). Likewise, \(\emm \) has a unique intersection with
\( \alpha _1 \), so the elements in a good \( \epsilon\)-class \(E\) for
\( (\E, \spi _1 \# \spi _2 ) \) correspond to elements of 
\begin{equation*}
 \TTT
_{\hat{\alpha}_1} \cap \TTT_{\hat{\beta}}  \cong 
 (\TTT
_{\hat{\alpha}_1^1} \cap \TTT_{\hat{\beta}^1} ) \times  (\TTT
_{\hat{\alpha}_1^2} \cap \TTT_{\hat{\beta}^2} ).
\end{equation*}
Thus \( E \) is naturally identified with \( E_1 \times E_2 \).

To define the Alexander grading, we orient \(\emm_1\), \(\emm_2\),
 and \(\emm\) in a consistant fashion. Then
the additivity of the Alexander grading implies that 
\begin{equation*}
A(\bfy _1 \times \bfy _2 ) = A (\bfy _1 ) + A (\bfy _2 ).
\end{equation*}

We claim that the homological grading on \( \cfs (\E) \) 
is additive as well. To prove this, it  suffices to show that 
\begin{equation*}
\gr (\bfy _1 \times \bfy_2) - \gr (\bfz _1 \times \bfy_2) =
\gr (\bfy _1) - \gr (\bfz _1) .
\end{equation*}
Let \( \sigma \) be a point in the region separating the two
summands (see Figure~\ref{Fig:KnotSum}), and choose \( \phi \in \pi
_2 ((\bfy_1, \bfz _1)) \) with \( n_{\sigma} (\phi) = 0\). Then \(
\phi \) extends to a disk \( \phi ' \in \pi _2 (\bfy_1 \times \bfy_2, \bfz_1
\times \bfy_2) \) whose domain in the region labeled \(\E_2 \) is
trivial. It is easy to see that that \( \mu
(\phi) = \mu (\phi ') \). Since we chose the orientation of the
meridians to be consistant,  \( n_{z_s}(\phi) = n_{z_s} (\phi ' ) \)
as well. This proves the claim. 
Thus \( \cfs (\E, \spi_1 \# \spi _2 ) \)
 looks like the tensor product complex
 \( \cfs (\E_1, \spi _1) \otimes \cfs (\E_2, \spi _2 ) \) as far as the
 generators and Alexander grading are concerned. In fact, the
 two complexes are isomorphic:

\begin{prop}
 \( \cfs (\E) \cong \cfs (\E_1) \otimes \cfs (\E_2 ) \) as filtered
 complexes. 
\end{prop}

To prove this fact, we borrow an argument from the treatment of
ordinary connected sums in Theorem 12.1 of \cite{OS2}.

\begin{proof}
Consider the Heegaard triple \( ( \Sigma , \alpha , \gamma ,
\bbar ) \), where \( \gamma = \{ \emm, \alpha ^1 _2, \ldots,
\alpha ^1 _{g_1}, \beta ^2 _1, \ldots, \beta ^2 _{g-1} \} \). 
This triple describes a cobordism from \( (\Nbar _1 \# ^{g_2}
(S^1 \times S^2)) \coprod (\Nbar _2 \# ^{g_1-1} (S^1 \times S^2)) \) to \(
\overline{N_1 \#_m N_2} \). In fact, it is not difficult to see that
the Heegaard diagram \( (\Sigma , \alpha , \gamma )\) is
just the Heegaard diagram for \( \Nbar_2 \) 
 connect summed with \( g_1-1\) copies of the genus-one Heegaard
diagram for \( S^1 \times S^2\). Similarly, 
\(  ( \Sigma ,\gamma , \beta ) \) is the Heegaard diagram
obtained by cancelling \( \emm_1 \) with \( \alpha ^1_1 \) in the Heegaard
diagram for \( \Nbar _1 \)  and connect summing
with \( g_2 \) copies of \( S^1 \times S^2 \). 

By counting holomorphic triangles, we
get a chain map 
\begin{align*}
 \Phi : (\cfhat (\Nbar _1, \E_1) \otimes  \cfhat
(\#^{g_2} S^1 \times S^2)) 
\otimes (\cfhat (\Nbar _2, \E_2) & \otimes 
\cfhat ( \# ^{g_1-1} S^1 \times S^2))  \\ & \to \cfhat (
\overline{N_1 \#_m N_2}, \E _1 \oplus \E_2 ). 
\end{align*}
The chain complexes \( \cfhat(\#^{g_2} S^1 \times S^2)) \) and 
\( \cfhat ( \# ^{g_1-1} S^1 \times S^2)) \) are isomorphic to
\(H_*(T^{g_2}) \) and \( H_*(T^{g_1-1})\), with trivial differential;
we let \(\theta _1 \) and \( \theta _2 \) be their top degree
generators. We define a map
\begin{equation*}
 \Phi_1 : \cfhat (\Nbar _1, \E_1)  
\otimes \cfhat (\Nbar _2, \E_2)  \to \cfhat (
\overline{N_1 \#_m N_2}, \E _1 \oplus \E_2 ). 
\end{equation*}
by setting \( \Phi _1 (\bfy _1, \bfy_2) = \Phi (\bfy _1 \times \theta
_1 , \bfy _2 \times \theta _2 ) \). We claim \( \Phi _1 \) is an
isomorphism of filtered chain complexes. 

First, we check that \( \Phi_1 \) is an isomorphism as a map of
groups. We argue (as in the proof of Theorem 12.1 in \cite{OS2}) that 
\( \Phi _1 = \Phi_0 + \alpha \), where \( \Phi _0 (\bfy_1 \otimes
\bfy_2) = \bfy _1 \times \bfy _2 \) and \( \alpha \) is of lower
order in the area filtration. Indeed, it is easy to see that there is
a small holomorphic triangle connecting \( \bfy _1 \times \theta _1
\), \(\bfy _2 \times \theta _2 \), and \( \bfy_1 \times \bfy_2
\). By choosing all the periodic domains to have very small area, we
can ensure that any other holomorphic triangle has larger area than
these small triangles, and thus that \( \alpha = \Phi _1 - \Phi _0 \)
is strictly lower in the area filtration. It follows that \( \Phi _1
\) is an isomorphism of groups, and thus of chain complexes. 

To show that \( \Phi _1 \) is filtration decreasing, we make the usual
argument in which we ``collapse parallel sides,'' and compare the
domain of \( \psi \in \pi _2 (\bfy_1 \times \theta _1 , \bfy_2
\times \theta _2, \bfz_1 \times \bfz_2) \) to the domain of a
corresponding \( \phi \in \pi _2 (\bfy_1 \times \bfy_2 , \bfz_1
\times \bfz_2) \). For \( \psi \) to be a positive domain, \( \phi \)
must be positive as well, which implies that \(A(\bfy_1 \times \bfy
_2 ) \geq A(\bfz_1 \times \bfz_2) \).  
\end{proof}

\noindent Proposition~\ref{Prop:ConnectSum} is an immediate consequence.

%% file: figs/Sum1.pstex_t
\begin{picture}(0,0)%
\includegraphics{figs/Sum1.pstex}%
\end{picture}%
\setlength{\unitlength}{1579sp}%
\begingroup\makeatletter\ifx\SetFigFont\undefined%
\gdef\SetFigFont#1#2#3#4#5{%
  \reset@font\fontsize{#1}{#2pt}%
  \fontfamily{#3}\fontseries{#4}\fontshape{#5}%
  \selectfont}%
\fi\endgroup%
\begin{picture}(14402,7217)(1200,-7877)
\put(2836,-1996){\makebox(0,0)[lb]{\smash{\SetFigFont{12}{14.4}{\familydefault}{\mddefault}{\updefault}{\color[rgb]{0,0,0}\({\mathcal H}_1\)}%
}}}
\put(13480,-1891){\makebox(0,0)[lb]{\smash{\SetFigFont{12}{14.4}{\familydefault}{\mddefault}{\updefault}{\color[rgb]{0,0,0}\( {\mathcal H}_2\)}%
}}}
\put(5105,-6811){\makebox(0,0)[lb]{\smash{\SetFigFont{12}{14.4}{\familydefault}{\mddefault}{\updefault}{\color[rgb]{0,0,0}\({\mathcal H}_1\)}%
}}}
\put(10885,-6781){\makebox(0,0)[lb]{\smash{\SetFigFont{12}{14.4}{\familydefault}{\mddefault}{\updefault}{\color[rgb]{0,0,0}\( {\mathcal H}_2\)}%
}}}
\end{picture}

%% file: triangle.tex
\section{Other Surgeries}
\label{Sec:ExactTriangle}

Until now, we have restricted our attention to the large \(n\)
surgeries on \(N\) ---  manifolds of
the form \(N(n,1) \) with \(n \gg 0 \). In this section, we explain
how to use the exact triangle to derive information about the Floer
homology of other surgeries on \(N\) from these
surgeries. We focus on the case of knots in \(S^3\), but much of what
we do applies to more general three-manifolds. 

Throughout this section, we take all of our homology groups with
coefficients in a field. 
We will sometimes find it more convenient to work with \( \hfm \)
rather than \( \hfp \). To get from one to the other, we use the
following elementary fact: 

\begin{lem}
If \(Y\) is a rational homology sphere, then \( \hfp (Y, \spi) \cong
\hfred (Y, \spi) \oplus \Qu \) and \( \hfm (Y, \spi) \cong
\hfred (Y, \spi) \oplus \Q[u] \). The grading of \( \hfred \subset
\hfm \) is one less than the grading of \( \hfred \subset \hfp \), and
the grading of \( 1 \in \Q[u] \) is two less than the grading of \( 1
\in \Qu \). 
\end{lem}

\begin{proof}
The long exact sequence 
\begin{equation*}
\begin{CD}
@>>> \hfm (Y, \spi) @>>> \hfi (Y, \spi) @>>> \hfp (Y, \spi) @>>> 
\end{CD}
\end{equation*}
splits to give short exact sequences
\begin{equation*}
\begin{CD}
0 @>>> \hfred (Y, \spi) @>>> \hfm (Y, \spi) @>>> \Q[u] @>>> 0 
\end{CD}
\end{equation*}
and 
\begin{equation*}
\begin{CD}
0 @>>> \Qu @>>> \hfp (Y, \spi) @>>> \hfred (Y, \spi) @>>> 0 
\end{CD}
\end{equation*}
It is easy to see that these sequences split.
\end{proof}

\subsection{Exact triangles}
The exact triangles of \cite{OS2} provide a powerful tool for
computing the Ozsvath-Szabo Floer homology.
 The formulation which is most useful to
us is described in  Theorem 10.19 of \cite{OS2}. For the reader's
convenience, we recall the statement here:

\begin{prop} (Theorem 10.19 of \cite{OS2}).
\label{Prop:ExactTriangles}
For \(n > 0 \), there are long exact sequences
\begin{equation*}
\begin{CD}
@>g_k^+>> \hfp (\Nbar, \spi) @>{\oplus f_i^+}>> \bigoplus _{i \equiv k (n)}
 \hfp (N(0,1), \spi_i) @>>> \hfp(N(n,1), \spi _k) @>>> 
\end{CD}
\end{equation*}
and 
\begin{equation*}
\begin{CD}
 @>>> \hfp(N(-n,1), \spi _k) @>>> \bigoplus _{i \equiv k (n)}
 \hfp (N(0,1), \spi_i) @>>> \hfp (\Nbar, \spi)@>>>
\end{CD}
\end{equation*}
The maps \(g_k^+\) and \(f^+_i\) (and their analogues in the second triangle) are induced by the appropriate surgery cobordisms. 
%(True even for the third map?)
\end{prop}

The exact triangle is  very helpful in studying the map a cobordism
 induces on \( \hfp\). 
In many situations, however,  we are more interested in the map
induced on \( \hfm \). For example, if we start with \( B^4\) 
and attach a zero-framed
two-handle along a knot \(K \), the relative invariant of the
resulting four-manifold is determined by the induced map on \( \hfm
\). To facilitate such computations, it would be nice if there were
exact triangle for \( \hfm \) as well as for \( \hfp
\). Morally, such a triangle should be the ``dual'' of the exact
 triangle for \( \hfp \). Since things are infinitely generated,
 however, we need to be careful. For example, if we try to replace the
 map  \( \hfp (N,1, \spi_k ) \to \hfp (\Nbar, \spi) \) with an
 analogous map on \( \hfm\), we get something with infinitely many
 nonzero terms! To get around this problem, we must work with a power series
 ring. Let
\( \hfmt (\Nbar, \spi) = \hfm (\Nbar, \spi) \otimes _{\Q[u] } \Q((u))
 \). 

\begin{prop}
\label{Prop:MinusTriangle}
Suppose \(n > 0 \) and that \(n \) does not divide \(k\).
Then there is a long exact sequence 
\begin{equation*}
\begin{CD}
@>>> \hfmt (\Nbar,\spi) @>{\oplus f_i^-}>> \bigoplus _{i \equiv k (n)}
 \hfred (N(0,1), \spi _i) @>>> \hfmt(N(n,1), \spi _k)@>>>
\end{CD}
\end{equation*}
whose maps are induced by the appropriate surgery cobordisms.
\end{prop}

\begin{proof}
We begin with the exact triangle for \(-n\) surgery on \(-N\): 
\begin{equation*}
\begin{CD}
 @>>> \hfp(-N(-n,1), \spi _k) @>>> \bigoplus _{i \equiv k (n)}
 \hfp (-N(0,1), \spi_i) @>>> \hfp (-\Nbar, \spi)
\end{CD}
\end{equation*}
Since the dual of an exact sequence of vector spaces is exact, we get
an exact sequence 
\begin{equation*}
\begin{CD}
 (\hfp (-\Nbar, \spi ))^* @<<< (\hfp(-N(-n,1), \spi _k))^* @<<<
 \bigoplus _{i \equiv k (n)}
 (\hfp (-N(0,1), \spi_i))^* .
\end{CD}
\end{equation*}
At the level of complexes, we have dualities
\begin{align*}
 (\cfp (-\Nbar, \spi ))^* & \cong \cfm(\Nbar, \spi) \\
(\cfp(-N(-n,1), \spi _k))^* & \cong \cfm(N(n,1), \spi _k) \\
(\cfp (-N(0,1), \spi_i))^* & \cong \cfm (N(0,1), \spi_i).
\end{align*}
All of these complexes are infinitely generated, so we cannot apply
the universal coefficient theorem directly. Instead, we exhaust
them by finite subcomplexes. Let \(C\) denote any of the \( \cfp \)
complexes above. We consider the increasing sequence of subcomplexes 
\(C_0 \subset C_1 \subset C_2 \ldots \), where \( C_n \) is generated
by \( \{ [\bfy, i] \ts | \ts i \leq n \}\). 
 Let \( \ldots \to Q_2 \to Q_1
\to Q_0 \) be the dual sequence of quotient complexes of the dual
complex \(Q\). Then by the
universal coefficient theorem, \( (H_*(C_i))^* \cong H_*(Q_i) \). 
Since \(C = \dlim C_i \), \(H_*(C) = \dlim H_*(C_i)\). Thus
\begin{equation*}
(H_*(C))^* \cong (\dlim H_*(C_i))^* \cong \ilim (H_*(C_i)^*) \cong 
\ilim H_*(Q_i)
\end{equation*}
Now when \(C \) is integer graded, as is the case for \( \Nbar \) or
\( N(n,1) \), the filtration \(Q_i \) has the property that all the
generators of a given degree are contained in some fixed
\(Q_i\). Using this fact, it is not difficult to see that 
\( \ilim H_*(Q_i) \cong \ilim D_i \), where \(D_i = \{ x \in H_*(Q)
\ts | \ts
\gr (x) \geq i \} \). Under this process, the inverse limit of a
torsion module \( \tmod _n\) is again \(\tmod _n \), while the inverse
limit of a free module \( \Q[u] \) is the power series ring \(\Q((u))
\). Thus \((H_*(C))^* \cong H_*(Q) \otimes_{\Q[u]}\Q((u))\). 
This identifies the terms in the exact sequence involving \(N(n,1)\)
and \(\Nbar \) with the corresponding terms in the statement.

\begin{lem}
For \(i \neq 0 \), \( \hfred (N(0,1), \spi _i) \cong \hfp (N(0,1),
\spi _i)\) and \( \hfm(N(0,1), \spi_i) \cong \hfred (N(0,1), \spi
_i) \oplus \hfi (N(0,1))\).
\end{lem}

\begin{proof}
For \(i \neq 0 \), \(\spi_i \) is a nontorsion \( \spinc \) structure. By
Proposition 11.5 of \cite{OS2} and the subsequent discussion,
 we know that \( \hfi ((N(0,1), \spi _i)) \)
 is a sum of modules of the form  \( \Q[u]/(u^m-1) \) for some \(m\).
 Since every element of 
\(\hfp (N(0,1), \spi _i)\) is killed by some power of \(u\), the map
\( \pi \) in the exact sequence 
\begin{equation*}
\begin{CD}
@>>>  \hfm (N(0,1), \spi _i) @>>>  \hfi (N(0,1), \spi _i) @>\pi>> 
 \hfp (N(0,1), \spi _i) @>>>
\end{CD}
\end{equation*}
must be trivial. 

To prove the second claim, we use the structure theorem for modules
over a PID to write \( \hfm (N(0,1), \spi _i) \cong \oplus T_i \),
where \(T_i \cong  Q[u] / (p_i) \). \(\hfred ((N(0,1), \spi _i)) \)
is the sum of those \(T_i\) for which \(p_i = u^{d_i} \), and \( \hfi
(N(0,1)) \) is the rest. 
\end{proof}

\begin{lem}
For \(i \neq 0 \), 
\((\hfred (-N(0,1), \spi _i))^* \cong \hfred (N(0,1), \spi _i) \).
\end{lem}

\begin{proof}
We still have \( \hfp (N(0,1), \spi _i) = \dlim H_*(C_i) \cong \dlim
 (H_*(Q_i))^* \). Now \(  \dlim (H_*(Q_i))^*  \) has an obvious map
 into \( (H_*(Q))^* \), and it is easy to see that this map is
 injective. Thus
 \( \hfp (N(0,1), \spi _i)\) injects into \( (\hfm (-(N(0,1), \spi_i))^*\). 

Suppose 
\( x \in \hfp (N(0,1), \spi _i) \), while \( y \in \hfi (-N(0,1), \spi _i) \subset \hfm (-N(0,1), \spi_i)\). Then 
\begin{equation*}
 \langle x, y \rangle = \langle x, u^{ni} y \rangle = \langle u^{ni} x, y \rangle = 0 
\end{equation*}
since any \(x\) is killed by a high enough power of \(u\). Thus we actually have 
\begin{equation*}
\hfred (N(0,1), \spi _i) = \hfp (N(0,1), \spi _i) \subset (\hfred (-(N(0,1), \spi_i))^*. 
\end{equation*}
Applying the same argument to \(- N \), we see that 
\begin{equation*}
\hfred (-N(0,1), \spi _i)  \subset (\hfred ((N(0,1), \spi_i))^*. 
\end{equation*}
which proves the lemma. 
\end{proof}

Since \( n \) does not divide \( k \), \(\spi _0 \) does not appear
 in the direct
sum. Thus the lemmas imply that 
\begin{equation*}
\bigoplus _{i \equiv k (n)} (\hfp (-N(0,1), \spi_i))^* \cong
\bigoplus _{i \equiv k (n)} \hfred (N(0,1), \spi_i).
\end{equation*}
This gives the exact sequence of the statement.

It remains to check that the maps in this sequence are
 induced by the surgery cobordisms. Consider the map 
\(g_k^-: \hfm (N(n,1), \spi _k) \to \hfm (\Nbar) \)
 induced by the cobordism \(W\) from \( N(n,1) \) to \(\Nbar \). Now the
 corresponding map in the  triangle  we dualized was
 induced by \(-W\). Thus  \(g_k^-\) is dual to
 \(g_k^+ \) at the chain level, and the 
 induced maps on homology are dual as well (we are using field coefficients). Since the map in the exact sequence  was
 constructed as the dual of \(g_k^+\),  it is \(g_k^- \) . A similar argument shows that  the map \( \hfm (\Nbar,\spi) \to \oplus _{i \equiv k (n)} \hfred (N(0,1), \spi _i)\) is given by \( \oplus f_i^- \).

%The hypothesis that \(n \) does not divide \(k\) implies that the \(
%\spi _i \)'s in the direct sum are all nontorsion. It follows
%from ??? of \cite{OS4} that \( \im f_i \subset \hfred(N(0,1)) \).
\end{proof}

\subsection{Integer surgeries}
Suppose \(K\) is a knot in \(S^3\). Let us try to compute \(\hfp
(K(m,1), \spi_k) \) for an arbitrary integer \(m\). To keep things
simple, we assume for the moment that \(k \) does not divide \(m\). 

We start by studying the exact triangles  of
Propositions~\ref{Prop:ExactTriangles} and \ref{Prop:MinusTriangle} 
for \( n \gg 0 \). Since \( \hfp (K(0,1)) \) is supported in a finite
number of \( \spinc \) structures, we
can take \(n\) to be so large that there is only one non-zero term in the
direct sum. Substituting the known homology of
 \( K(1,0) = S^3 \), the sequences reduce to 
\begin{equation*}
\begin{CD}
@>>> \Qu @>0>> \hfred(K(0,1), \spi_k) @>>> \hfp (K(n,1), \spi _k)
 @>g_k^+>> \Qu
\end{CD}
\end{equation*}
and
\begin{equation*}
\label{Eq:MinusSeq}
\begin{CD}
@>>> \Q((u)) @>f_k^->> \hfred(K(0,1), \spi _k)
 @>>> \hfmt (K(n,1), \spi_k) @>g_k^->> \Q((u)) 
\end{CD}
\end{equation*}

Recall that \( \hfmt(K(n,1)) \cong \Q((u)) \oplus \hfred(K(n,1))
 \). The map \(g^- _k\) takes  the \( \Q((u)) \)  summand injectively 
to \( \hfmt (S^3) \) and is trivial on \( \hfred \). Thus it may be
 described by a single number: the rank of its
 cokernel.

\begin{defn}
If \(K\) is a knot in \(S^3\), the local \(h\)-invariant \(h_k(K)\) is
the rank of \( \coker g_k^-  \cong \im f_k^- \).
\end{defn}

We have chosen the phrase ``local {\it h}-invariant'' because of the
close  connection between this number and Fr{\o}yshov's {\it
  h}-invariant in Seiberg-Witten theory  \cite{FrLect}. The
analogous object in \OS theory is the ``correction term'' \(d\)
defined in \cite{OS4}.

We make a few elementary observations regarding the \(h_k\). First, 
the maps \( f_k^- \) and \( f_{-k}^- \) are conjugate symmetric, 
  so  \( h_k(K) = h_{-k}(K) \). 
Second, note that we  could also have defined \( h_k \) in terms of 
\( \hfp \).
Indeed \( \hfp (K(n,1), \spi_k) \cong  \Qu \oplus \hfred(K(n,1),
\spi_k) \), and \(g_k^+\) maps the  \( \Qu \) summand  onto \( \hfp (S^3) \). 
It is not difficult to see
 (either by counting Euler characteristics or by looking at the
 gradings) that \(h_k \) is the rank of the kernel of the map \( \Qu
 \to \hfp(S^3) \). 

Together, the \(h_k\)'s and the groups \( \hfp (K(n,1), \spi _k) \)
for \(n \gg 0 \) encapsulate almost  everything there is to know about 
positive integer surgeries on \(K\). To see this, we return to the 
 ``\(-\)'' exact triangle for large \(n\), which splits into a short exact
 sequence
\begin{equation*}
\begin{CD}
0@>>> \tmod _{h_k} @>>> \hfred (K(0,1), \spi _k) @>>> \hfred (K(n,1),
\spi _k) @>>> 0 
\end{CD}
\end{equation*}
This enables us to recover the Betti numbers of \( \hfred (K(0,1),
\spi _k)\), although there is still some ambiguity about the \(\Q[u]\)
module structure. 

Next, we consider  the ``\(-\)''  exact triangle
for arbitary values of \(m>0\):
\begin{equation*}
\begin{CD}
@>>> \hfm (S^3) @>{\oplus f_i^-}>> \bigoplus _{i \equiv k (m)}
 \hfred (K(0,1), \spi _i) @>>> \hfm(K(m,1), \spi _k) @>{g^-_{m,k}}>> 
\end{CD}
\end{equation*}
 The form of this sequence is  determined by the single number 
\begin{equation*}
 h_{m,k} = \rank \coker g_{m,k}^- = \rank \im (\oplus f_i^-).
\end{equation*} 
The maps \( f^-_i \) are induced by the surgery cobordism from \(
S^3 \) to \(K(0,1) \), so they are independent of \(m\). It follows
that 
 \( h_{m,k} =
\max \{ h_i(K) \ts | \ts i \equiv k \mod (m) \} \), and there is   another short
exact sequence 
\begin{equation*}
\begin{CD}
0@>>>\tmod _{h_{m,k}} @>>> \oplus _{i\equiv k (m)}
\hfm (K(0,1), \spi_i) @>>>  \hfred (K(m,1), \spi _k)  @>>> 0
\end{CD}
\end{equation*}
which enables us to compute the Betti numbers of \( (K(m,1), \spi _k) \). (Again, there may be some ambiguity in the structure as a \( \Z[u] \) module. When \(K\) is perfect, both exact sequences split for grading reasons, so the \(u\) action is completely determined. See \cite{2bridge}.)
It is also easy to find the correction term \(d\) for 
\( K(n,1) \): since the
exact triangle respects the absolute grading, we have \( d(K(n,1),
\spi _k) = d(L(n,1), \spi _k) - 2h_{n,k} \). 

When \(m<0\), the Floer homology of \( K(m,1) \)
can be computed in a similar way, using the second exact sequence of
Proposition~\ref{Prop:ExactTriangles}. In this case, the 
\(h_k (K)\)'s are replaced by  \( h_k^- (K) = h_k
(-K) \).

The same strategy can be used  in the case
of a more general three-manifold \(N\). Things are more complicated,
however, because the image of \( \hfm (N, \spi_k) \) in \( \hfm
(\Nbar, \spibar) \) can no longer be specified by a single number. 

%\vskip0.1in
%\noindent{\bf Example:} Let \(K\) be the right-handed \((2,n) \) torus
%knot. Then 
%\begin{equation*}
% \hfred (K(0,1), \spi _k) \cong \tmod _{\lceil n/2 - |k|
%  \rceil}
%\end{equation*}
%The map \(f_k^- \) is surjective, so \( h_k(K) =\lceil n/2 - |k|
%  \rceil \) and   \( \hfm (K(n,1), \spi _k)
%\cong \Z[u] \) when \(n \gg 0 \). For the left-handed torus knot
%\(-K\), we have \(\hfp (-K(0,1), \spi _k)  \cong \hfp (K(0,1), \spi
%_k) \), but the map \( f_k^- : \hfm (S^3) \to \hfp (-K(0,1), \spi _k)
%\) is trivial, so \( \hfm (-K(n,1)) \cong \Z [u] \oplus \tmod _{\lceil
%  n/2 - |k| \rceil} \).  (Comment on Sein Structures?)

\subsection{Torsion \(\spinc \) structures}
Up to this point,
 we have avoided the case when \(k\) divides \(m\)
because it involves the torsion \( \spinc \) structure \( \spi _0 \)
on \( K(0,1) \). 
 To treat this \( \spinc
\) structure  on an equal footing with the others, we must use twisted
coefficients. We digress briefly to describe some basic facts about twisted coefficients for knot complements. 

\begin{prop}
\label{Prop:TwistedKnot}
If \(K\) is a knot in \(S^3\), 
\begin{equation*}
\underline{\hfp}(K(0,1), \spi _k ) \cong
\hfp(K(0,1), \spi _k ) \otimes \Qt 
\end{equation*}
if \(  k \neq 0 \). For \( k = 0 \) 
\begin{equation*}
\underline{\hfp}(K(0,1), \spi _k ) \cong 
\underline{\hfred}(K(0,1), \spi _k )  \oplus \Qu 
\end{equation*}
and there is a \( \Q \)-module \(A_0 \) such that 
\begin{equation*}
\underline{\hfred}(K(0,1), \spi _k ) \cong A_0 \otimes \Qt. 
\end{equation*}
\end{prop}

\noindent {\bf Remark:}
The module \(A_0 \) is always finitely generated. It is
 a quotient of \( \hfp (K(0,1), \spi _k)
 \) and has \( \hfred (K(0,1), \spi_0) \) as a quotient. In many respects,
\(A_0 \) is a more natural object of study than the latter group ---
 for example, its Euler characteristic is given by the same formula
 as the Euler characteristics of \( \hfp (K(0,1), \spi _k) \) for \( k
 \neq 0 \). The Seiberg-Witten analogue  of \(A_0 \) is the group
 obtained by making a small nonexact perturbation of the Seiberg-Witten
 equations to completely eliminate the reducibles.

\begin{proof} 
We use the exact triangle for twisted coefficients (Theorem 10.23 of
\cite{OS2}), which gives us a long exact sequence 
\begin{equation*}
\begin{CD}
@>g>> {\hfp} (S^3) \otimes \Qt
@>{ f}>> 
 \underline{\hfp} (K(0,1), \spi_k) @>>>
 {\hfp}(K(n,1), \spi _k) \otimes \Qt
\end{CD}
\end{equation*}
When \( k \neq 0 \), \(g \) is a surjection, so 
\( \underline{\hfp} (K(0,1), \spi_k) \cong \ker g \) is a free 
 \( \Qt \) module. To put it another (more complicated) way, 
\( \underline{\hfp} (K(0,1), \spi_k) \) is the homology of the complex
\begin{equation*}
\begin{CD}
0 @>>> {\hfp}(K(n,1), \spi _k) \otimes \Qt @>g>> {\hfp} (S^3) \otimes
\Qt @>>> 0 . 
\end{CD}
\end{equation*}
The advantage of this point of view is that if we substitute \(T=1\),
we get a complex whose homology is \( \hfp (K(0,1), \spi _k )
\). Applying the universal coefficient theorem, we see that
\(\hfp (K(0,1), \spi _k) \cong \underline{\hfp} (K(0,1), \spi_k)
\otimes _{\Qt} \Q \). To arrive at the statement of the proposition, we
tensor both sides with \( \Qt \) and use the identity
\begin{equation*}
 (A \otimes _B C) \otimes _C B  \cong A \otimes _B (C \otimes _C B)
 \cong A \otimes _B B  \cong A.
\end{equation*}

If \(k = 0 \), \(g \) is no longer onto, and we have a short exact sequence
\begin{equation*}
\begin{CD}
0@>>> \coker g @>>> \underline{\hfp} (K(0,1), \spi_k) @>>> \ker
g @>>>0.
\end{CD}
\end{equation*}
We claim that \(  \coker g \) is the submodule \(N_{\infty} \) of
\( \underline{\hfp} (K(0,1), \spi_k) \) which pushes down from
infinity. Indeed, Theorem 11.12 of \cite{OS2} tells us that 
\( \underline{\hfi} (K(0,1), \spi_k) \) is a \((T-1)\)-torsion module,
so \( N_{\infty} \) is a \((T-1)\)-torsion module as well. 
But since \( \ker g \) is a submodule of a free  \( \Qt \) module, it is
free. Thus the image of \(N_{\infty} \) in \( \ker g \) is trivial, which
implies  \( N_{\infty} \subset \coker g \). On the other hand, \(
\coker g \) is a quotient of 
\begin{equation*}
\Qu \otimes \Qt / g(\Qu \otimes \Qt) \cong \Qu
\end{equation*}
so everything in \( \coker g \) pushes down from infinity. 

Since \(  \ker g \cong \underline{\hfred} (K(0,1), \spi_k)  \) is
free, we can split the exact sequence to get the stated direct sum
decomposition. Finally, we set \(A_0 = \underline{\hfred} (K(0,1),
\spi_k) \otimes _{\Qt} \Q \) and use the same trick as above to see that 
\( \underline{\hfred} (K(0,1), \spi_k) \cong A_0 \otimes \Qt \).
\end{proof}

To define \(h_0 (K) \), we consider the map 
\begin{equation*}
g_k \: \hfp (S^3) \otimes \Zt \to \hfp (-K(-n,1), \spi _k ) \otimes \Zt.
\end{equation*}
in the twisted exact triangle for large \(-n\)-surgery on \(-K\).
In analogy with the untwisted case, we  
 define \(h_k (K) \) to be the rank of \( \ker g_k \) over \( \Qt \). 
Of course, before we do so, we should check that 
 this is consistent with our former definition when \( k \neq 0\).

\begin{prop}
For \(k \neq 0 \), the rank of \( \ker g_k \) is equal to \(h_k (K) \). 
\end{prop}

\begin{proof}
As in the untwisted case, 
\begin{equation*}
 \chi (\underline{\hfp }(K(0,1), \spi _k))
  = \chi (\hfred (-K(-n,1), \spi _k)) + \rank \ker g_k.
\end{equation*}
Since \( \chi (\underline{\hfp }(K(0,1), \spi _k)) = \chi
 (\hfp (K(0,1), \spi _k) \),  the ranks of the kernels of
  the twisted and untwisted maps must be equal.
\end{proof}

With these facts in place, we can repeat the calculations
of the preceding section using twisted coefficients. The net result is
a short exact sequence 
\begin{equation*}
\begin{CD}
0 @>>> \hfred (-K(-m,1), \spi _k) @>>> \bigoplus _{i
  \equiv k \ (m)} A_i   @>>> \tmod_{h_{m,k}}  @>>>0
\end{CD}
\end{equation*}
from which we can easily deduce the Betti numbers of 
\( \hfred (K(m,1), \spi _k) \), regardless of whether \(k\) divides  \(m\). 
In addition, we see that the formula  \(d(K(m,1), \spi _k ) = 
d(L(m,1), \spi _k) - 2h_{m,k} \) holds for any value of \(k\). 

%The number \(h_0 \) can also be described in terms of \( \hfp (K(0,1),
%\spi _0 ) \).
%
%\begin{prop}
%\begin{equation*}
%\hfp (K(0,1), \spi _0 ) \cong  \hfred (K, \spi _0) \oplus \Zu \oplus \Zu
%\end{equation*}
%If \(1-o \) and \( 1_e \) are in the summands with even and odd
%absolute grading, then \( \gr 1_e - \gr 1_o = 2 h_0 \).
%\end{prop}
%
%\begin{proof}
%??????
%\end{proof}
%
%
%\subsection{Fractional surgeries}
%We can also use twisted coefficients to compute the Floer homology of
%\( 1/q \) surgery  on \(K\). To do so, we use the long exact
%sequence for fractional surgeries (Theorem 10.??? of \cite{OS2})
%together with the following 
%
%\begin{lem}
%\begin{align*}
%\hfp (K(0,1), \spi _k; \Z/q) \cong 
%\begin{cases}
%\hfp (K(0,1), \spi _k ) \otimes \Z [\Z/q] & \text{if} \  k\neq 0 \\
%A_0 \otimes \Z [\Z/q] \oplus H_*(S^1)
%\otimes\Zu  & \text{if} \  k\neq 0
%\end{cases}
%\end{align*}
%\end{lem}
%
%\begin{proof}
%This follows from Proposition~\ref{Prop:TwistedKnot} and the universal
%coefficient theorem. FALSE AS U MODULES!
%\end{proof}
%
%\begin{prop}
%There is a short exact sequence ????????????
%\end{prop}
%
%In general, we expect that it should be possible to compute the
%\OS Floer homology of any surgery on \(K\) in terms of the \( \hfp
%(K(n,1), \spi _k ) \) for \(k \geq 0 \) and the \( h_k \)'s. We hope
%to discuss this matter in a future expository paper. 
%
\subsection{A model for the exact triangle}
\label{Subsec:Model}
Of course, there is an important gap in the 
preceding discussion. How can we compute the \( h_k (K) \)'s? 
Although we are not able to provide a complete answer to this
question, we can reduce it to a more tractable problem about the complex \(
\cfp (K, \spi _k) \).
Recall from Corollary~\ref{Cor:Seq1} that there is a short exact
sequence 
\begin{equation*}
\begin{CD}
 0 @>>> C_{\spi_k}(\E) @>>> \cfp (\E,\spi_k) @>>> \cfps (\E) @>>> 0
\end{CD}
\end{equation*}
where \( C_{\spi_k} \) is the complex shown in Figure~\ref{Fig:C_k}.
It is interesting to compare the resulting long exact sequence
\begin{equation*}
\label{Eq:TriangleModel}
\begin{CD}
@>>>H_*(C_{\spi_k}) @>>> \hfp (N, \spi _k) @>f>> 
\hfp (\Nbar, \spibar) @>>>
\end{CD}
\end{equation*}
with the exact triangle for \(n\gg 0 \). This comparison is
 particularly easy when \(N\) is a knot complement. 
 Since \( C_{\spi _k} \) is finite, 
 \(f\) must map the \(\Qu \) summand of \(\hfp (K, \spi
 _k) \) onto \( \hfp (S^3) \cong \Qu \).  Let \(c _k \) be the rank of
 the kernel of this map.

\begin{prop}
\label{Prop:Seq=Tri}
For \(k \geq 0 \),  \( c_k = h_k (K) \) .
\end{prop}

\begin{proof}
We have short exact sequences 
\begin{equation*}
\begin{CD}
0 @>>> \tmod _{c_k} @>>> H_*(C_{\spi_k}) @>>> \hfred (K, \spi_k) @>>> 0 
\end{CD}
\end{equation*}
and for \(k \neq 0 \)
\begin{equation*}
\begin{CD}
0 @>>> \tmod _{h_k(K)} @>>> \hfp (K(0,1), \spi _k)
 @>>> \hfred (K, \spi_k) @>>> 0 .
\end{CD}
\end{equation*}
For \(k > 0 \), the Euler characteristic of \(C_{\spi_k}\) is  
 the Euler
characteristic of \(\hfp (K(0,1), \spi _k)\). Indeed, 
\(\chi  (\hfp (K(0,1), \spi _k)) \) was computed in Theorem 9.1 of
\cite{OS2} by
identifying it with \( \chi (C_{\spi_k}) \). Thus we must have \(c_k=
h_k(K) \).

 For \(k= 0 \), there is a similar
argument using twisted coefficients. By
Proposition~\ref{Prop:TwistedKnot} and the subsequent discussion,
 we still have 
\begin{equation*}
 \chi (\underline{\hfp} (K(0,1), \spi _0)) = \chi (\tmod _{h_0}) +
\chi ( \hfred (K, \spi _0)) .
\end{equation*}
 In addition, the filtration argument
used to prove Theorem 9.1 of \cite{OS2} can be extended to the case 
\(k=0 \) by using 
twisted coefficients. Thus \( \chi (\underline{\hfp} (K(0,1), \spi _0))
= \chi (C_{\spi_0}) \), which proves the claim for \(k=0 \) as well. 
\end{proof}
 \noindent Note that the statement of the proposition is definitely
{ false} when \( k < 0 \).

\begin{cor}
\label{Cor:hFormula}
Denote the element of lowest degree in  the \( \Qu \)
summand of \( \hfp (K, \spi_k) \) by \(1_k\). Then for 
\( k \geq 0 \), \( h_k (K) = (\gr (1_s) - \gr (1_k))/2 \). (Both gradings are taken in 
\( \cfi (K) \).) 
\end{cor}

\noindent Combined with  the calculations of section~\ref{Sec:Quasiperfect},
this implies

\begin{cor}
If \(K\) is perfect, \( h_k (K) = \max \{ \lceil s(K) - k \rceil / 2,
0 \} \). 
\end{cor}

Proposition~\ref{Prop:Seq=Tri} is very useful to us. Rather than 
studying the map 
\begin{equation*}
 \hfp (S^3) \to \hfp (K(0,1), \spi_i),
\end{equation*}
 which is a complicated problem 
involving holomorphic triangles, we can think about the map \( \hfps
(K)  \to H_*(C_{\spi _i}) \). This map is defined purely in terms of
the  stable complex, and we can often use formal properties of the
Alexander grading to prove things about it. (See the next section for
 some examples.) 

We think of the long exact sequence
 derived from the Alexander filtration  as a {\it model} for the
exact triangle for surgery on \(K\).  The maps in the triangle are not
the identical to the maps in the long exact sequence, but the fact that they have isomorphic kernels and cokernels is good often enough. 
 There are similar 
 sequences which serve as models for each of the other three
 large-\(n\) surgery triangles associated to a knot. These are all
 derived in an obvious fashion from the one we have considered, so we
 leave it to the reader to work out the details.

It would be useful  to know if there are other, more general
 situations in which the long exact sequence derived from the
 Alexander filtration is a model
 for the exact triangle. We will return to this question at the end of
 the section.

\subsection{Properties of \(h_k(K)\)}
We collect some basic facts about  the local \(h\) invariants here.

\begin{prop}
\label{Prop:hBounds}
For \( k > 0 \), we have \( h_k(K) +1 \geq h_{k-1}(K) \geq h_k(K) \).
\end{prop}

\begin{figure}
\includegraphics{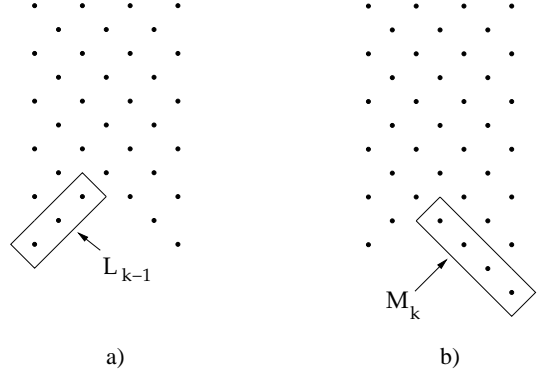}
\caption{\label{Fig:Seq2} The subcomplexes \( L_{k-1} \) and \(M_k\)
  of Proposition~\ref{Prop:hBounds}}
\end{figure}

\begin{proof}
There is a short exact sequence 
\begin{equation*}
\begin{CD}
0 @>>> L_{k-1} @>>> \cfp(K, \spi _{k-1}) @>>> \cfp (K, \spi _k) @>>> 0 
\end{CD}
\end{equation*}
as shown in Figure~\ref{Fig:Seq2}a. In the associated long exact sequence,
 the \( \Qu \) summand of \(\hfp(K, \spi _{k-1}) \) maps onto the \( \Qu \)
 summand of  \( \hfp (K, \spi _k) \), so \( \gr \ 1_{k-1} \leq \gr \ 1_{k} \).
By Corollary~\ref{Cor:hFormula}, \( h_{k-1}(K) \geq h_k(K)
\). For the other direction, we use the short exact sequence
\begin{equation*}
\begin{CD}
0 @>>> M_{k} @>>> \cfp(K, \spi _{k}) @>>> \cfp (K, \spi _{k-1})\{1\} @>>> 0 
\end{CD}
\end{equation*}
of  Figure~\ref{Fig:Seq2}b. \(\cfp (K, \spi _{k-1})\{1\} \) denotes
 the complex \(\cfp (K, \spi _{k-1}) \) ``shifted up by 1'' inside of
 \(\cfi (K) \). 
 This time the \( \Qu \) summand of \(\hfp(K, \spi _{k}) \)   maps onto
the \( \Qu \) summand of \( \hfp (K, \spi _{k-1})\{1\} \), but
 because of the shift in degree we only get 
 \( h_k(K) +1 \geq h_{k-1} \). 
\end{proof}

\begin{cor}
\( h_{n,k} (K) = h_{k_0} (K) \), where \(k_0 \) is the number
congruent to \(k \mod n \) with smallest absolute value.
\end{cor}

\begin{cor}
\label{Cor:hGenusBound}
\( h_k (K) \leq g^*(K) - k \), where \(g^*(K) \) is the slice genus of
\(K\). 
\end{cor}

\begin{proof}
By the adjunction inequality,
 the map \( f_k : \hfm (S^3) \to \hfred (K(0,1),
\spi _k) \) is trivial for \( k \geq g^*(K) \). Thus \(h_{g^*(K)} =
0\). 
\end{proof}

It would be interesting to know what values of \(h_k \) can be
achieved by actual knots. The preceding results place some basic
constraints on these values. Unfortunately, the number of cases in
which we know how to compute \(h_k (K) \) is small. Beyond the
perfect knots, for which the \(h_k\)'s are completely determined by the
number  \(s(K) \), the only real class of examples for which we know
 \(h_k \) is knots with lens space surgeries. Ozsvath and Szabo have
 computed \(h_k\) for such knots in \cite{OS2}. They show that if
 \(K\) has a positive lens space surgery (like a right-hand torus
 knot), then 
\begin{equation*}
h_k(K) = \sum _{i=1}^\infty i a_{k+i}
\end{equation*}  
where the \(a_i \) are the coefficients of the Alexander polynomial. 
In other words, \(h_k \) is the Seiberg-Witten invariant associated
to the \(k\)th \(\spinc \) structure for 0-surgery on \(K\).
If \(K\) has a negative lens space surgery (like a left-hand torus
knot), then the \(h_k(K)\)'s are all zero. 

Although it would be unwise to conjecture anything based on so little
evidence, these examples raise some interesting questions. 
First, the bound provided by Corollary~\ref{Cor:hGenusBound} is probably not
sharp. Indeed, Fr{\o}yshov has shown in \cite{Froyshov} that his
Yang-Mills version of the \(h\)-invariant satisfies \( h (K) \leq
\lceil g^*(K)/2 \rceil \). In light of this fact, 
it seems likely that the perfect knots have the
 ``maximum possible
growth rate'' for the \(h_k\)'s, in the sense that  there should be an
inequality  \(h_k(K) \leq \lceil \frac{g^*(K) - k}{2} \rceil \).
It is
definitely not true, however, that \( h_{k-2}(K) \leq h_{k}(K) + 1
\) --- there are plenty of examples of torus knots for which \( h_k
\) increases twice in succession as \(k\) decreases.

There is also a relation between \(h_k(K) \) and the invariant \(s(K)
\): 

\begin{figure}
\includegraphics{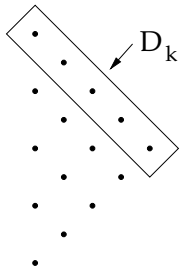}
\caption{\label{Fig:D_k} The quotient complex \(D_k\).}
\end{figure}

\begin{prop}
\label{Prop:s=Top}
Let \(K\) be a knot in \(S^3\).
Then \( h_k(K) > 0  \) whenever \( |k| < s(K) \).
\end{prop}

\begin{proof}
By the conjugation symmetry, it suffices to show that \(h_{k} (K) > 0
\) whenever \(s(K) > k \geq 0 \). 
Let \( \alpha \in \cfs (K) \)
 be a chain representing the generator of \(\hfs (K)
\). We filter \( \cfis (K) \) in the usual way, so that the filtered
quotients are all isomorphic to \( \cfs (K) \). The resulting  spectral
sequence converges at the \(E_2 \) term to \( \Z[u, u^{-1}] \). The
class \(u^n \) is represented by the chain \( \alpha _n = 
[\alpha, -n] + u^n \eta \),
where \( \eta \in \hfms (K)\).  Let \( p_k \: \cfis (K) \to \cfp (K, \spi _k) \) be the
projection. Then \( [p_s ( \alpha _n)] = 0 \) when \(n > 0 \) and is
nonzero when  \( n \leq 0 \). By
Corollary~\ref{Cor:hFormula},   \(h_k (K) > 0 \) if and only if \([ p_k
(\alpha _1)] \neq 0 \) in \( \hfp (K, \spi _k) \).

Observe  that \(p_k(\alpha _1)\) is contained in the subcomplex
\(C_{\spi_k} \) of   \( \cfp (K, \spi_k) \). We claim that if \(s(K)
> k \), \(p_k(\alpha_1)\) cannot be exact in \(C_{\spi_k} \). To see this,
consider the quotient complex \(D_k\) shown in
Figure~\ref{Fig:D_k}. If \(p_k(\alpha_1) \) is exact in \(C_k\), its image in
\(D_k\) must be exact there. But \(D_k \) is isomorphic to the quotient complex 
\( \cfs (K) / \{\bfz \ts | \ts A(\bfz) \leq k \} \), and the image of
\(p_k(\alpha_1) \) in \(D_k\) is the image of \( \alpha \) in this
quotient. By the
definition of \(s(K) \), the image of \(\alpha \) in \(D_k\) is nonzero
precisely when \( s(K) > k \). This proves the claim.

Since the boundary map \( \hfps (K) \to H_*(C_{\spi_k}) \) is always
trivial, \([p_k(\alpha_1)] \) is nonzero in \(\hfp (K,
\spi_k)\), which is what we wanted to prove.  
\end{proof}

Let \(W\) be the surgery cobordism from \(S^3 \) to \(K(0,1)
\). The proposition says that the induced map \( \cfm (S^3) \to \cfm
(K(0,1), \spi_k) \) is nontrivial whenever \(|k| < s(K) \). In fact,
\Ozsvath and \Szabo have shown in \cite{OS10} that the map is
nontrivial if and only if \(|k| < s(K) \). In particular, for any knot
\(K\) in \(S^3\), either all the \(h_k(K)\)'s or all the 
\(h_k(\overline{K})\)'s must vanish. 

\subsection{Cancellation}
\label{Subsec:Cancellation}
Proposition~\ref{Prop:s=Top} indicates the geometric significance of the invariant
\(s(K) \). Although we have no general method of computing  this
invariant, we describe in this section an algorithm for finding  the {\it homological} grading of
the generator of \( \hfs (K) \). When combined with information about the
stable complex, this is often sufficient to determine \(s(K) \). If
\(K\) is perfect, for example, \(s(K) \) is equal to the homological
grading of the generator. 

\begin{figure}
\includegraphics{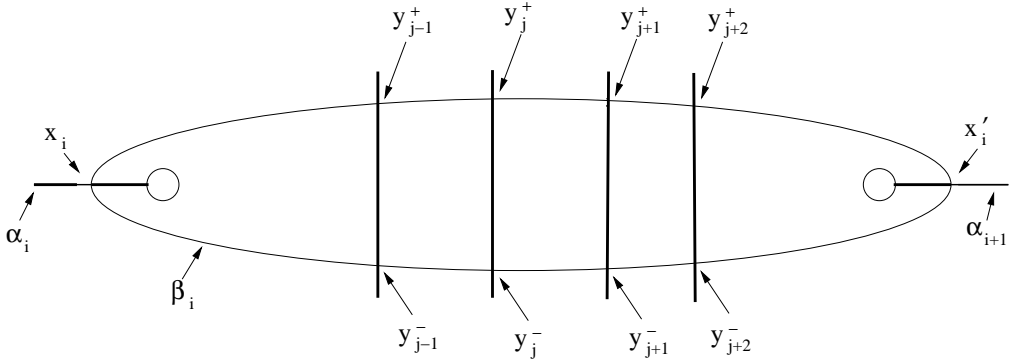}
\caption{\label{Fig:Beta} A \( \beta \) handle in a Heegaard splitting
  of \( \Ko\).}
\end{figure}

Recall from section~\ref{Subsec:HeegaardKnot} that a knot complement
admits a Heegaard splitting whose \( \beta \) handles are the
boundaries of regular neighborhoods of all but one of the
overbridges. To label the \( \alpha \) and \( \beta \) handles of 
this splitting, we start at some underbridge, which we label
\( \alpha _1 \), and travel along the knot, labeling the successive
overbridges and underbridges \( \beta _1, \alpha _2, \beta _2 \), {\it
  etc.} We omit  \( \beta _g \), and twist up around   \( \alpha _g
\).  The resulting Heegaard diagram may be divided into three regions: the
tubes, the {\it interior regions} bounded by the \( \beta \)'s in the
plane of the diagram, and the {\it exterior regions}, which make up
everything else in the plane of the diagram. We  orient
 \(m\) so that the basepoint \(z_s \) lies in the exterior region.

 As shown in Figure~\ref{Fig:Beta}, most of the
intersections of \( \beta _i \) with the \( \alpha \)'s fall
naturally into pairs \(y_j^\pm \).
 The exceptions are the points at either end of
the overbridge. We label these points \(x_i\) (the intersection
between \( \beta _i \) and \( \alpha_i \)) and \( x_i'\) (the
intersection between \( \beta _i \) and \( \alpha _{i+1} \)). 

\begin{prop}
\label{Prop:Cancellation}
The homological grading of the generator of \(\hfs(K) \) is the same as
the grading of  \( {\bf x} = \{ x_1,x_2,\ldots, x_g \} \).
\end{prop}

Before giving the proof, we pause to discuss the philosophy
behind it. Like all Floer theories, the \OS Floer homology is
 a Morse theory at heart. The Morse functional on \( \cfhat \)
 is  the area
function on \( \Sigma \); {\it i.e.} if \( \phi \in \pi_2(\bfy, \bfz)
\) with \(n_z(\phi) = 0\), \( \AAA(\bfy) - \AAA(\bfz) \) 
is the signed area of \( \DDD (\phi ) \).
% In
%particular, if \(A(\bfz) > A(\bfy) \), there can never be a
%differential from \( \bfy \) to \( \bfz \). (Proof: if \( \phi \in
%\pi_2(\bfy, \bfz) \) some region of \(\DDD (\phi) \) must have
%negative coefficients.)
 Our
ability to deform the metric on \( \Sigma \) gives us great freedom to
deform this Morse functional. 

The set of generators for \( \cfs(\E) \)
may be  complicated, but  their homology is always very simple:
  \( \hfhat (S^3) \cong \Z\). Now  there is a
well-known trick to turn a simple Morse function \(f\) into a more complicated
one. By making small perturbations, we can introduce 
cancelling pairs of critical points \(v_i^\pm\), such that  \(f(v_i^+) -
f(v_i^-)  \) is very small and there is a unique differential from \(v_i
^+ \) to \(v_i ^- \). For an appropriate choice of
metric on \( \Sigma \), \( \cfs (\E) \) exhibits precisely this sort
of behavior.

\begin{figure}
\includegraphics{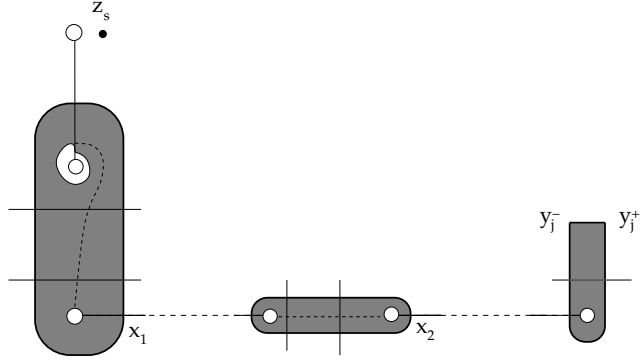}
\caption{\label{Fig:ThinDomain} Domain of the differential from \(
  \bfy_i^+ \) to \( \bfy_i^-\). The number of tubes in the domain can
  be anywhere from \(0\) to \(g-2\).} 
\end{figure}

To be specific, the cancelling pairs will be of the form
\begin{equation*}
\bfy_i^\pm = 
\{ y^\pm_j, x_1, x_2, \ldots, x_n,z_1,\ldots,z_{g-n-1} \},
\end{equation*}
 where 
 \( y^{\pm}_j\) is a pair of intersections like that of Figure~\ref{Fig:Beta}
 and the \(z_i\)'s are other points. 
 The form of the domain of \( \phi \in \pi_2(\bfy_i^+,
\bfy_i^-) \) is sketched in Figure~\ref{Fig:ThinDomain}. It is
supported on the interior region and tubes of the Heegaard
diagram. If we want \(\bfy_i^\pm \) to behave like a 
pair of newly introduced critical points, the difference \(\AAA(\bfy_i^+) -
\AAA( \bfy_i^-) \) should be  small. We can accomplish this for all
of the pairs in question by taking the areas of the interior region
and the tubes to be very small --- much smaller than the area of any
region in the exterior of the Heegaard splitting. 

Next, we need to know that  there is
 a  differential from \( \bfy_i^+ \) to \( \bfy_i^- \). For the
 moment, we simply state this as a fact, leaving the proof to the appendix.

\begin{lem}
\label{Lem:AnnularDiff}
There is  a unique differential from \(\bfy_i^+ \) to \(\bfy_i^-\)
 in \( \cfs (E) \). 
\end{lem}

\noindent We refer to the differentials of this form as 
``short differentials.''

We now define an equivalence relation on the set of generators  \( \TAB \) by requiring that
points connected by short differentials be equivalent. 
Every equivalence class contains \(2^m \) elements
 of  the form \( \{y_{j_1}^\pm,
y_{j_2}^\pm, \ldots, y_{j_m}^\pm, z_{m+1}, z_{m+2}, \ldots, z_{g-1} \}
\),  where
the \(z_i\)'s  are ``unpaired points'' ---  either \(x_j\) or \(x_j'\). 
For each equivalence class \( \eqy_i \), we let
  \(I(\eqy_i) \) be the smallest interval containing \(\AAA(\bfy) \) for
  all \( \bfy \in \eqy_i \). We say that \( \eqy_i\)  is {\it isolated }
if \( I(\eqy_i) \) is disjoint from all of the other \( I (\eqy_j ) \). 
Since the difference  in \(\AAA\) between two elements of the same
equivalence class is much smaller than the difference between two
elements in different classes, we expect that the \( \eqy _i \) should
be isolated. This is indeed the case:

\begin{lem}
\label{Lem:Isolated}
For a   
generic choice of large areas for the exterior regions of \(\Sigma \),
all equivalence classes are isolated.
\end{lem}

\begin{proof}
First, we observe that the length of \( I(\eqy) \) is bounded by
\((g-1)\) times the area of the interior and tubular regions.
 Thus by choosing these areas to be sufficiently small,
 we can assume that all of
the \( I(\eqy_i ) \) have length less than \(1\).

Suppose \( \bfy _1, \bfy _2 \in \TAB \), and let \( \phi \in \pi_2
(\bfy_1, \bfy_2 ) \). We claim that if \(
\DDD(\phi) \) is supported in the interior and tubular regions of \(
\Sigma \), then \(\bfy_1 \) and \(\bfy_2 \) are in the same equivalence
class. Indeed, the components of the \( \alpha \) curves in \( 
\partial \DDD(\phi) \) must be supported in the interior regions. For
homological reasons, they cannot go over the tubes. Thus they
are either trivial or join some \( y^+_j \) to the corresponding 
\(y^-_j\). This  implies the claim.

Next, we choose representative elements \(\bfy_i \in \eqy _i \) and
topological disks \( \phi _{ij} \in \pi_2 (\bfy _i, \bfy _j) \), and  
consider the differences \(\AAA_{ij} = \AAA(\bfy _i) - \AAA(\bfy _j) =
\AAA(\phi_{ij}) \).
 The \(
\alpha \) and \( \beta \) curves divide the exterior region into
various components, each of whose areas we can vary separately. If
\(a_k \) is the area of the \(k\)th component, we have 
\begin{equation*}
\AAA_{ij} = \sum_k n^k_{ij} a_k + C_{ij}
\end{equation*}
where \(n^k_{ij} \) is the multiplicity of \(\phi _{ij} \) over the
\(k\)th component and \(C_{ij} \) is the area of
\( \DDD (\phi_{ij}) \) over the interior and tubular regions.
If \(i \neq j \), \(y_i \) and \(y_j \) are in different equivalence classes,
so  \(n^k_{ij} \neq 0 \) for some \(k\). 

To show that all of the \(\eqy _j \) are isolated, it 
suffices to show that \(|A_{ij}| > 2 \) for all \(i\neq j \). The set
 of \(a_k\)'s for which \(|A_{ij}| < 2 \) is a width two
neighborhood of some hyperplane in the space of all possible \(a_k\)'s.
 A finite union of such
neighborhoods cannot fill up the entire positive quadrant, so the
lemma is proved.
\end{proof}

Since differentials decrease \(\AAA\), we have a filtration on \( \cfs
(E) \) whose filtered subquotients \(C^i \) are generated by the elements
of  \( \eqy _i \). 

\begin{lem}
\label{Lem:Acyclic}
If \(\eqy _i \) contains more than one element, the 
 complex \(C^i \) is acyclic.
\end{lem}

\begin{proof}
So far the only assumption we have made about the interior and tubular
 regions is
that their total area is  small. We now place some more specific
requirements on the distribution of area within these
regions. Without loss of generality, we assume that their total area
should be bounded by 1. 

First, we require that the total area of the interior regions 
should be very small --- less than \(4^{-g}\). Most of the area will be
 concentrated on the tubes. We assign the tube between \( \beta _{g-1}
 \) and \( \beta_{g-2}\) an area of \(4^{1-g} \), the
 tube between \( \beta _{g-1} \) and \( \beta_{g-2}\) 
 an area of \( 4^{2-g}\), and so on until we reach the  tube between \(
 \beta _2 \) and \( \beta _1 \), 
 which has an area of \(4^{-2} \). 

Fixing a particular equivalence class \( \eqy _i \), we say that the
handle \( \beta _k \) is {\it active} if the points in \( \eqy _i \)
contain some \( y_j^\pm \) on \( \beta _k\), and {\it inactive} if
they contain \( x_k \) or \( x_k' \). If \( \beta _{k_1}, \beta
_{k_2}, \ldots, \beta _{k_m} \) is a list of the active \( \beta\)'s
in order of increasing \( k_i\), we can specify
an element \(  \bfy \in \eqy_i \)  by a sequence of \(m\)
\(+\)'s and \(-\)'s, where the first \(\pm \) determines whether  
\( \beta _{k_1}\)  contains \( y^+ \) or \(y^- \), 
  the second determines whether \( \beta _{k_2} \)
contains \( y^+ \) or \(y^- \) on the next smallest, 
 and so forth. With these conventions, it is easy to check
that the ordering which \(\AAA\) induces on the elements of \( \bfy \) is
the same as the lexicographic ordering on the sequence of \(+\)'s and
\(-\)'s. 

If \(\eqy _i \) has more than one element,
 we can use this area functional to filter \(C^i\) so that each graded
subquotient is generated by two elements of the form \(w+ \) and
\(w- \), where \(w\) is some sequence of \(+\)'s and \(-\)'s. By
Lemma~\ref{Lem:AnnularDiff}, there is a unique differential from \(w+
\) to \(w- \), so the filtered subquotients are all acyclic. Applying the
reduction lemma, we see that \(C^i \) is chain homotopy equivalent to
a complex with no generators, so it is acyclic too. 
\end{proof}

\begin{proof}(of Proposition~\ref{Prop:Cancellation})
Together with the reduction lemma, the previous result shows that 
 \( \cfs (E) \) is chain homotopy equivalent to a
complex generated by those \(\bfy \in \TAB \) which are the unique
elements in their equivalence classes. We claim that \( \{ x_1,x_2,
  \ldots x_g \} \) is the only such \( \bfy \). Indeed, since we
  twisted up around \( \alpha _g \), \(\bfy \) cannot contain \(x_g'
  \). Thus it must contain \(x_g \). But this implies it cannot contain
  \(x_{g-1}' \), so it must contain \(x_{g-1} \), and so forth. 
\end{proof}

\subsection{A Model for Negative Surgeries}
We conclude our discussion of the exact triangle by returning to the
question raised at the end of section~\ref{Subsec:Model}. As is often
the case, it turns out to be more convenient to consider the exact
triangle for {\it negative} \(n\)-surgery (\(n \gg 0\)): 
\begin{equation*}
\begin{CD}
@>>>\hfp (N(0,1), \spi_k) @>>> \hfp (\Nbar, \spi) @>f>>
\hfp (N(-n,1), \spi _k) @>>> .
\end{CD}
\end{equation*}
\begin{figure}
\includegraphics{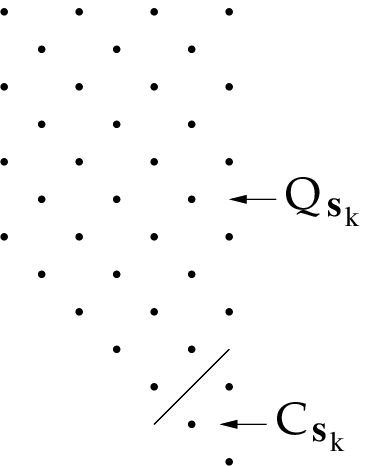}
\caption{\label{Fig:MinusSurg} The short exact sequence \(0 \to
  C_{\spi_k} \to  \cfpr (N, \spi) \to Q_{\psi_k} \to 0 \).} 
\end{figure}
The short exact sequence which models this triangle is shown in
Figure~\ref{Fig:MinusSurg}. The dual of the quotient complex
\(Q_{\spi_k}\) is
the subcomplex of \( \cfir (-N, \spi) \) which computes \( \hfm (-N,
\spi_k) \), so \(H_*(Q_{\spi_k}) \cong \hfp (N(-n,1), \spi_k) \) for \(n \gg 0
\). Thus we have a long exact sequence
\begin{equation*}
\begin{CD}
@>>>H_*(C_{\spi_k}) @>>>\hfp (\Nbar, \spi) @>{\pi_{s*}}>>
 \hfp (N(-n,1), \spi_k) @>>> .
\end{CD}
\end{equation*}
We will prove 
\begin{thrm}
\label{Thm:Model}
Suppose that \(\spi \) is a \( \spinc \) structure on \( \Nbar \)
which is fixed by the conjugation symmetry, and that this symmetry
acts trivially on \( \hfhat (\Nbar, \spi) \). Then
for \(k > 0 \), 
\( \ker f \cong \ker \pi_{s*} \)  and
\( \coker f \cong \coker \pi_{s*} \) as  \( \Q[u]\) modules.
\end{thrm}
 \noindent In particular, the theorem applies to any knot
 in \( \#^n(S^1 \times S^2) \). 

The map \( f \) is induced by the surgery cobordism \(W\) from 
\( \Nbar \)  to \(N(-n,1) \). To be precise, let \(x\) be the generator
 of \(H^2(W) \) such that the pullback of \(x\) to \(N(-n,1)\) is
 \(\spi_{i+1}-\spi_i\),  and denote by \( \vspi _i \) the \( \spinc \)
 structure on \(W\) with \( c_1(\vspi _i) = ix\). (Note that \( \vspi _i
 \) is only defined for \( i \equiv n \pod 2\).) The restriction of \(
 \vspi _i \) to \(N(n,1) \) is \(\spi_{(i-n)/2}\).  (The easiest way to see
 this to take \(n\) odd and consider \(\vspi _{\pm n}\), which induce
 the same \( \spinc \) structure  on \(N(n,1)\). Since the \( \vspi
 _{\pm n} \)  are conjugate symmetric, this \(\spinc \) structure must be \(
 \spi _0 \).)
Then if  \( f_i\: \hfp
(\Nbar ) \to  \hfp (N (-n,1), \spi _{(i-n)/2})  \)
 is the map defined by the pair \((W, \vspi _i)\),
\begin{equation*}
f = \sum_{i\equiv n + 2k \ (2n)} f_i .
\end{equation*}

It would be difficult to study all of the terms in this sum. When \(n\) is
large, however, we really only need to worry about 
 \(f_{\pm n +2k} \), which we call \(f_\pm\) for short.

\begin{lem}
\label{Lem:PlusMinus}
\( \ker f = \ker \ts (f_- + f_+) \) and {\em \( \coker f \cong \coker  
(f_- + f_+)\)}.
\end{lem}

\begin{proof}
If \(x \in \hfp (\Nbar) \) is homogenous in grading, the grading formula
tells us that 
\begin{equation}
\label{Eq:GrDrop}
\gr \ts f_-(x) - \gr \ts f_{an+2k}(x) 
\geq \frac{ - (-n+2k)^2 + (an+2k)^2}{4n} = \frac{a^2-1}{4}n +(a+1)k
\geq 0
\end{equation}
since \(n \gg k > 0 \). Thus \(f_-\) has the smallest drop in absolute
grading among all the \(f_i\) appearing in the sum.
Let us identify the relative gradings on 
\(\hfp (\Nbar, \spi) \) and \(\hfp (N(-n,1), \spi _k ) \) in such a way
that \(x\) and \(f_-(x) \) have the same grading. 

With this choice, we claim that there is some \(m\) independent of 
\(n\) such that
\(f\) is an isomorphism in degrees greater than \(m\).
 Indeed, it suffices to choose \(m\) so large that 
 all elements in \(\hfp(N(-n,1), \spi_k) \) and \(\hfp(\Nbar) \) with grading
 greater than \(m\) are in the image of \( \hfi \). Then the map 

\begin{equation*}
  f_{\geq m} \: \hfp
(\Nbar)_{\geq m} \to \hfp(N(-n,1), \spi_k) _{\geq m} 
\end{equation*} is the same as 
\begin{equation*}
  f_{\geq m}^{\infty} \: \hfi
(\Nbar)_{\geq m} \to \hfi(N(-n,1), \spi_k) _{\geq m}.
\end{equation*}
  Since  \( f_{\geq m}^\infty = f_-^{\infty} +
\epsilon \), where \(f_-^{\infty}\) is an isomorphism and  \(
\epsilon \) is strictly lower in degree, the claim holds.
%If \(f_{\leq M} : \hfp (N(n,1), \spi _k )_{\leq M} \to 
%\hfp (\Nbar, \spi)_{\leq M}\) be the restriction of \(f\) to elements of
%degree \( \leq M\). 

We now observe that \( \ker f = \ker f_ {\leq m} \)
and \( \coker f \cong \coker f_{\leq m} \). The first equality is
clear. To prove the second, we consider the composition
\begin{align*}
\coker f_{\leq m} = \frac{\hfp (N(-n,1), \spi_k)_{\leq m}}{ \im
 f_{\leq m}} & \into
 \frac{\hfp (N(-n,1), \spi_k)}{  \im f_{\leq m}} \\  & \to \frac{\hfp (
N(-n,1), \spi_k)}{  \im
 f}  = \coker f.
\end{align*}
It is easy to check this map is an isomorphism. 

 On the other hand, when
\(n\) is very large, equation~\ref{Eq:GrDrop} implies that
 \( f_{\leq m} = (f_- + f_+)_{\leq m} \).
Since \( f_- + f_+ \) is
also an isomorphism on elements of degree \( > m \), we have
\begin{gather*}
\ker f = \ker f_{\leq m} = \ker (f_+ + f_-)_{\leq m} = \ker (f_+
+ f_-) \\
\coker f \cong \coker f_{\leq m} = \coker (f_+ + f_-)_{\leq m} \cong  
\coker (f_+ + f_-).
\end{gather*}
\end{proof}

\noindent Let \( \pi_s \: \cfpr (\Nbar )  \to  \cfpr (N, \spi _k) \) be the
 projection.

\begin{prop}
\label{Prop:fMinus}
Up to filtration preserving isomorphism of the target, \(f_- = \pi_{s*} \).
\end{prop}

\begin{figure}
\includegraphics{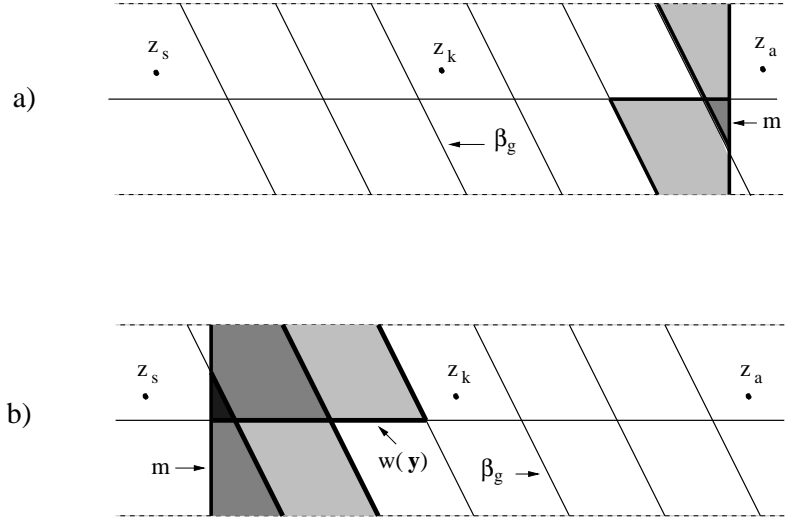}
\caption{\label{Fig:TriangleSpiral} The Heegaard diagrams used to
  compute {\it a)} \(f_-\) and {\it b)} \(f_+\). In each figure, the 
 domain of a typical triangles \( \Delta_\bfy \) is shaded. }
\end{figure}

To prove this,  we study the map \( \fmi \: \cfi (\Nbar)  \to
\cfi (N(-n,1))\), which is defined by counting holomorphic triangles in
the relevant Heegaard triple diagram.
 In  Figure~\ref{Fig:TriangleSpiral}a, we show a portion of the spiral
 region for this diagram. (It lies on the cylinder obtained by
 identifying the top
 dotted line with the bottom one.)  We position the attaching handle
 \(\emm \) so that \(w(\bfy) \) is on the stable side of \( \emm \)
 for any generator \(\bfy \) of \(\cfs (N(-n,1))\). 
 Then each such \( \bfy \) has  a unique small triangle \( \Delta _\bfy \)
 connecting it to the
corresponding generator \( \bfy ' \)   of \( \cfhat (\Nbar) \). The
part of the domain of this triangle which lies inside the spiral
region is shown in the figure. It is easy to see  that 
\( \mu (\Delta_\bfy) = 0 \) and \( \# \MMM (\Delta_\bfy) = \pm 1 \). 

\begin{lem}
The \( \spinc \) structure
\( \spi_{z_k} (\Delta_\bfy)\)  determined by \( \Delta _{\bfy} \) is
\(\vspi _- \).
\end{lem}

\begin{proof}
Clearly \(\spi_{z_k} (\Delta_\bfy) \) restricts to \( \spi _k \) on
\(N(-n,1) \). Let us take both \(n\) and \(k\) to be very large (but
with \(k < n\)). 
Since \( \# \MMM (\Delta_\bfy) = \pm 1 \) and \(n_{z_k}(\Delta_\bfy) =
0\),  \( \spi_{z_k} (\Delta_\bfy) \) 
  gives a nonzero map on \( \cfhat \). When \(n\) is very large, the
proof of Lemma~\ref{Lem:PlusMinus} shows that this can happen only if 
\(\spi _{z_k} (\Delta _\bfy) \) is equal to \( \vspi_+ \) or \( \vspi _- \).
 Since \(k\) is very large as well, the same argument shows that the
degree of \(f_+ \) is much less than that of \(f_-\), so it must be
zero on \( \cfhat \) as well. Thus \(\spi _{z_k}(\Delta _\bfy) = \vspi_- \)
when \( k\) is large. The result for arbitrary values of
 \(k\) follows from the
analogue of Lemma 2.12 of \cite{OS1} for Whitney triangles. 
\end{proof}

\begin{lem}
The map \(\fmi \) respects the bifiltration on \( \cfi (\Nbar) \cong
 \cfi (N(-n,1) )\). 
\end{lem}

\begin{proof}
For the stable filtration, this is just the fact that \(n_{z_k}(\psi)
\geq 0 \) for any holomorphic triangle \( \psi \). To prove it for the
antistable filtration, suppose \(\bfy \) is a generator of \(\cfhat
(N(-n,1)) \) and that \(\bfz ' \) is a generator of \( \cfhat (\Nbar)
\) with the same \( \Z/2 \) grading as \( \bfy \). (If the \( \Z/2
\) grading is different, \( \bfy \) can never map to \( \bfz'\).)
Let \( \phi \in \pi_2(\bfy', \bfz') \) be the unique disk with \( \mu
(\phi ) = 0\). Then \( \psi = \Delta _\bfy \# \phi \in \pi_2(\bfy,
\theta, \bfz ')\) has \( \mu
(\psi) = 0 \) and \( \spi_{z_k}(\psi) = \vspi _-\). Thus
 \( \psi \) is the triangle which determines the \(\bfz ' \) component of \(
\fmi(\bfy) \). Since the \( \beta_g \) component of \( \psi \) is
supported inside the spiral, the usual argument shows that \(\fmi\)
preserves the Alexander grading. 
\end{proof}

\begin{proof} (of Proposition~\ref{Prop:fMinus})
Since \( f_- \) respects the bifiltration, we can reduce it to get a
map \( \cfpr  (\Nbar) \to \cfpr (N(-n,1), \spi_k) \), 
which we continue
to denote by \( f_-\). To prove the proposition, it suffices to show that 
 \(\fmi \) induces an isomorphism on
 corresponding copies of \( \scx^2_j \). (The ``dots'' in the
 diagram.)   This follows from a  standard argument
using the area filtration. We make the domains of the triangles 
\( \Delta _y \) very small, so that their area is less than that of
any other positive triangle. Then the map \( \scx^1_j(
\Nbar) \to \scx ^1 _j (N(-n,1)) \) has the form \( 1 + \epsilon \),
where \( \epsilon \) decreases the area filtration. This proves the claim. 
\end{proof}

\noindent {\bf Remark:} This argument should be compared with the
proof of Theorem 9.1 in \cite{OS2}. In fact, the domains of the triangles
used to define \( f_- \) are closely related to the domains of the
disks which define the differential \( \delta _1 \) described in that paper.
 We conjecture that \(f_-\) and \( \delta _1 \) are actually the same
 map. 

It seems difficult to compute the map \(f_+ \) from the Heegaard data of
Figure~\ref{Fig:TriangleSpiral}a. The domains of the triangles
defining this differential are related the domains of the disks which
define the differential \( \delta _2 \) of \cite{OS2}. In
particular, they all have very large multiplicity.

 Fortunately, there
is another way to understand \(f_+\). We use the slightly different
Heegaard diagram of Figure~\ref{Fig:TriangleSpiral}b. The
only thing that has  changed in this figure is the position of the  attaching
handle  \(\emm \), which has been shifted so that all the
\(w(\bfy)\)'s are on its
antistable side. Thus the chain complex \(\cfp (N(-n,1),
\spi _k) \) remains the same. When we go to compute the Floer homology
of \( \Nbar \), however, the basepoint \( z_k\) is on the opposite
side of \( \emm \), so we get the complex \( \cfa (N) \) rather than
\( \cfs (N) \).
 Let \( \pi_a \:  CF^+_{ar}(\Nbar) \to \cfpr (N(-n,1), \spi_k)
 \) be the projection. Then the same arguments used to prove
Proposition~\ref{Prop:fMinus} imply

\begin{prop}
Up to filtered isomorphism of the target, \(f_+ = \pi_{a*} \). 
\end{prop}

%\begin{figure}
%\caption{\label{Fig:TriangleSpiral2} Heegaard diagram used to compute
%  \(f_+\). Part of the domain of a small triangle is shown. }
%\end{figure}

 Thus both \( f_-\) and \( f_+ \) are induced by quotient maps at
the chain level, but the maps are from {\it a priori} different
complexes, both of which compute \( \hfp (\Nbar)\). 
 By Proposition~\ref{Prop:Symmetry}, however,
 there is an isomorphism \( \iota \: \cfpr(N) \to \cfpa (N) \), for which
the  induced map \( \iota_* \: \hfp(\Nbar, \spi) \to \hfp (\Nbar, \spi)
\)  is the
conjugation symmetry. We will use this fact to prove the theorem.

\begin{proof} (Of Theorem~\ref{Thm:Model})
By hypothesis,  \( \iota _*\) acts trivially on 
\(\hfp (\Nbar, \spi) \), so 
\begin{equation*}
  f_- + f_+ =  \pi_{s*} +  (\pi_{a} \iota)_*.
\end{equation*}
 If \(A = \ker \pi_{s}\) and \(B = \ker \pi_a \iota\), then the fact
 that \(k > 0 \) implies that \(A \subset B \).

We show that \( \ker \pi_{s*} = \ker (\pi_{s*}+(\pi_{a} \iota)_*)
\). Indeed, if \(\pi_{s*}(x) = 0\), \(x\) can be represented by a
chain in \(A\). Since \(A \subset B\), \((\pi_{a} \iota )_*(x) = 0 \) as
well. Conversely, suppose \( \pi_{s*}(x)+(\pi_{a} \iota)_*(x) = 0\). Write
\(x = x_i + x_{i-1} + \ldots \), where \(x_i\) is homogenous of degree
\(i\). Since \(k> 0 \), \(( \pi_{a} \iota)_* (x_j) \)
 has degree strictly lower
than \( \pi_{s*}(x_j) \). Thus if we consider the equation
 \( \pi_{s*}(x)+(\pi_{a}\iota)_*(x) = 0 \) one degree at a time, we obtain a
 series of equations 
\begin{align*}
\pi_{s*}(x_1) &= 0  & \pi_{s*}(x_2) & = 0  & \ldots  & &
\pi_{s*}(x_{i-k}) &= 0 \\
\pi_{s*}(x_{i-k-1})& = - (\pi_{a} \iota)_*
(x_1) & \pi_{s*}(x_{i-k-2})& = -
(\pi_{a}\iota )_*(x_2) & \ldots &
\end{align*}
Since \( \pi_{s*}(x_1) = 0\), \(x_1\) can be represented by a chain in
\(A\), so \((\pi_{a} \iota) _*(x_1) = 0\) as well. Thus \(
\pi_{s*}(x_{i-k-1})=0\). Repeating, we see that \(\pi_{s*}(x_j) = 0 \)
for every \(j\). This proves the claim. 

To show that \(\coker \pi_{s*} \cong \coker (\pi_{s*} + (\pi_{a}
 \iota)_* ) \),
consider the map \( g \: \cfpr(N, \spi_k) \to \cfpr (N, \spi_k) \) 
given by \(g = 1 + i \pi
\), where \( \pi \) is the projection which kills the remaining
generators of \(B\), and \( i \) is the reflection
isomorphism. Since \(k > 0\), \(i \) strictly reduces the Alexander
grading. This implies that 
 \(g \) is an isomorphism, so \(g_*\) is  an isomorphism
on homology. Clearly \(g(\pi_s(x)) = \pi_s(x) + \pi_a \iota (x) \), so \(g_*\)
defines an isomorphism from \(\coker \pi_{s*} \) to \( \coker 
(\pi_{s*} + (\pi_{a} \iota )_*) \). 

\end{proof}

%\begin{figure}
%\includegraphics{figs/AB.eps}
%\caption{\label{Fig:AB} The subcomplexes \(A = \ker \pi_s \) and \( B
%  \ker \pi_a \iota \) }
%\end{figure}

%% file: alternating.tex
\section{Alternating Knots}
\label{Sec:Alternating}

Among knots with small crossing number, perfect knots are
very common. Our main result in this direction is
Theorem~\ref{Thm:Alternating=QP}, 
which states that alternating knots which satisfy a certain
``smallness'' condition are perfect. Small alternating knots are
significantly more common than two-bridge knots; we have checked by
computer that there is only one non-small alternating knot with fewer than 11
crossings (\(10_{123} \) in Rolfsen's tables) and only 6 non-small
alternating knots with 11 crossings. In fact, the techincal condition
of smallness is unnecessary: using a different Heegaard splitting from
the one considered here, \Ozsvath and \Szabo have shown that all
alternating knots are perfect. 

Combining these facts  with a few direct
calculations, we have

\begin{thrm}
All  alternating knots with 10 or fewer crossings are perfect
with \( s(K) = \sigma (K) /2 \). The same is true for nonalternating
knots with 9 or fewer crossings, except for \(8_{19}\) and 
 \(9_{42} \). 
\end{thrm}

The proof of  Theorem~\ref{Thm:Alternating=QP} is the main topic of
this section. Afterwards,
we briefly discuss computations for nonalternating knots and describe some
interesting correspondences between \( \cfr \) and Khovanov's Jones
polynomial homology. 

\subsection{Heegaard splittings and Generators} 
\begin{figure}
\includegraphics{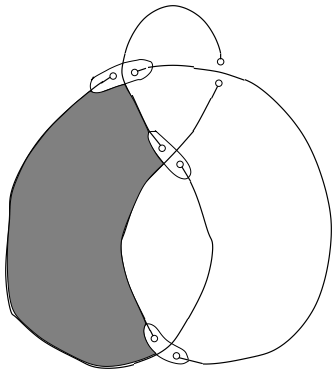}
\caption{\label{Fig:Alternating} A Heegaard splitting of the
  complement of the figure 8 knot. The shaded region is the domain of
  a disk differential}
\end{figure}

Let \(K\) be an
alternating knot. To compute \( \cfr (K) \), we use a Heegaard
splitting  \( (\Sigma, \alpha, \beta ) \)
of \(\Ko \) derived from a pair \((\DDD, c_1) \), where \( \DDD \) is an
 alternating diagram of \(K\) and \(c_1 \) is a crossing of \( \DDD
 \). The basic form of this splitting is shown in 
Figure~\ref{Fig:Alternating}. The genus of \( \Sigma \) is the
crossing number of \(K\). 
Each edge of \( \DDD \) corresponds to an \(\alpha \) handle, and each
crossing \(c_i\) with \(i >1 \) corresponds to a \( \beta \)
handle \( \beta_i \). {\it Warning:} such diagrams can be slightly misleading, since the ``knot diagram'' obtained by omitting the \( \beta\)'s is actually a diagram of \( -K \)! (Remember, the \( \beta\)'s are tubular neighborhoods of the overbridges.) Also, one must resist the temptation to attach the  ``tubes''  to the two circles enclosed by a \( \beta\) handle, rather than to the two ends of an edge.

 Since the diagram is alternating,  the \( \alpha _i\)'s and \( \beta _j \)'s  all have the same local models. Each \( \beta _j \)
 has two {\it end intersections} \( x_j^\pm \) and two {\it middle
   intersections} \(y_j ^ \pm \). The same is true for all of the \(
 \alpha _i\)'s, except that there is no \( \beta _1 \), so one
  \( \alpha _i \) (call it \(\alpha _1\)) is missing its middle
intersections, and two others are missing one of their end intersections. 

We twist up around \( \alpha_1 \) and study the generators of 
 the resulting complex \( \cfs (\E) \). These generators are in
 one-to-one correspondence with the elements of \( \TABO \). 
To describe them, we introduce the following

\begin{defn}
A marked partial resolution (or MPR for short)
 of \(K\) consists of the following data:
\begin{enumerate}
\item To each crossing \(c_i\) with \(i>1\), the assignment of either
  a resolution \( \langle \smoothing \rangle \) or \( \langle
  \hsmoothing \rangle \) or a sign \(+\) or \(-\).
\item A set of  components of the resulting partially resolved diagram
(the ``marked components'') such that 
\begin{enumerate}
\item No marked component contains any unresolved crossings.
\item Exactly one marked component passes through each resolved
  crossing.
\end{enumerate}

\end{enumerate}
\end{defn}

\begin{prop}
The generators of \( \cfs (\E) \) are in one-to-one correspondence
with the set of marked partial resolutions of \(K\).
\end{prop}

\begin{proof}
Suppose we are given a point \( \bfy \in \TABO \). The corresponding
marked resolution is constructed as follows. For each \(i>1\), \( \bfy
\) contains a unique point on \( \beta _i\), which is either a middle
intersection or an end intersection. If it is a middle intersection,
we label \(c_i
\) with a \(\pm\), according to the sign of \( y_i^\pm \in \bfy
\). If it is an end intersection --- call it \(x\) ---, we resolve \(c_i \). 
Let \( \alpha_j \)
be the overbridge passing through \(c_i\). Note that since \(i \neq
1\), \(j \neq 1 \) as well. Thus \( \bfy \) contains a point on \( \alpha
_j \). But \( \alpha _j\)'s  middle intersections
are not in \( \bfy \), so exactly one if its end intersections must
be. Call this  end intersection \(x'\).
 We resolve \(c_i \) so that the overbridge containing \(x\) 
 connects with the overbridge containing \(x'\), as illustrated in
Figure~\ref{Fig:Resolution}. Finally, we mark the component containing
\(x \) and \(x'\). 

\begin{figure}
\includegraphics{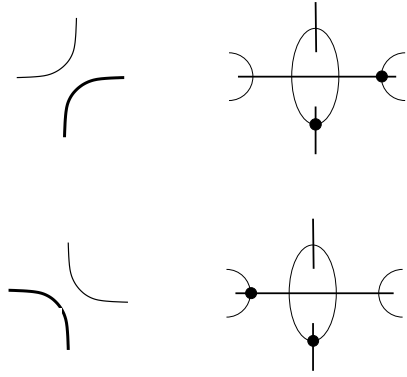}
\caption{\label{Fig:Resolution} Resolving a crossing with an end
  intersection. The marked component is shown by the heavy line.} 
\end{figure}

Let us check that the resulting diagram satisfies the definition of an
MPR. The basic point is that we can ``follow'' a
marked segment. If we start at a crossing \(c_i\) containing an end
intersection \(x\), the marked
segment takes us to a new end intersection \(x'\) in \(c_{i'}\). 
Thus \(c_{i'}\) must also be a resolved crossing. When we resolve it,
the new marked segment continues the previous one. If we
repeat this process,  we must eventually return to our original
end intersection \(x\). 
It is clear from this description that 
 a marked component does not contain any unresolved
crossings. As for the second condition, a resolved crossing \(c_i\) must
contain one marked component by definition. It cannot contain more
than one, since then both end intersections of \(
\beta _i \) would be in \( \bfy \). 

Now suppose we are given a marked partial resolution of \(K\). The
construction of the corresponding \( \bfy \in \TAB \) is
straightforward. Each \(c_i \) with \(i>1\) contributes an element of
\(\bfy\). If \(c_i \) is  marked \(\pm \), \( \bfy \) contains
\(y_i^{\pm} \). If \(c_i \) is resolved, \( \bfy \) contains whichever
end intersection is on the marked component. We need to check that 
each \( \alpha _j \) with \(j>1\) contains a unique point in \( \bfy
\).  If the  crossing \(c_{i}\)
through which \( \alpha _j \) passes is  marked \( \pm \), one of
\(\alpha_j\)'s middle intersections is in \(\bfy \), and neither end
intersection is. (Such an end intersection would have to lie on a
marked component passing through \(c_i\).)
On the other hand, if \(c_i \) is resolved, exactly one end
intersection of \( \alpha _j \) is on a marked component.
\end{proof}

The generators corresponding to diagrams in which no crossing is
resolved are of special interest to us:

\begin{defn}
The {\it base generators} of \(  \cfs (\E) \) are the \(2^{g-1}\)
elements of \( \TABO \) corresponding to MPRs in which every
crossing (except \(c_1 \)) is labeled with a \(+\) or a \(-\).
\end{defn}

\noindent Recall that in section~\ref{Sec:Filtration} we defined a grading
\( \tilde{\mu} = \gr - A \) on \( \cfs (\E) \). 

\begin{lem}
\label{Lem:BaseGenerate}
All of the base generators have the same value of \( \tilde{\mu} \).
\end{lem}

\begin{proof}
Let \(\bfy ^+ = \{y_1^\pm, y_2^\pm, \ldots, y_i^+, \ldots,  y_g^{\pm}
\} \) and \(\bfy ^- = \{y_1^\pm, y_2^\pm, \ldots, y_i^-, \ldots,  y_g^{\pm}
\} \) be a pair of base generators whose signs differ at a single
crossing. Then
it is not difficult to see that \(\bfy ^+ \) and \( \bfy ^- \) are
connected by a punctured annular differential of the sort described in
Lemma~\ref{Lem:AnnularDiff}. Any two generators connected by such a
differential have the same value of \(  \tilde{\mu} \). 
\end{proof}

To prove Theorem~\ref{Thm:Alternating=QP}, we will use the reduction
lemma to compare \( \cfs (\E) \) with a complex generated by a subset of
the basepoint generators. The first step in this process 
 is to divide the generators of \( \cfs (\E) \) into groups.
To every \( \bfy \in \TABO \), we  assign a base generator \(b(\bfy) \)
 by leaving all the labeled crossings of the MPR
 of \(\bfy \)  unchanged and replacing  each resolved crossing with a \(+\)
or a \(-\), depending on whether the marked component contains
\(y_i^+\) or \(y_i^-\). This assignment divides \(E\) into \(2^{g-1}\)
equivalence classes, each containing a single base generator. We
denote the class containing the base generator \( \bfy \) by \( \eqy
(\bfy) \). 

\subsection{Disk Differentials} We now describe  a special class of
differentials  between elements in the same equivalence class. An
example of such a differential for the figure-eight knot
is shown in  Figure~\ref{Fig:Alternating}. The following lemma
guarantees that the shaded domain actually has \( \# \MMM(\phi) = \pm 1
\). 

\begin{lem}
\label{Lem:DiskyDiffs}
Suppose \( \DDD (\phi) \) is a multiplicity one disk with all interior
corners. Then \( \mu (\phi) =1 \) and  \( \# \MMM(\phi) = \pm 1 \). 
\end{lem}

\noindent The proof is given in the appendix.

\begin{figure}
\includegraphics{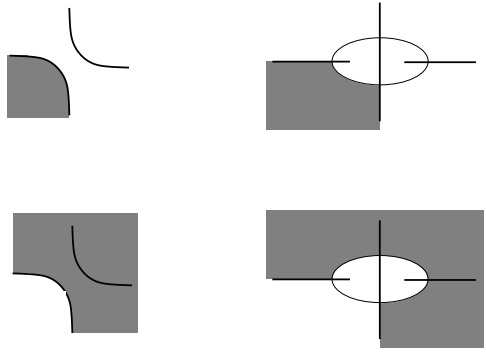}
\caption{\label{Fig:Corners} \( \DDD( \phi ) \) in the neighborhood of
  a crossing.}
\end{figure}

More generally, suppose \(C\) is a marked 
component of the MPR for \(\bfy \in E\). \(C\)
bounds two disks in the \(S^2\) of the diagram. We say that \(C\) is
{\it small} if one of these disks does not have any crossings in its
interior.  In this case, \(C\) bounds a disk  \( \DDD (\phi) \)
in \(\Sigma \). Near a crossing, \( \DDD (\phi) \) looks
like the local models shown in Figure~\ref{Fig:Corners}. 
 \( \DDD (\phi) \)  satisfies the conditions of
Lemma~\ref{Lem:DiskyDiffs}, so we get a differential between 
 \( \bfy \) and another generator \( \bfy'\), whose MPR is obtained
 by ``forgetting'' \(C\), {\it i.e} replacing each end intersection on
 \(C\) by the adjacent middle intersection. Clearly \( \bfy \) and
 \(\bfy'\) are in the same equivalence class. 

Motivated by this fact, we make the following 
\begin{defn}
 An alternating
 knot \(K\) is small if there is some diagram \(\DDD \) of \(K\) and choice
 of  crossing \(c_1\) such that every marked component of
 every MPR of \((\DDD, c_1) \) is small.
\end{defn}

From now one, we assume that \(K\) is a small knot, and that \( (\DDD, c_1) \) is the corresponding small diagram.
We call the differentials described above {\it disk
  differentials}. Since their domains are 
  supported away from the spiral region, disk differentials preserve the
  Alexander grading. Thus all the elements of \(\TABO \) in a given
  equivalence class have the same Alexander grading. Moreover, any
  element in an equivalence class can be connected to its base
  generator  by a series of disk differentials, so it is easy to work
  out the homological grading on \( \cfs (\E) \). (These statements
  still hold even if \(K\) is not small.)

All the disk differentials in a given equivalence class
are compatible with each other in the following sense:

\begin{lem}
\label{Lem:Compatible}
Suppose \( \bfy \) is a base generator, and let 
\begin{equation*}
\CCC(\bfy) = \{ C
\ts | \ts C \ \text{is a marked component of an MPR} \ \bfz \
\text{with} \ b( \bfz) = \bfy \}.
\end{equation*}
Then no two elements of \( \CCC (\bfy) \) pass through the same
crossing of \(K\). 
\end{lem}

\begin{proof}
Two elements of \( \CCC (\bfy) \) which share an edge must be the
same. Indeed, if we are given an edge of \( C \in \CCC (\bfy) \), we
can reconstruct \(C\) by following it along from one crossing to the
next. Each time we come to a new vertex, \(\bfy \) tells us which way
the resolution has to go. 
Suppose \(C_1, C_2 \in \CCC (\bfy) \) are two different
components coming into the same crossing. They must leave by the same
edge, which contradicts the fact that they are different.
\end{proof}

\begin{cor}
\label{Cor:PowerSet}
The elements of \( \eqy (\bfy) \) are in bijective correspondence with
subsets of \( \CCC (\bfy) \).
\end{cor}

\begin{proof}
 Given an MPR  \( \bfy ' \in \eqy ( \bfy )\),
 we map it to its set of marked components. 
In the other direction, given  \(\{ C_1,C_2, \ldots, C_m \} \subset \CCC (\bfy)
\), we start with the MPR
for \( \bfy \) and resolve all the crossings contained in the \(C_i\).
By Lemma~\ref{Lem:Compatible}, no crossing is contained in more than
one of the \(C_i\). Thus we can resolve each crossing to get the
\(C_i\)'s, and then mark them to get a valid MPR. It is easy to check that these two procedures define inverse maps.
\end{proof}

\subsection{Cancellation} The remainder of the proof is similar
to the proof of Proposition~\ref{Prop:Cancellation}. We use the area
functional to filter \( \cfs (\E) \) so that all the generators in an
equivalence class are in the same filtered subquotient. Next, we use a
cancellation argument to show that any equivalence class with more
than one element is acyclic. Finally, we  apply the reduction lemma to
prove the theorem. 

\begin{lem}
\label{Lem:Isolated2}
For an appropriate choice of metric on \( \Sigma \), the \(\eqy (\bfy)
\) are all isolated in the sense of Lemma~\ref{Lem:Isolated}.
\end{lem}

\begin{proof}
We choose a metric on \( \Sigma \)
 in which the area of all the exterior regions
put together is much less than the area of any one of the tubes. 
 Any \(\bfy \in E\)  is connected to \( b(\bfy ) \) by a sequence of
 disk differentials whose domains contain only exterior regions. Thus it
 suffices to show that there is a choice of areas on the tubes for
 which the base generators are widely separated. Since the base
 generators are connected by annular differentials, this can be
 accomplished by an argument like that of Lemma~\ref{Lem:Acyclic}.
\end{proof}

We take such a metric and use the resulting area functional to filter \( \cfs
(\E) \) so that each filtered subquotient is generated by the elements
of some \( \eqy (\bfy) \). 

\begin{lem}
\label{Lem:Acyclic2}
If \( \eqy (\bfy ) \) contains more than one element, the complex it
generates  is acyclic.
\end{lem}

\begin{proof}
By Corollary~\ref{Cor:PowerSet}, the elements of  \(\eqy (\bfy) \) 
correspond to  subsets of \( \CCC (\bfy ) \).
If \( \CCC (\bfy) \) has \(n\) elements,
  we can represent each  \(\bfy ' \in \eqy (\bfy ) \)
by a sequence \( {s_i} \) of \(n\)  \(1\)'s and \( 0\)'s, where \(s_i \)
is \(1\) if \(C_i \in \CCC (\bfy) \) is marked, and \(0\) if it is
not.
 With this
convention, it is easy to see that 
\begin{equation*}
A(\{s_i\}) = \sum _{i=1}^n  s_i a_i
\end{equation*}
where \(a_i\) is the signed area of of the disk domain bounded by
\(C_i\). (The sign of \(a_i\) is determined by whether the
differential goes to or from the MPR with the marked component.)

We now argue as in Lemma~\ref{Lem:Acyclic}. We are free to vary the
\(a_i\)'s, so long as we keep them small, so we take {\it e.g.} \(a_i
= 4^{1-i} a_1 \). With this choice,  the area
filtration splits the complex generated by \( \eqy (\bfy )\) into
\(2^{n-1}\) subquotients, each containing two generators connected by the disk
differential corresponding to \(C_n\). By the reduction lemma, \(C
(\bfy ) \) is acyclic. 
\end{proof}

\begin{proof} (Of Theorem~\ref{Thm:Alternating=QP}.)
The area filtration restricts to a filtration on the complex \( \scx
_1 ^j (\E) \). Since disk differentials preserve the Alexander
grading, all of the elements of a given \( \eqy (\bfy ) \) belong to
the same  \( \scx_1 ^j \). Then
Lemma~\ref{Lem:Isolated2} shows that  we can filter \(\scx _1 ^j ( \E)  \) 
so that each filtered subquotient is generated by some \( \eqy (\bfy)\). By
Lemma~\ref{Lem:Acyclic2} and the reduction lemma, \( \scx_1 ^j (\E)\) 
is chain homotopy equivalent to a complex generated by those
base generators  \(\bfy \) such that \( A(\bfy) = j \) and  \( \eqy(\bfy ) \)
contains no elements other than \( \bfy\).
 But by Lemma~\ref{Lem:BaseGenerate}, all
base generators have the same value of \( \tilde{\mu} \). Thus \( \cfr
(K) \) is perfect. 
\end{proof}

\subsection{Small knots}

Although the condition of smallness is seems  quite restrictive,
there are some infinite families of small knots.

\begin{lem}
Any two bridge knot is small. Any alternating three-strand pretzel
knot is small.
\end{lem}

\begin{figure}
\includegraphics{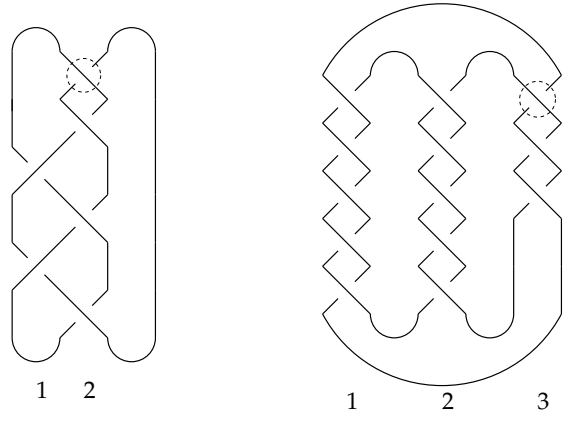}
\caption{\label{Fig:Twobridge} A two-bridge knot and an alternating
  three-bridge pretzel knot. The special crossings \(c_1\) are circled.}
\end{figure}

\begin{proof}
Any two bridge knot has an alternating diagram like that shown in 
Figure~\ref{Fig:Twobridge}a. If we choose \(c_1 \) as shown, any marked
component of an MPR of \(D\) 
 must be supported in the columns labeled 1 and 2, so it cannot
have any crossings in its interior. A very similar argument applies
to the three-strand pretzel knot shown in part {\it b} of the
figure. 
\end{proof}

\noindent Thus we recover the result  that two-bridge knots
are perfect. 

Unsurprisingly, alternating knots with small crossing number tend to
be small. 
Using a computer, we have checked that the only nonsmall alternating
knots with fewer than 12 crossings are \(10_{123} \), \(11^a_{236}\),
\(11^a_{327}\),\(11^a_{335}\),\(11^a_{356}\),\(11^a_{357}\), and 
\(11^a_{366}\) .
(Knots with 10 or fewer crossings are labeled according to Rolfsen's
tables in \cite{Rolfsen}; knots with more crossings are identified by their
{\it Knotscape} number.) As the crossing number increases, however, it
seems clear that the fraction of alternating knots which are small will
tend to zero.

%% file: appendix.tex
\section{Appendix on Differentials}

In this appendix, we collect  some explicit results about the
differentials in the \OS complex.
The fundamental question in this area is ``which domains correspond to
differentials?'' In other words, given the domain \( \DDD (\phi) \)
of a topological disk \( \phi \), we would like to be able to
determine \( \mu (\phi) \) and, if \( \mu (\phi) = 1\), 
 \( \# \MMM (\phi) \) from the
combinatorial information contained in the Heegaard diagram. Below, 
we give a  formula for \( \mu (\phi) \). This formula  has a purely
combinatorial interpretation for a wide class of domains, and it
should be possible to extend this interpretation to all domains. In contrast,
\( \# \MMM (\phi) \) 
 may not be determined  by the combinatorics of the
diagram. (See Lemma 5.7 of \cite{OS1} for an example.) Still, there
are many useful cases in which it is. We describe two such cases --- disky
differentials and punctured annular differentials --- 
which were  used  in sections~\ref{Sec:ExactTriangle} and  
\ref{Sec:Alternating}.

\subsection{Domains}

 We begin by reviewing the
definition of \( \DDD (\phi) \) (Lemma 3.6 of \cite{OS1}.) Given
 \( \phi \in \pi_2 (\bfy, \bfz) \), we choose some map \(u \: {
  D}^2 \to s^g\SSS  \) representing \(\phi\) and consider the
 incidence correspondence
\begin{equation*}
S = \{ (x,y) \in {
  D}^2 \times \Sigma \ts | \ts y \in u(x) \}.
\end{equation*}
 Without loss of
generality, we can choose \(u\) transverse to the diagonal in \(
s^g\SSS \), so that the projection \( \pi_1 \: S \to { D}^2 \)
 is a \(g\)-fold branched covering map.

 The projection \(\pi_2 \: S \to \SSS _g \) defines a
singular \(2\)--chain in \( \Sigma _g \), and thus (since \( \partial S\)
maps to \( \alpha \cup \beta \)), a cellular \(2\)--chain
\(\DDD(\phi) \) in the cellulation
of \( \Sigma _g \) defined by \( \alpha \) and \( \beta
\).   \( \DDD (\phi) \) depends only on \( \phi \).
Indeed, its multiplicity on a given cell is the intersection number
\(n_w(\phi)\)  between  \( \phi \) and \({w}\times s^{g-1} \SSS_g
\), where \(w\) is a point in the interior of the cell.

 We mark the two points on \(\partial D^2 \)
 which map to \(\bfy\) and
\( \bfz\). Their inverse images under \( \pi_1\) give rise to \(2g\)
marked points on \( \partial S\). Since \( \pi_1 \: \partial S \to
\partial D^2 \) is a covering map, \(\pi_2 \) sends these marked points 
 alternately to the \(y_i\)'s and 
and the \( z_j\)'s. The segments between them map alternately to \(
\alpha \) and \( \beta \). As a result, the cellular chain \(\partial
\DDD (\phi) \) is composed of \(2g\) arcs, one on each of the \(
\alpha \) and \( \beta \)--circles. Each \(y_i\) and \(z_j\) is the
endpoint of one such arc and the starting point of another. 
Conversely, it is not difficult to see that
any null-homologous 1-chain satisfying these conditions bounds \(
\DDD (\phi) \) for some \(\phi \in \pi_2 (\bfy, \bfz) \). 

\begin{figure}
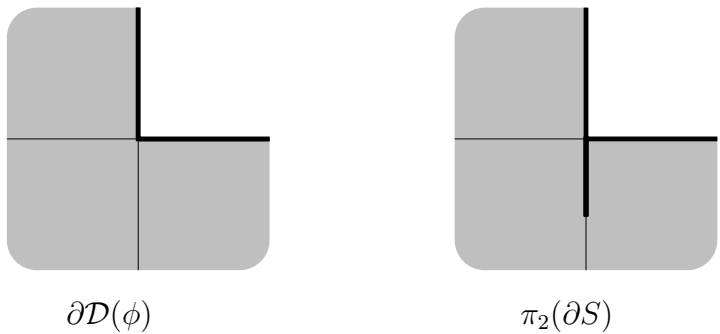
\caption{\label{Fig:InnerDiff} A typical example in which 
 \( \partial \DDD (\phi) \) and  \( \pi_2 (\partial S)\) differ. Varying
 the length of the ``tail'' in \( \pi_2(\partial S)\),  changes the
 modulus of \(S\). } 
\end{figure}

We refer
to the \(y_i\)'s and \(z_j\)'s collectively  as the
 {\it corner points} of \( \DDD (\phi) \). If \(y_i = z_j\) for some
 \(i\) and \(j\), the corresponding corner point is {\it degenerate.}
An arc joining a degenerate corner to itself may be degenerate as
well. If both arcs are degenerate, we say that the corner point is
{\it isolated}. For example, the corner points of  the trivial disk \(
\theta \in \pi_2(\bfy, \bfy) \) are all isolated. 

It is worth pointing out  that even when \(u\) is holomorphic, 
the singular chain \( \partial S\) need not be equal to the cellular
chain \( \partial \DDD (\phi) \), but only homologous to it. 
 This fact will be relevant when we study the Maslov index,
since varying the singular chain will change the complex structure on
\(S\). A typical example of this phenomenon
 is shown in Figure~\ref{Fig:InnerDiff}.

\subsection{ The Maslov index}
\label{Subsec:Maslov}

In this section, we describe how to compute
 \( \mu (\phi ) \) from the domain of \( \phi \). We state the basic
 formula here, but postpone its proof to the next section. 

\begin{figure}
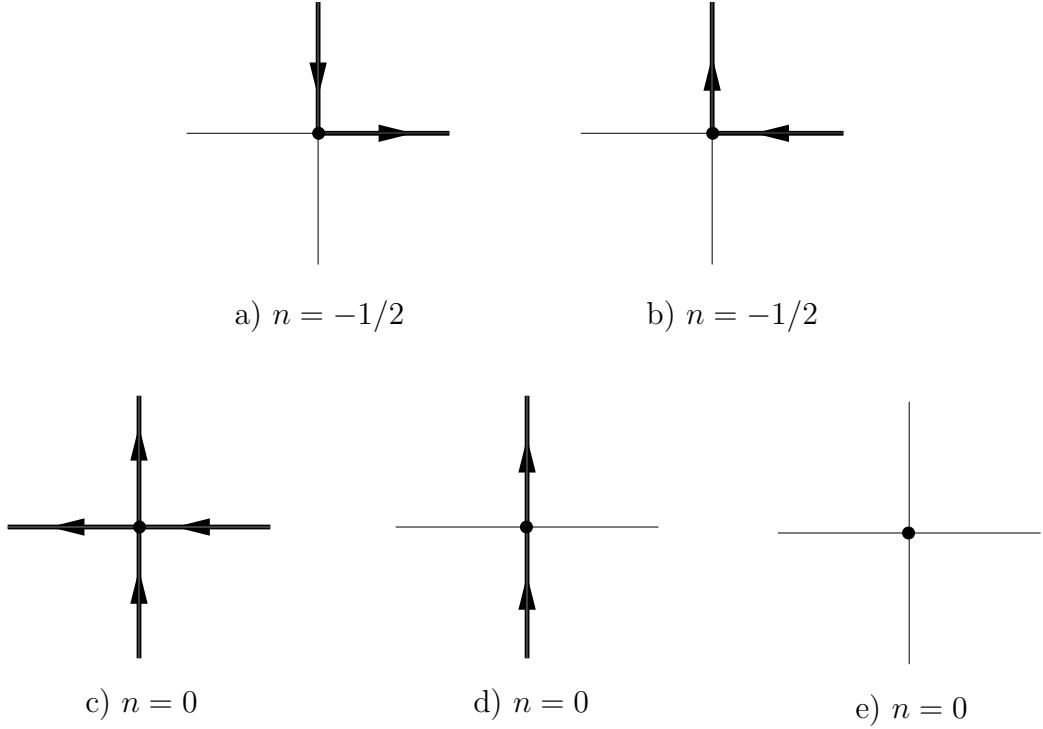
\caption{\label{Fig:LocalMult} The local multiplicities for the
  various possible corner points. The two intersecting lines are \(
  \alpha _i \) and \( \beta _j \); \( \partial \DDD (\phi) \) is shown
  in bold. Types c), d), and e) are degenerate;
  d) is isolated.} 
\end{figure}

\begin{thrm}
\label{Thm:Maslov}
\begin{equation}
\label{Eq:Maslov}
\mu (\phi) = 2\chi(\DDD(\phi)) + \phi \cdot \Delta  + \sum
_{i=1} ^{2g} n(c_i).
\end{equation} 
\end{thrm}

 The terms in this formula are to be interpreted as
follows. \( \chi (\DDD (\phi)) \) is the Euler characteristic of
\(\DDD (\phi) \) as {\it a cellular chain}. In other words, it is \(
\chi (T) \), where \(T\) is a surface with boundary \( \partial \DDD
(\phi) \) which maps to \( \Sigma \) as an unbranched cover with
multiplicities given by \( \DDD (\phi) \). (This definition extends to
domains with negative multiplicities, via the relation
\( \chi (\DDD(\phi ) + \Sigma ) = \chi (\DDD (\phi) ) + 2- 2g \).) 
Note that \( \partial T \) need not map to isolated corner points of
\( \phi \). For example, if \(\theta \in \pi_2 (\bfy, \bfy) \) is the
trivial  disk, \(T\) is the empty set. 
The second term is the intersection number of \( \phi \) with the
diagonal \( \Delta \) in \( s^g\SSS_g \). Since  \(\TTT_\alpha \) and
\( \TTT _\beta \) are disjoint from \( \Delta \), this number is well
defined. Finally, the numbers \(n(c_i) \) are certain {\it local
  multiplicities} of the \( 2g \) corner points of \( \DDD (\phi)
\). These numbers depend only on the geometry of \( \partial \DDD
(\phi) \)  in a
neighborhood of the corner point. The local multiplicities of the
different possible corners 
are shown in Figure~\ref{Fig:LocalMult}. Note that the local
multiplicity of any degenerate corner point is \(0\).

The first and last terms in Equation~\ref{Eq:Maslov} are easily
computed from \( \DDD (\phi) \), but the middle term is not as
obvious. We explain how to evaluate it  in 
 the special case when \( \partial \DDD (\phi) \)
is {\it planar}, that is, when we can surger \( \SSS _g \) along \(g\)
circles disjoint from \( \partial \DDD (\phi) \) to get a 
 sphere. This procedure associates
to each domain \( \DDD (\phi) \) in \(\Sigma \) a new domain 
\( \DDD  (\phi') \) in \(S^2\), together with
\(2g\) marked points \(p_1^\pm, p_2^\pm, \ldots, p_{g}^\pm \) (the
 centers of the surgery disks.) It is not difficult to see that 
\(\DDD (\phi ') \) determines a  homotopy class of disks \(\phi ' \)
 in \( s^g S^2 \cong {\bf CP}^g \).

Let us consider the simpler problem of finding the intersection
number of \( \phi '\) with the diagonal \( \Delta_{S^2} \) in \( s^g S^2 \). 
Suppose there is a point \(w \in S^2 \) with \( n_w(\phi ') =0 \). We
put \( w \) at infinity, so that \(\phi \) can be represented by
a disk in \(s^g \C = \C^g \). In this case, we can use  the
well-known identification of \( \pi_1(\C^g - \Delta_{\C}) \) 
with the braid group \( B_g \). We briefly recall this correspondance here.
 Suppose we are given a  loop in 
\( \C^n - \Delta_{\C} \). By the usual incidence construction, we get
 maps \( \pi_1 \colon C \to S^1 \) and \( \pi_2 \colon C \to \C \), where
 \(\pi_1 \) is a \(g\)-fold covering map. Let  \(\pi_1^{-1}(1) = \{a_1,
\ldots, a_g\}\). If
\(g \colon [0,1] \to S^1 \) is the standard generator of \(\pi_1(S^1) \) and
\(\tilde{g}_i \colon [0,1] \to C \) is the lift of \(g\) which starts at
\(a_i\), we set \( \gamma _i = \pi_2 \circ g_i \).
Then the  curves \( \Gamma_i (t) = (\gamma _i (t), t)) \) define the
associated  braid in \( \C \times \R \cong \R^3\).

 To apply this construction to \( \partial \DDD(\phi ') \), we
  take  \( \gamma _i \) to be the curve which starts at 
\(y_i \), travels along an \( \alpha\)-segment of \( \partial \DDD
(\phi')  \) until it reaches some \(z_j\), and then
continues along a \( \beta \)-curve to  \(y_k\). We denote the
resulting braid by \( \beta (\phi') \).

\begin{figure}
\input{figs/braid.pstex_t}
\caption{\label{Fig:BigelowCalc} \( \DDD (\phi) \) and \(\beta (\phi)
  \). }
\end{figure}
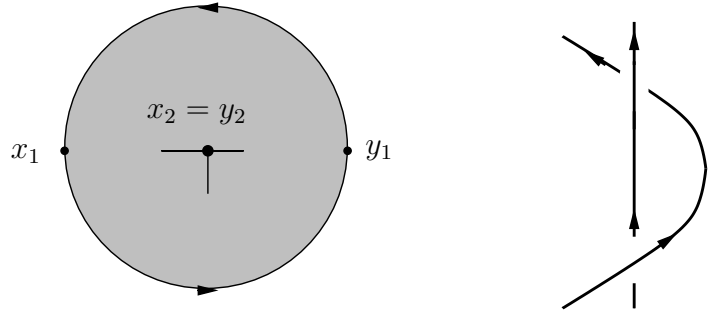

It is easy to see that the abelianization of \(B_g \) is \(\Z \) and
that the quotient homomorphism \(t\colon B_g \to \Z \) is the map which
sends a braid to its writhe. On the other hand, the map \( \pi_1 (\C^g
- \Delta _\C)  \to \Z \) which sends a loop to its linking number with
\( \Delta \) factors through the abelianization and has \(1\) in its
image. (Just take a small loop around a smooth point of \( \Delta_\C
\).) Thus it agrees with \(t\) up to a factor of \( \pm1\).
Checking the sign, we see that  \( \phi ' \cdot \Delta _\C
= t(\beta (\phi')) \). 
\vskip0.05in
\noindent{\bf Example:} Suppose \( \DDD(\phi ') \) has the form shown in
Figure~\ref{Fig:BigelowCalc}a. The corresponding braid has writhe two,
so \(S\) is a double branched cover of \(D^2\) branched at two points;
{\it i.e.} an annulus. The map \( \pi_2\) pinches one component of
\( \partial S\) down to a T--shaped loop supported on \( \alpha \)
and \( \beta \).  
\vskip0.05in

We now drop the hypothesis that \(n_w(\phi') = 0 \) and put an
arbitrary point of \(S^2 - \alpha - \beta \) at infinity. Let \( L \in
\pi_2 ({\bf CP}^g ) \) be the class of a line. Then \(n_w (L) = 1\)
for any \(w\) and \( L \cdot \Delta _{S^2} = 2g-2 \) 
 (since \( \DDD (L) = S^2 \) is a branched \(g\)-fold
cover of \(S^2\)).  The connected
sum \( \phi_1 ' = \phi ' \# (-n_{\infty}(\phi')L ) \) has
\(n_\infty(\phi_1') = 0 \), so
\begin{align*}
\phi_1' \cdot \Delta_{S^2} & = t(\beta(\phi_1')) = t(\beta(\phi '))
 \quad \text{and} \\
 \phi ' \cdot \Delta
_{S^2} & = (2g-2) n_\infty (\phi') + t(\beta (\phi')).
\end{align*} 

Finally, we return to the problem of determining \( \phi \cdot \Delta
\). Let \(S \subset (\SSS)^g\) and \(S'\subset (S^2)^g \)  be the
surfaces determined by \(\phi \) and \( \phi ' \). We would like to
understand the relationship between them. If
\begin{equation*}
D_i^\pm =  \{(x_1, \ldots x_g) \subset (S^2)^g \ts | \ts x_j = p_i^{\pm}
\text{ for some } j \},
\end{equation*}
 a generic \(S'\) will intersect \(D_i^\pm \) transversely at smooth
points of \(D_i^\pm\). The surface \(S\) is by obtained by 
embedding \( (S^2)^g
- \cup D_i^\pm \) in \( (\Sigma _g)^g \) in the obvious way. The image
of \( S' \) under this embedding   will be a punctured surface.
 Each puncture created by   \(D_i^+
\) will have a corresponding puncture created by 
\(D_i^-\), and  \(S\)
is obtained by filling in each pair of  punctures with a copy of the
appropriate annulus in \( \SSS_g \). 
To compute the effect of this procedure on \( \phi \cdot \Delta \),
 note that  each  added annulus results in
  two more branch points of \( \pi _1 \), and thus in an additional
two intersections with the diagonal. To determine the sign of these
intersections, we use the following lemma, whose proof is a
straightforward calculation in local coordinates. 

\begin{lem}
If \( \pi_1 \: S \to D^2 \) is an orientation preserving branched
cover and \(x \in S\) is a branch point of \( \pi_1 \), the sign of the
corresponding intersection with \( \Delta \) is the local degree of \( \pi_2
\) at \(x\). 
\end{lem} 

It follows that the total contribution to the intersection number
coming from \(V_i\) is twice the multiplicity of \(\DDD (\phi)
\) over \(V_i\), which is in turn equal \( n_{p_i^+}( \DDD (\phi '))
+ n_{p_i^-}(\DDD (\phi ') )\). 

To summarize, we see that we have proved the following

\begin{prop}
If \( \partial \DDD (\phi) \) is a planar curve in \( \SSS_g \) and
\(\DDD (\phi ')\) is its associated domain in \(S^2\), then  
\begin{equation*}
\phi \cdot \Delta = (2g-2) n_\infty (\DDD (\phi')) +t(\beta (\phi')) +
 \sum_{i=1}^{g} [ n_{p_i}^+  (\phi') +   n_{p_i}^-(\phi ')]. 
\end{equation*} 
\end{prop}

\noindent{\bf Remark:} When \( \partial \DDD (\phi)
 \) is not a planar curve, it should
be possible to evaluate  \(\phi \cdot \Delta \) by thinking
about the braid defined by \( \phi \) in the braid group of the
punctured surface \( \SSS_g - pt \). The abelianization of this group
is \( \Z \oplus \Z ^{2g} \), with the first summand corresponding to
the linking number  of \( \partial \phi \) 
with \( \Delta\) and the second summand
corresponding to the class of \( \partial \phi \) in \( H_1
(s^g\SSS_g) \cong H_1 (\SSS_g) \). Thus the problem is to find
an approriate generalization of the writhe for the braid group of a
 punctured surface.

\subsection{Proof of Theorem~\ref{Thm:Maslov}}

Before starting the proof, we give an alternate formulation of
equation~\ref{Eq:Maslov} together with a heuristic argument explaining
why this formula should be true. 

\begin{figure}
\input{figs/LocalMult2.pstex_t}
\caption{\label{Fig:NuMult} Possible values of \( \nu \). Darker
  shading indicates regions of greater multiplicity. }
\end{figure}
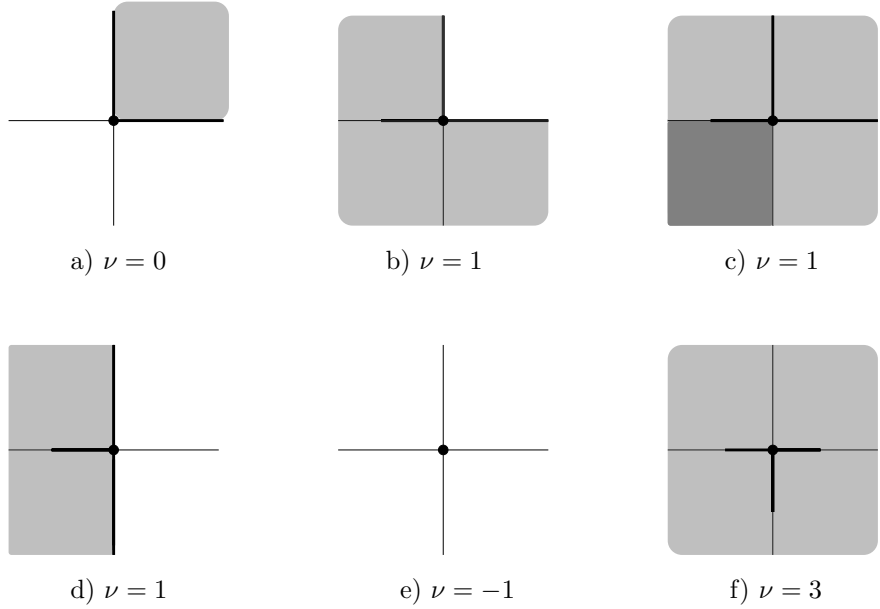

Suppose for the moment that \( \phi \) is actually represented by a
holomorphic map \(u\), so that the maps \( \pi_1 \: S \to D^2 \) and 
 \( \pi_2 \: S \to \Sigma _g \)  are branched covers with \( b_1 = g - \chi
 (S)\) and \(b_2 \) branch points. The surface \(S\) is
equipped with two pieces of data: a complex structure and a system of
\(2g \) marked points on its boundary. We can reconstruct this data
either from the branched cover over the disk or from the branched
cover over \( \Sigma _g \): in order for there to be a holomorphic
disk,  the two reconstructions must agree. Over the disk, \(S\) is
determined by the location of the \(b_1 \) branch
points. Over \( \Sigma _g \), we can vary the location of \(b_2\)
branch points, but we may also be able to vary the complex structure
on \(S\) by
``cutting'' the boundary of \( \DDD (\phi) \) near the corner
points, as illustrated in Figure~\ref{Fig:InnerDiff}. Each direction we
can cut in gives us an additional real parameter. 
We denote the number of parameters we get at the corner \(c_i\) by \( \nu
(c_i)\).  The possible values of
\(\nu\) are shown in Figure~\ref{Fig:NuMult}. 
 Thus we have a total of \(2 b_1 + 2 b_2 + \sum \nu (c_i)
\) free parameters. (Each degenerate corner  contributes only
once to the sum, not twice.) On the other hand, the space of complex 
structures
on \(S\) has dimension \(-3 \chi (S) \), and the marked points are
determined by another \(2g \) real parameters. Subtracting, we get

\begin{align*}
\mu (\phi) & = 2 (g- \chi (S))  + 2b_2 + \sum  \nu (c_i) -(2g- 3
\chi(S)) \\
& = \chi (S) + 2b_2 + \sum  \nu (c_i).
\end{align*}

To see that this is equivalent to  equation~\ref{Eq:Maslov}, 
 we would like to write \(
 \chi (S) = b_2 +  \chi (\DDD(\phi)) \).
This is not quite correct, however, since boundary components of \(S\)
will collapse to isolated corner points. We can account for this fact
 by adding a correction term \( \sum m(c_i) \) to the right-hand
 side, where \(m\) is \(+1 \) for a corner point of type {\it e)} (which
 results in an extra disk, and \(-1\) for a point of type {\it f)} (which
 results in an extra puncture.) Then we have
\begin{align*}
\chi (S) + 2b_2 + \sum  \nu (c_i) &= 2 \chi (D) - \chi (S) + 
 \sum [ \nu (c_i) + 2m(c_i)] \\ &= 2 \chi (D) + \phi \cdot \Delta - g + 
\sum  [\nu (c_i) + 2m(c_i)] \\ &= 2 \chi (D) + \phi \cdot \Delta + 
\sum n(c_i)
\end{align*}
where in the last step we have distributed the factor of \(-g \) over
the sum, adding a \(-1/2\) to each ordinary corner and a \(-1\) to
each degenerate corner. We leave it to the reader to check that the
resulting contributions agree with the values of \(n(c_i) \) given
in section~\ref{Subsec:Maslov}. 

To actually prove the theorem, we attempt to
 represent \( \phi \) by a map \(u\)
which is similar to a holomorphic map, in the sense that 
\( \pi _1 \) and \( \pi _2 \) are orientation preserving branched
covers. We will not always achieve this goal, but we can control the
situations in which it fails to hold. 

\begin{lem}
If equation~\ref{Eq:Maslov} holds for those \( \phi \) with
\( \DDD (\phi) \geq 0 \), it
holds for all \( \phi \). 
\end{lem}

\begin{proof}
Let \(E \in
\pi_2(s^g\SSS) \) be the positive generator. Then \( \DDD (E) =
\SSS \) is a \(g\)--fold cover of \(S^2 \) branched at \( E \cdot \Delta \)
 points, so   \(E \cdot \Delta = 4g - 2\).
 Given \( \phi \in \pi_2(\bfy, \bfz) \), consider \( \phi '  = 
\phi \# n E \). Then \( \DDD (\phi ' ) = \DDD (\phi) + n \SSS \geq 0
\) for \(n \gg 0 \), so
\begin{align*}
\mu (\phi ' ) &= 2 \chi (\DDD (\phi') ) + \phi ' \cdot \Delta + \sum
n(c_i) \\
 \mu (\phi  ) + 2n &= 2 \chi (\DDD (\phi)+ n(2-2g) ) + \phi  \cdot
 \Delta + n(4g-2) + \sum
n(c_i) \\
\mu (\phi  ) &= 2 \chi (\DDD (\phi) ) + \phi  \cdot \Delta + \sum
n(c_i). 
\end{align*}
\end{proof}

\begin{defn} Suppose we are given \(\pi_2 \colon S \to \SSS_g\). A collapsed
  tube of \(\pi_2 \) is an annulus \( A\cong S^1\times[-1,1] \subset S
  \) such that \(\pi_2|_A\) factors through the projection to \([-1,1]
  \). 
\end{defn}

\begin{lem}
\label{Lem:Nodal}
If \( \DDD (\phi ) \geq 0 \), then there is a representative \(u\) of
\(\phi \) for which \(S\) is a  nodal surface, \( \pi _1 \colon
S \to D^2 \) is an orientation preserving  branched cover, and 
 \( \pi _2 \colon S \to \SSS _g \) is an orientation preserving
 branched cover away from a finite  union of collapsed tubes.
\end{lem}

\begin{proof}
Our first step is to construct a
branched covering map \(p_2\colon T  \to \DDD (\phi) \) such that \(
\partial T\) maps to the path \( \partial \DDD (\phi) \).
To do so, we start with \(n_{z_i} (\phi) \) copies of the cell \( \mathcal{C}
_i \), where \(z_i \in   \mathcal{C}_i \). We glue these cells
together along their edges in a maximal fashion, so that when we are
done, there are no unglued edges which project to the same edge in \(
\SSS \) with opposite orientation. (Of course, there are many such
gluings; we just choose one.). Thus we obtain a surface \(T_0\) with a
projection \(p_0\: T_0 \to \SSS \). The map \(p_0 \) is a covering map
except at the vertices of the cellulation, where it may be branched. 
Since the gluing was maximal, \( p_0( \partial T_0) \) is equal to \(
\partial \DDD (\phi) \) as a chain. The two may not be equal as paths,
however, since  \( p_0( \partial T_0) \) may ``turn a corner'' at a
vertex where  \(\partial \DDD (\phi) \) goes straight.
 To solve this problem, we slit the two components of \(T_0 \) near
 the vertex in question
 and glue them back in the opposite way. This introduces
 another branch point in the projection. After making these
 modifications, we obtain the desired surface \(T\) and map \( p_2\: T
 \to \SSS \). 

Next, we  let \(T'\) be a connected surface obtained by
  taking the connected sum of all components of \(T\). Then we can find
  a map \(p_2' \: T' \to \SSS_g \) which agrees with \(p_2\) away from
  the connected sum regions, and for which each connected sum region
  is a collapsed tube. Since \(T'\) is connected, there is
 a \(g\)-fold covering
  map \(p_1'\: T' \to D^2 \) which has precisely two branch points in
  each collapsed tube. We immerse \(T' \) in \(D^2
  \times \SSS _g \) via the map \(p_1' \times p_2' \). When \( p
  _1' \) and \( p_2' \) are generic, the image will have isolated
  double points, none of which are contained in the collapsed
  tubes. We define \(S\) to be this image and \( \pi_1\) and \(
  \pi_2 \) to be the projections onto the two factors. With these maps,
   the hypotheses of the lemma are clearly satisfied. 
\end{proof}

\begin{proof}(of Theorem~\ref{Thm:Maslov}) 
The Maslov index of \(\phi \) is defined as follows
\cite{Viterbo}. Choose a complex structure \(J\) on \(s^g\SSS\).
 Let \(p\:I^2 \to D^2\) be the map from unit square to
the disk which collapses the two vertical sides, and choose any map 
\(u\: D^2 \to s^g\SSS \) which represents \(\phi \). 
 Then \(E = (up)^*Ts^g\SSS \) is a
(necessarily trivial) bundle of complex vector spaces over the
 square. We define a family  \(L\) of Lagrangian subspaces on \(
 E|_{\partial I^2}\) by pulling back the tangent spaces of \(
\TTT_\alpha \) and \( \TTT _\beta \). To be precise, we 
identify the tangent spaces
at \( \bfy \) and \( \bfz \) with \(\C^g\) in such a way that 
\(T\TTT_\alpha = \R^n \) and \( T \TTT_\beta = i \R^n \), and set
\begin{align*}
L_{(s,0)} & = (up)^*T\TTT_\alpha \\
L_{(s,1)} & = (up)^*T\TTT_\beta  \\
L_{(*,t)} & = e^{i\pi t/2} \R ^n \quad (*=0,1).
\end{align*}
If we choose a trivialization \(E \cong I^2 \times \C ^g \), 
 \(L\)  defines a loop in the space  \( \Lambda _g \) of totally real
subspaces of \(\C ^g \), and \( \mu (\phi) \) is the class of this loop in
\( \pi _1 (\Lambda_g ) \cong \Z \). Equivalently, it is the obstruction to
extending \(L\) to a family of totally real subspaces on all of
 \(I^2\). (Note that this definition differs from the usual one in
 that we  use totally real spaces rather than
 Lagrangians. The reason  is that 
the symplectic form on \( s^g \SSS_g \) does not extend across the
 diagonal.  Since the space of totally real subspaces deformation
 retracts to the space of Lagrangians, this has no effect
 on the definition.)

Without loss of generality, we may assume \( \DDD (\phi) \geq 0 \) and
that the map \(u\) is of the form described in Lemma~\ref{Lem:Nodal}. 
To compute \( \mu(\phi) \), we consider a slightly different family  of
Lagrangians on the square. To define this family, we  round off the
corners of \(S\) to get a smooth manifold \(S'\). Then 
 \(T(\partial S') \) is a natural family of Lagrangian subspaces
of \(TS'|_{\partial S'} \). \(L'\) is defined to be the push-forward of
this family. More precisely, we choose \(\pi_1'\:S' \to I^2\) and 
\( \pi_2'\:S' \to \SSS_g\) which
agree with \(\pi_1\) and \( \pi_2 \) except near the corners. Then if
\(x_1,\ldots, x_g\) are the preimages of \((s,t) \in \partial I^2 \)
under \( \pi_1 '\), we can identify \(E_{(s,t)}\) with \((T(\SSS)^g)
_{(x_1,x_2,\ldots,x_g)} \). Under this identification, 
\(L'_{(s,t)} \) is  the space spanned by the vectors \((0,\ldots,d\pi_2
(T(\partial S')_{x_i}),\ldots, 0)\), where the nonzero entry is in
the \(i\)th position.

By construction, \(L'\) and \(L\) agree on the horizontal sides of
the square, but they may differ on the vertical sides. To compute 
 \( \mu (L) - \mu (L') \), we note that on the vertical
sides, \(L\) and \(L'\) decompose into  a direct sum in the sum of the
tangent spaces near  each corner: \(L = \oplus_{c_i} L(c_i) \) and
 \(L' = \oplus 
_{c_i} L'(c_i) \), where \(L_{c_i} \) is one dimensional if \(c_i\) is
an ordinary corner and two-dimensional if \(c_i\) is degenerate.
Thus we can write
\begin{equation*}
 \mu (L) - \mu (L') = \sum_{c_i} \mu (L(c_i)) - \mu (L'(c_i)) =
\sum_{c_i} n(c_i)
\end{equation*}
where the numbers \(n(c_i)\) depend only on the local geometry of the
corner \(c_i\). 

We now concentrate on finding \(\mu (L') \). To do so, we try
to write down an extension of \(L'\) to the entire square, and then
study the places where this extension fails to exist. 
The induced orientation on \(
\partial S'\) determines a nonvanishing vector field \({\bf v}\) on 
\( \partial S'\). We pick an extension of 
\({\bf v} \) to \(S'\)  satisfying the following
conditions. Near each collapsed tube, \({\bf v} \) should be
chosen so that it can be extended to a nonzero vector field on the
surface obtained by surgering the tube. (Thus \( {\bf v} \) defines a
vector field on the surface \(T\) of Lemma~\ref{Lem:Nodal}.) Elsewhere,
\( {\bf v} \) should vanish at isolated points, none of which are nodes of
\(S\). There will be \( \chi (T)\) such points, counted with
sign.  Away from
these points, \({\bf v} \) determines a family of Lagrangians on
\(S'\). We try to extend \(L'\) to the entire square by pushing this
family of Lagrangians forward in the same way we defined \(
L'\). There are four sorts of places where this extension might fail
to exist: 
\begin{enumerate}
 \item Branch points of \( \pi_1\) not on the nodes or collapsed
 tubes. 
\item Branch points of \( \pi_2 \) not on the nodes or collapsed tubes.
\item Nodal points of \(S\).
\item The collapsed tubes. 
\item  Points  where {\bf v} vanishes. 
\end{enumerate}

Without loss of generality,  we can assume that all branch points of
type {\it 1)}  and {\it 2)}
  are simple and that the images of each collapsed tube
is supported in a small neighborhood of some point in the
square. Moreover, we can assume that 
 all these points have disjoint images under \(
\pi_1\).  Let \(\{x_i\} \subset I^2\) 
be the set of these images, and let \( \gamma _i \) be a small circle
centered at \(x_i\). Then
\begin{equation*}
\mu(L'|_{\partial I^2}) = \sum \mu (L'|_{\gamma_i}).
\end{equation*}

Computing  \( \mu (L'|_{\gamma_i}) \) is  a local problem, which we
can study individually for each of the types listed above.

\begin{lem}
\label{Lem:Type1}
 Suppose
 that \(x_i\) is a point of type {\it 1)}. Then \( \mu
 (L'|_{\gamma_i}) = 1 \).
\end{lem}

\begin{proof}
Let  \(U_i\)
be a small disk centered at \(x_i\) and containing \( \gamma _i\). 
 Without loss of
generality, we can choose \( {\bf v} \) to be a constant vector field
on each component of \( \pi_1^{-1} (U_i)\).
 Then \(L'|_{\gamma_i} = L_1 \oplus L_2 \), where \(L_1\) is a
 \(2\)--dimensional Lagrangian coming from the component with the branch
 point, and \(L_2\) is a constant Lagrangian coming from the other
 \(g-2\) components. Thus  \(\mu(L'|_{\gamma_i}) = \mu (L_1) \).

 To evaluate \( \mu (L_1) \),
 we use the following local model: \(\gamma_i \) is the unit circle in
 \(\C\), \(\pi_1\: D^2 \to D^2 \) is the map \( z \to z^2 \), and
\(\pi_2 \: D^2 \to \C \) is the identity map. The
 corresponding \(u\: D^2 \to s^2 \C \) sends \( z \to
 \{\sqrt{z},-\sqrt{z}\}
 \). Identifying \( s^2 \C \) with \( \C^2 \) by the map which sends
 \(\{u,v\} \to (u+v,uv) \), we see that  \(u\: D^2 \to \C^2 \) is given by \(
 u(z) = (0,-z) \).  Then a tangent vector \(\alpha \) at the point \(
\sqrt{z} \in S \) pushes forward to the tangent vector \( (\alpha, -\sqrt{z}
 \alpha) \in T\C^2 |_{(0,-z)} \). Since \({\bf v} \) is a constant
 vector field, it  pushes forward to define 
\begin{equation*}
L_1 = \{ ((\alpha + \beta){\bf v}, (\beta - \alpha) \sqrt{z}{\bf v}) \subset
u^*(T\C^2)|_{z} \ | \ts \alpha, \beta \in \R \}.
\end{equation*}
Thus \( L_1 \) is the direct sum of a trivial family of Lagrangians in
the first factor with a family with Maslov index 1 in the second
factor, so \( \mu (L_1) = 1 \). 
\end{proof}

\begin{lem}
 Suppose
 that \(x_i\) is a point of type {\it 2)}. Then \( \mu
 (L'|_{\gamma_i}) = 2 \).
\end{lem}

\begin{proof}
As in the previous lemma, we choose \({\bf v} \) to be constant on
each component of \( \pi_1 ^{-1} (U_i) \). Then \( L'|_{\gamma_i} \)
splits as a direct sum, and the only contribution to the Maslov index
will come from the summand containing the branch point. In this case,
there is a local model in which \(\pi_1 \: D^2 \to D^2 \) is the
identity map and \( \pi_2 \: D^2 \to D^2 \) is the map \(z \to z^2
\). Then \(u\:D^2 \to D^2\) is also given by \( z \to z^2 \), and the
constant vector field \({\bf v} \) pushes forward to the vector field
\(z{\bf v} \). It is easy to see that the Maslov index of the
associated family of Lagrangians is 2. 
\end{proof}

\begin{lem}
\label{Lem:Type2}
If \(x_i\) is a point of type {\it 3)}, then {\em \( \mu(L'|_{\gamma _i})
= 2 \ts  \text{sign} \ts x_i \)}, where {\em \( \text{sign} \ts
 x_i \)} is the sign of the
intersection in \(D^2 \times \SSS_g \). 
\end{lem}

\begin{proof}
As before, we see that \(
\mu(L'|_{\gamma _i}) =  \mu(L_1) \), where \(L_1\) is a
two-dimensional family of Lagrangians pushed forward from the two
components of \(S\). 
Near \(x_i\), a local model for the two sheets of \(S\) is given by
\begin{equation*}
S=\{(\bfx, A \bfx ) \ts | \ts \bfx \in \R^2\} \cup
\{(\bfx, B \bfx )\ts | \ts \bfx \in \R^2\}
\end{equation*}
where \(A\) and \(B\) are invertible \(2\times2\) matrices.
The maps \(\pi_1\) and \( \pi_2 \) are the projections onto the first
and second coordinates, and the  map 
\( u\: D^2 \to \C^2 \) is  given by \( u(\bfx) = (A\bfx + B
\bfx, (A\bfx)(B\bfx)) \), where the parentheses indicate complex
multiplication under the identification \(\R^2 = \C \). If we choose
\({ \bf v}\)  to be given by constant vectors \({\bf v}_1 \) and 
\({\bf v}_2 \) on the two components of \(S\), we see that 
\begin{equation*}
L_1 = \{(\alpha A{\bf v}_1 + \beta B {\bf v}_2, \alpha (A{\bf v}_1) (B
\bfx) + \beta (B{\bf v}_2) (A \bfx) ) \ts | \ts \alpha, \beta \in \R\}.
\end{equation*}
To evaluate \( \mu (L_1)\), we choose \({\bf v}_1 \) and 
 \({\bf v}_2 \) so that \( A{\bf v}_1 = B{\bf v}_2 = {\bf w} \). Then
 \(L_1 \) becomes 
\begin{equation*}
\{ ((\alpha + \beta){\bf w}, (\alpha A\bfx + \beta B \bfx ){\bf w}) \}.
\end{equation*}
This family is totally real at all points on \(
\gamma_i \)  whenever \( A-B \) is nonsingular, so we can homotope
 \(L_1\) to a similar family \(L_1'\) for which  \(A' = A-B\) and 
\(B' = 0\). It follows that  \( \mu (L_1) = \mu (L_1')  = 2 \det (A-B) \). 

On the other hand, the sign of \(x_i\) as an intersection point is
 the sign of 
\begin{equation*}
\det \left( \begin{array}{cc} I & A \\ I & B  \end{array}\right) = \det (A-B).
\end{equation*}
This proves the claim.
\end{proof}

\begin{lem}
If \(x_i\) is a point of type {\it 4)}, then \( \mu(L'|_{\gamma _i})=0\).
\end{lem}

\begin{proof}
If \(u\) is a representative for \(\phi \), we can always make another
representative \(u'\) of \(\phi \) by adding a collapsed tube to
\(S\). \(u'\) and \(u\) have the same points of types {\it 1)--3)},
and they obviously have the same Maslov index, so adding a collapsed
tube must have no effect on \( \mu \). 
\end{proof}

\begin{lem}
If \(x_i\) is a point of type {\it 5)}, then {\em \( \mu(L'|_{\gamma _i})=
2 \ts \text{ind} _{x_i} {\bf v} \)}.
\end{lem}

\begin{proof}
As in Lemmas~\ref{Lem:Type1} and \ref{Lem:Type2},
 we can decompose \( L'|_{\gamma_i} = L_1
\oplus L_2 \), where \(L_1\) is a constant Lagrangian and \(L_2\) is
the one-dimensional family of Lagrangians determined by \({\bf v}\)
near the point at which it vanishes. It is then easy to check that \(
\mu (L_2) = 2 \ts \text{ind} _{x_i} {\bf v}\).
\end{proof}

 Since the sum of \( \text{ind} _{x_i} {\bf v} \) over
all points of type {\it 5)} is \(\chi (T) \), we have 
\begin{align*}
\mu (L'_{\partial_{I^2}}) & = b_1+2b_2 + 2n+ 2 \chi (T).
\end{align*}
where \( b_1\)  and \(b_2\)
are the number of points of types {\it 1)} and {\it 2)} and \(n\) is
the signed number of nodes of \(S\). We claim that \( b_1 + 2n = \phi
\cdot \Delta.\) To see this, recall that intersections
 of \( \im u \) with \( \Delta \)
correspond to branch points of \( \pi_1 \), and the sign of the
intersection is the sign of the map \( \pi_2 \). By construction,
branch points of type {\it 1) } have multiplicity one and positive
sign. Next a local calculation like that of Lemma~\ref{Lem:Type2}
shows that each node of \(S\) corresponds to a point of tangency with \(
\Delta \), and the sign of the intersection corresponds to the sign of
the intersection at the node. Thus the nodes contribute a total of \(2
n \) to the intersection number. Finally, the two
branch points on each collapsed tube contribute with opposite
signs. To see this, we perturb the map \( \pi_2 \) a bit, so the
collapsed tube becomes a flattened tube, with the two branch points on
opposite sides. This proves the claim. 

Next, we observe (as  in the heuristic argument) that \( b_2 + \chi (T) =
\chi (\DDD (\phi) ) \) up to correction factors of \( \pm 1 \)
introduced by isolated corners. We lump these factors in with the
\(n(c_i) \) to obtain a formula like that of
equation~\ref{Eq:Maslov}.

 Finally, we must determine the numbers \(n(c_i)\). For nondegenerate
 corners, this is a direct computation with one dimensional
 Lagrangians. For corners of type {\it a)}, \(n(c) = -1 \) if \(c\) is
 one of the \(y_i\) and \(0\) if \(c\) is one of the
 \(z_i\). Similarly, \(n(c) \) is either \(0\) or \(1\) for corners of
 type {\it b)}. Since the total number of \(y_i\)'s and \(z_i\)'s is
 equal, we can symmetrize these figures to arrive at the values shown
 in Figure~\ref{Fig:LocalMult}. Corners of type {\it c)} are direct sums
 of a corner of type {\it a)} and a corner of type {\it b)}, so they
 have \(n=0\). Finally, we evaluate \(n \) for corners of type {\it
 d)} and {\it e)} by applying the formula to domains for which \(\mu
 \) is known.
\end{proof}

\subsection{The Localization Principle}

We now turn our attention to the problem of determining whether a
given domain admits a differential.
There are two basic principles to keep in mind here. First, \(\DDD
(\phi)\) must be positive  whenever \( \# \MMM (\phi) \neq 0 \). 
Second, \( \# \MMM (\phi) \) is a {\it local} quantity: it
depends only on the geometry in a small neighborhood of the
domain. More precisely, suppose \( \DDD (\phi) \) is a positive
domain. Define the support of \( \phi \) (written
\(\overline{\DDD} (\phi) \)) to be the union of all cells  with
nonzero multiplicities in \( \DDD (\phi) \) and all  corner points
of \(\DDD (\phi) \). Then we have

\begin{lem}
If \(U\) is any open neighborhood of \( \overline{\DDD} (\phi)\), and
we choose a generic path of almost complex structures \(J_s \) which is
sufficiently close to a complex structure \(J\) induced by a complex
structure  on \( \SSS \), then the
image of any \(J_s\)-holomorphic representative of \( \phi \) is contained in
\(s^gU\). 
\end{lem}

\begin{proof}
 Suppose that the statement
 is false. Then we can choose a sequence \(J^i_s\)
 of generic \(J_s\)'s
converging to \(J\), for which there are \(J_s\)--holomorphic
 representatives \(u^i\)
of \( \phi \) whose image  is not contained in \(s^gU\).
 By passing to a subsequence, we can assume there is a
 sequence of points \(x^i \not \in s^gU\) with \(x_i \in \im u^i \)
 and \(x^i \to x \) for some \( x \not \in s^gU \). Then by Gromov
 compactness, we obtain a (possibly degenerate) \(J\)-holomorphic
 representative \(u\) of  \( \phi \) with \(x \in \im u \). 
In this case, the corresponding map 
\( \pi_2 \: S \to \SSS  \) is holomorphic, so it
is either open and orientation-preserving or constant on each
irreducible component of \(S\).
 The image of each connected component of \(S\)  must contain a corner
point, so  \( \pi_2 (S) = \overline{\DDD} (\phi) \subset U \). But this
contradicts the fact that \(x \not \in s^gU \). 
\end{proof}

As a first application of this {\it localization principle} (\Ozsvath
and Szab{\'o}'s phrase), we show that to study differentials, it suffices
to consider connected domains. More precisely,
suppose that \( \overline{\DDD} (\phi) \) can be divided into two
disjoint domains \(\overline{\DDD}(\phi_1) \) and \( \overline{\DDD}
 (\phi_2) \)
 containing \(g_1\) and \(g_2\) corner points. We let \(
 \TTT_{\alpha_1} \subset s^{g_1}\SSS_g\) be the product of those 
\( \alpha_i \)'s containing corner points of \( \DDD (\phi_1) \), \(\bfy
_1 \in s^{g_1} \SSS_g \) be the set of those \(y_i\)'s in \( \DDD
(\phi _1)\), {\it etc.} Then \( \DDD (\phi_i) \) defines a homotopy class
 \( \phi _i \in \pi_2(\bfy_i, \bfz_i)\) of disks in \(s^{g_i}\SSS_g
 \),  and we have

\begin{prop}
\( \mu (\phi) \cong \mu (\phi_1) + \mu (\phi_2) \) and \( \MMM (
\phi) \cong \MMM (\phi _1 ) \times \MMM (\phi_2) \). 
\end{prop}

\begin{proof}
If \(U_1\) and \(U_2\) are disjoint open sets containing
\(\overline{\DDD}(\phi_1) \) and \( \overline{\DDD}(\phi_1) \), the
localization principle tells us that the image of any holomorphic
representative \(u\) of \( \phi \) is contained in \( s^g(U_1\cup U_2)
\). In fact, it must be in the connected component of this set
containing \( \bfy \) and \( \bfz\), which is \(s^{g_1}U_1 \times
s^{g_2}U_2 \). Similarly, it is not difficult to  see that \( \partial
u \) must actually lie on \( \TTT_{\alpha_1} \times \TTT _{\alpha_2}
\cup \TTT_{\beta_1} \times \TTT _{\beta_2} \). Thus we are really
studying maps of a holomorphic disk into a product. In this case, it
is well-known that the Maslov index and moduli spaces of
\(J\)-holomorphic representatives behave as stated. 

\end{proof}

\begin{cor}
If  \( \mu (\phi) = 1 \) and 
\( \overline{ \DDD } (\phi) = \overline{\DDD} (\phi _1 ) \coprod 
\overline{\DDD} (\phi _2)
\), then \( \MMM (\phi) \) is empty unless one component (say \(
\overline{\DDD}
(\phi_2)\)) is a union of isolated points, in which case \( \MMM (\phi) \cong
\MMM (\phi _1) \).  
\end{cor}

\begin{proof} This is the standard connected sum argument from gauge
  theory. Since \( \mu (\phi) = 1 \), we can assume without loss of
  generality that \( \mu (\phi_2) \leq 0 \). If it is negative, \(
  \MMM (\phi_2) \) is generically empty. If it is \(0\) and \( \MMM
  (\phi_2) \) is nonempty, the \(\R\) action on it is not free, so the
  corresponding map is constant. In this case \( \DDD (\phi_2)
  \) is a union of isolated points and \( \# \MMM (\phi_2) = 1 \). 
\end{proof}

\subsection{Examples of Differentials}
We conclude by  proving the existance of the differentials
used in sections~\ref{Sec:ExactTriangle} and
\ref{Sec:Alternating}. The basic idea is to use the localization
principle and the invariance of \(
\cfr \) to show that the domain in question  must give a differential. 
We start with the ``punctured annular differential'' \( \phi \)
whose domain is shown in Figure~\ref{Fig:ThinDomain}, {\it i.e.}
\( \phi \in \pi_2(\bfy_+, \bfy_-) \), where \(\bfy_\pm = \{x_1,
\ldots, x_n, y_\pm\} \). 

\begin{figure}
\includegraphics{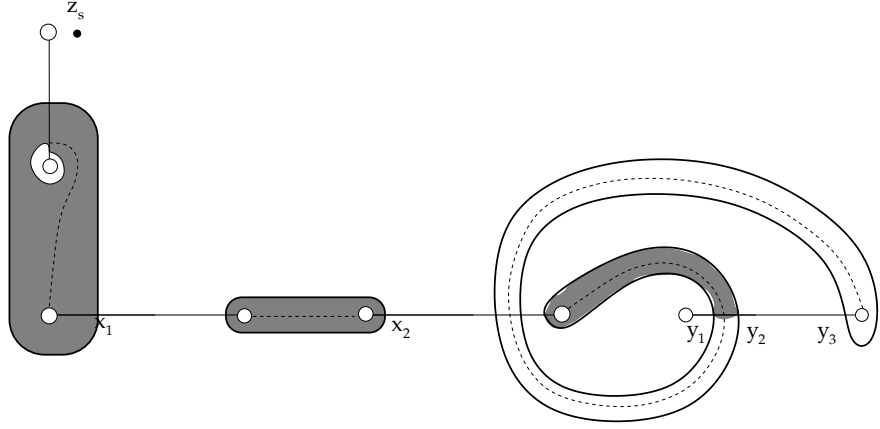}
\caption{\label{Fig:LongTrefoil} A Heegaard splitting of the trefoil
  (\(n=2\)). }
\end{figure}
%\begin{figure}
%\includegraphics{figs/GenAnulus.eps}
%\caption{\label{Fig:ThinDomain2} Domain of the differential from \(
%  \bfy_i^+ \) to \( \bfy_i^-\). } 
%\end{figure}

\begin{lem}
The disk \( \phi \) described above has \(\mu (\phi) = 1 \) and \(\#
\MMM (\phi) = \pm 1 \) for any \(J_s \) sufficiently close to a
complex structure \(J\). 
\end{lem}

\begin{proof}
Consider the Heegaard splitting of the  trefoil shown in
Figure~\ref{Fig:LongTrefoil}. It is easy to see that \( \cfs \) for
this diagram is generated by \(\bfy_i= \{x_1,x_2, \ldots, x_n, y_i\} \),
  and that \(A(y_i) = i\). The domain of the unique \( \psi \in
  \pi_2(\bfy_1, \bfy_2) \) with \(n_{z_s}(\psi) = 0 \) is shaded.
 We claim that \( \mu (\psi) = 1\) and \( \# \MMM (\psi) = \pm 1\).
 When \(n=1\), this is Lemma 8.4 of \cite{OS1}. The result for general
 \(n\) then follows from  the fact that \( \cfr \) is an invariant.
The domain of \( \psi \) is combinatorially identical to the domain of
\( \phi \), except that the latter domain may
be cut by other \( \beta\)--curves, as shown in
Figure~\ref{Fig:ThinDomain}. These
segments are {\it inaccessible}  --- they cannot be connected to the
corner points through a region of positive multiplicity. By the
localization  principle, they cannot be in the image of \(
\partial S\), and are thus irrelevant to the problem of determining \(
\# \MMM (\phi) \). 
By varying the conformal structure of the surface in
Figure~\ref{Fig:LongTrefoil}, we can make \(\DDD (\psi)\) conformally
equivalent to \( \DDD (\phi) \), so  \( \# \MMM (\phi) = \# \MMM (\psi) 
 = \pm 1 \).  

\end{proof}

Next, we consider ``disky differentials.''

\begin{lem}
If the domain of \( \phi \in \pi_2(\bfy, \bfz) \) is an embedded
disk with all interior corners, then \(\mu (\phi) = 1 \) and \( \# \MMM
(\phi) = \pm 1 \). 
\end{lem}

\begin{figure}
\includegraphics{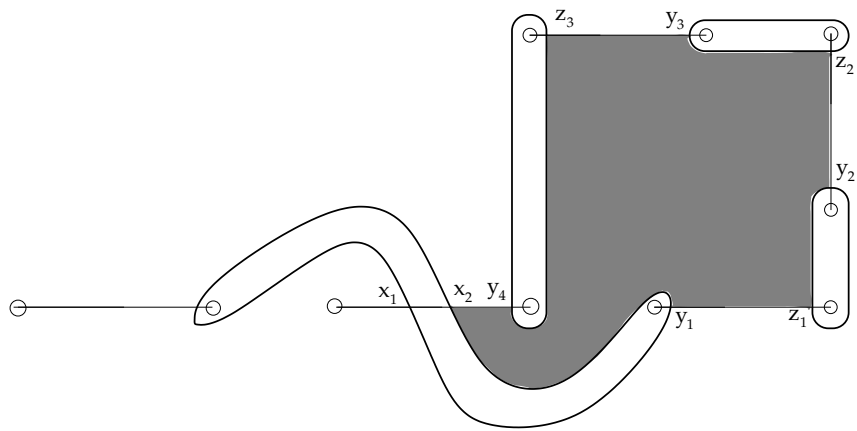}
\caption{\label{Fig:DiskyProof} A disky differential in a diagram for
  the unknot (\(n=4\)).}
\end{figure}

\begin{proof}
Consider the diagram of the unknot shown in
Figure~\ref{Fig:DiskyProof}. The stable complex of this knot has three
generators: \( \bfx_1 = \{x_1,z_1,z_2, \ldots, z_{n-1} \} \), \( \bfx_2 = 
 \{x_2,z_1,z_2, \ldots, z_{n-1} \} \), and \( \bfy = \{y_1,y_2,
 \ldots, y_n\} \). It's easy to see that \(A(\bfx_2) = A(\bfz) \neq
 A(\bfx_1) \), so by the invariance of \(\cfr\), there must be a
 cancelling differential between \( \bfx_2 \) and \(\bfz\).
 The unique \(\psi \in \pi_2(\bfz,\bfx_2) \) with \( n_{z_s}(\psi)
 = 0 \) has the shaded disk in the figure as its domain. Since the
 domain of the corresponding \( \psi ' \in \pi_2(\bfx_2, \bfz)\) is
 negative (in fact, \( \DDD (\psi ') = - \DDD (\psi) \)),
 we must have \( \mu (\psi) =
 1\) and \( \# \MMM (\psi) = \pm 1 \). By varying the metric on \(
 \Sigma \), we can make \( \DDD (\psi) \) conformally equivalent to
\( \DDD (\phi) \), which proves the lemma. 
\end{proof}

%% file: figs/cut.pstex_t
\begin{picture}(0,0)%
\includegraphics{figs/cut.pstex}%
\end{picture}%
\setlength{\unitlength}{1973sp}%
\begingroup\makeatletter\ifx\SetFigFont\undefined%
\gdef\SetFigFont#1#2#3#4#5{%
  \reset@font\fontsize{#1}{#2pt}%
  \fontfamily{#3}\fontseries{#4}\fontshape{#5}%
  \selectfont}%
\fi\endgroup%
\begin{picture}(9088,4225)(3139,-6535)
\put(9601,-6436){\makebox(0,0)[lb]{\smash{\SetFigFont{12}{14.4}{\rmdefault}{\mddefault}{\updefault}{\color[rgb]{0,0,0}\( \pi_2 (\partial S) \)}%
}}}
\put(3901,-6436){\makebox(0,0)[lb]{\smash{\SetFigFont{12}{14.4}{\rmdefault}{\mddefault}{\updefault}{\color[rgb]{0,0,0}\(\partial \DDD(\phi)\)}%
}}}
\end{picture}

%% file: figs/LocalMult.pstex_t
\begin{picture}(0,0)%
\includegraphics{figs/LocalMult.pstex}%
\end{picture}%
\setlength{\unitlength}{1973sp}%
\begingroup\makeatletter\ifx\SetFigFont\undefined%
\gdef\SetFigFont#1#2#3#4#5{%
  \reset@font\fontsize{#1}{#2pt}%
  \fontfamily{#3}\fontseries{#4}\fontshape{#5}%
  \selectfont}%
\fi\endgroup%
\begin{picture}(13041,9148)(97,-10055)
\put(10801,-9961){\makebox(0,0)[lb]{\smash{\SetFigFont{12}{14.4}{\rmdefault}{\mddefault}{\updefault}{\color[rgb]{0,0,0}e) \(n=0\)}%
}}}
\put(6001,-9886){\makebox(0,0)[lb]{\smash{\SetFigFont{12}{14.4}{\rmdefault}{\mddefault}{\updefault}{\color[rgb]{0,0,0}d) \(n=0\)}%
}}}
\put(1126,-9886){\makebox(0,0)[lb]{\smash{\SetFigFont{12}{14.4}{\rmdefault}{\mddefault}{\updefault}{\color[rgb]{0,0,0}c) \(n=0\)}%
}}}
\put(3001,-5011){\makebox(0,0)[lb]{\smash{\SetFigFont{12}{14.4}{\rmdefault}{\mddefault}{\updefault}{\color[rgb]{0,0,0}a) \(n=-1/2\)}%
}}}
\put(8176,-5011){\makebox(0,0)[lb]{\smash{\SetFigFont{12}{14.4}{\rmdefault}{\mddefault}{\updefault}{\color[rgb]{0,0,0}b) \(n=-1/2\)}%
}}}
\end{picture}

%% file: figs/braid.pstex_t
\begin{picture}(0,0)%
\includegraphics{figs/braid.pstex}%
\end{picture}%
\setlength{\unitlength}{1184sp}%
\begingroup\makeatletter\ifx\SetFigFont\undefined%
\gdef\SetFigFont#1#2#3#4#5{%
  \reset@font\fontsize{#1}{#2pt}%
  \fontfamily{#3}\fontseries{#4}\fontshape{#5}%
  \selectfont}%
\fi\endgroup%
\begin{picture}(14604,6477)(1651,-9415)
\put(1651,-6211){\makebox(0,0)[lb]{\smash{\SetFigFont{12}{14.4}{\rmdefault}{\mddefault}{\updefault}{\color[rgb]{0,0,0}\(x_1\)}%
}}}
\put(9076,-6136){\makebox(0,0)[lb]{\smash{\SetFigFont{12}{14.4}{\rmdefault}{\mddefault}{\updefault}{\color[rgb]{0,0,0}\(y_1\)}%
}}}
\put(4501,-5386){\makebox(0,0)[lb]{\smash{\SetFigFont{12}{14.4}{\rmdefault}{\mddefault}{\updefault}{\color[rgb]{0,0,0}\(x_2=y_2\)}%
}}}
\end{picture}

%% file: figs/LocalMult2.pstex_t
\begin{picture}(0,0)%
\includegraphics{figs/LocalMult2.pstex}%
\end{picture}%
\setlength{\unitlength}{1579sp}%
\begingroup\makeatletter\ifx\SetFigFont\undefined%
\gdef\SetFigFont#1#2#3#4#5{%
  \reset@font\fontsize{#1}{#2pt}%
  \fontfamily{#3}\fontseries{#4}\fontshape{#5}%
  \selectfont}%
\fi\endgroup%
\begin{picture}(13813,9465)(64,-9980)
\put(11326,-4711){\makebox(0,0)[lb]{\smash{\SetFigFont{10}{12.0}{\rmdefault}{\mddefault}{\updefault}{\color[rgb]{0,0,0}c)  \(\nu = 1 \) }%
}}}
\put(1051,-9886){\makebox(0,0)[lb]{\smash{\SetFigFont{10}{12.0}{\rmdefault}{\mddefault}{\updefault}{\color[rgb]{0,0,0}d)  \(\nu = 1 \) }%
}}}
\put(6226,-9886){\makebox(0,0)[lb]{\smash{\SetFigFont{10}{12.0}{\rmdefault}{\mddefault}{\updefault}{\color[rgb]{0,0,0}e)  \(\nu = -1 \) }%
}}}
\put(11401,-9886){\makebox(0,0)[lb]{\smash{\SetFigFont{10}{12.0}{\rmdefault}{\mddefault}{\updefault}{\color[rgb]{0,0,0}f)  \(\nu = 3 \) }%
}}}
\put(1051,-4711){\makebox(0,0)[lb]{\smash{\SetFigFont{10}{12.0}{\rmdefault}{\mddefault}{\updefault}{\color[rgb]{0,0,0}a)  \(\nu = 0 \) }%
}}}
\put(6001,-4711){\makebox(0,0)[lb]{\smash{\SetFigFont{10}{12.0}{\rmdefault}{\mddefault}{\updefault}{\color[rgb]{0,0,0}b)  \(\nu = 1 \) }%
}}}
\end{picture}

%% file: gridley.bbl
\begin{thebibliography}{10}

\bibitem{ExactT}
P.~Braam and S.~Donaldson.
\newblock Floer's work on instanton homology, knots, and surgery.
\newblock In {\em The Floer Memorial Volume}, pages 195--256. Birkhauser, 1995.

\bibitem{CrowFox}
R.H. Crowell and R.H. Fox.
\newblock {\em Introduction to Knot Theory}.
\newblock Ginn and Co., 1963.

\bibitem{FSSFHF}
R.~Fintushel and R.~Stern.
\newblock Instanton homology of {S}eifert-fibred homology spheres.
\newblock {\em Proc. London Math. Soc.}, 61:109--137, 1990.

\bibitem{Floer1}
A.~Floer.
\newblock Witten's complex and infinite dimensional {M}orse theory.
\newblock {\em J. Differential Geom.}, 30:207--21, 1989.

\bibitem{Froyshov}
K.A. Fr{\o}yshov.
\newblock Equivariant aspects of {Y}ang-{M}ills {F}loer theory.
\newblock {\em Topology}, 41:525--552, 2002.

\bibitem{FrLect}
K.A. Fr{\o}yshov.
\newblock Lectures at {H}arvard {U}niversity.
\newblock Spring 2001.

\bibitem{Kauffman}
L.~Kauffman.
\newblock {\em On Knots}.
\newblock Annals of Mathematics Studies, 115. Princeton University Press, 1987.

\bibitem{Khovanov}
M.~Khovanov.
\newblock A categorification of the {J}ones polynomial.
\newblock {\em Duke Math. J.}, 101:359--426, 2000.

\bibitem{Khovanov2}
M.~Khovanov.
\newblock Patterns in knot cohomology i.
\newblock math.QA/0201306, 2002.

\bibitem{KKR}
P.~Kirk, E.~Klassen, and D.~Ruberman.
\newblock Splitting the spectral flow and the {A}lexander matrix.
\newblock {\em Comment. Math. Helv.}, 69:375--416, 1994.

\bibitem{Kr}
P.~B. Kronheimer.
\newblock Embedded surfaces and gauge theory in three and four dimensions.
\newblock In {\em Surveys in Differential Geometry, Vol. III}, pages 243--298.
  Int. Press, 1998.

\bibitem{MKLect}
P.~B. Kronheimer and T.~S. Mrowka.
\newblock 1996 {I.P.} lecture.

\bibitem{ESL2}
E.~S. Lee.
\newblock {K}hovanov's invariants for alternating links.
\newblock math.GT/0210213, 2002.

\bibitem{Man}
C.~Manolescu.
\newblock Seiberg-{W}itten-{F}loer stable homotopy type of three-manifolds with
  \(b_1 = 0 \).
\newblock {\\} math.DG/0104024, 2001.

\bibitem{UsersGuide}
J.~McCleary.
\newblock {\em User's Guide to Spectral Sequences}.
\newblock Mathematics Lecture Series, 12. Publish or Perish Inc, 1985.

\bibitem{McMullen}
C.~McMullen.
\newblock The {A}lexander polynomial of a 3-manifold and the {T}hurston norm on
  cohomology.
\newblock {\em Ann. Sci. \'Ecole Norm. Sup. (4)}, 35:153--171, 2002.

\bibitem{MOY}
T.~Mrowka, P.~Ozsv{\'a}th, and B.~Yu.
\newblock Seiberg-{W}itten monopoles on {S}eifert-fibred spaces.
\newblock {\em Comm. Anal. Geom.}, 4:685--791, 1997.

\bibitem{OS2}
P.~Ozsv{\'a}th and Z.~Szab{\'o}.
\newblock Holomorphic disks and three-manifold invariants: properties and
  applications.
\newblock To appear in Annals of Mathematics. math.SG/0105202.

\bibitem{OS1}
P.~Ozsv{\'a}th and Z.~Szab{\'o}.
\newblock Holomorphic disks and topological invariants for rational homology
  three-spheres.
\newblock To appear in Annals of Mathematics. math.SG/0101206.

\bibitem{OS4}
P.~Ozsv{\'a}th and Z.~Szab{\'o}.
\newblock Absolutely graded {F}loer homologies and intersection forms for
  four-manifolds with boundary.
\newblock math.SG/0110170, 2001.

\bibitem{OS3}
P.~Ozsv{\'a}th and Z.~Szab{\'o}.
\newblock Holomorphic triangles and invariants for smooth four-manifolds.
\newblock {\\} math.SG/0110169, 2001.

\bibitem{OS8}
P.~Ozsv{\'a}th and Z.~Szab{\'o}.
\newblock Heegaard {F}loer homology and alternating knots.
\newblock math.GT/0209149, 2002.

\bibitem{OS7}
P.~Ozsv{\'a}th and Z.~Szab{\'o}.
\newblock Holomorphic disks and knot invariants.
\newblock math.GT/0209056, 2002.

\bibitem{OS10}
P.~Ozsv{\'a}th and Z.~Szab{\'o}.
\newblock Knot {F}loer homology and the four-ball genus.
\newblock math.GT/0301026, 2003.

\bibitem{2bridge}
J.~Rasmussen.
\newblock Floer homology of surgeries on two-bridge knots.
\newblock {\em Algebraic and Geometric Topology}, 2:757--89, 2002.
\newblock math.GT/020405.

\bibitem{Rolfsen}
D.~Rolfsen.
\newblock {\em Knots and Links}.
\newblock Mathematics Lecture Series, 7. Publish or Perish Inc., 1976.

\bibitem{Schu}
H.~Schubert.
\newblock Knoten mit zwei {B}r{\"{u}}cken.
\newblock {\em Math. Z.}, 65:133--70, 1956.

\bibitem{Turaev}
V.~Turaev.
\newblock Torsion invariants of \( \spinc \) structures on 3-manifolds.
\newblock {\em Math. Res. Lett.}, 4:679--695, 1997.

\bibitem{Viterbo}
C.~Viterbo.
\newblock Intersection de sous-vari{\'e}t{\'e}s lagrangiennes, fonctionelles
  d'action, et indice des syst{\`e}mes hamiltoniens.
\newblock {\em Bull. Soc. Math. France}, 115:361--390, 1987.

\end{thebibliography}
